\documentclass[11pt,a4paper]{article}
\usepackage{amsmath,amsfonts,amssymb,amsthm}

\newcommand{\ac}[2]{a_{#1 #2}}

\newcommand{\al}{\alpha}
\newcommand{\arcsinh}{\operatorname{arcsinh}}
\newcommand{\Aut}{\operatorname{Aut}}
\newcommand{\ba}{\begin{array}}
\newcommand{\bc}[2]{b_{#1 #2}}
\newcommand{\be}{\mathbf{e}}
\newcommand{\bet}{\beta}
\newcommand{\bof}{\mathbf{f}}
\newcommand{\bt}{\mathbf{t}}
\newcommand{\bv}{\mathbf{v}}
\newcommand{\bw}{\mathbf{w}}
\newcommand{\bx}{\mathbf{x}}
\newcommand{\by}{\mathbf{y}}
\newcommand{\C}{\mathbb{C}}
\newcommand{\cc}{\bar{c}}
 
\newcommand{\del}{\partial}
\newcommand{\delh}{\partial (h)}
 
\newcommand{\ea}{\end{array}}
\newcommand{\efct}{\mathrm{e}}
\newcommand{\Eig}{\mathrm{Eig}}
\newcommand{\eps}{\varepsilon}
\newcommand{\folgt}{\Rightarrow}
\newcommand{\ga}{\gamma}
\newcommand{\Ga}{\Gamma}
 \newcommand{\grad}{\operatorname{grad}} 
\newcommand{\half}{{\tfrac{1}{2}}}
\newcommand{\hs}{\varphi}
\newcommand{\I}{\mathrm{i}}
\newcommand{\Id}{I\!d}
\newcommand{\la}{\lambda}
\newcommand{\lap}{\Delta}
\newcommand{\lp}{\left<}
\newcommand{\M}{M^3}
\newcommand{\map}{L}
\newcommand{\mat}[9]{\begin{pmatrix} #1& #2& #3\\#4& #5& #6\\#7& #8& #9\end{pmatrix}}
\newcommand{\Min}{\R^3_1}
\newcommand{\p}{\partial}
\newcommand{\pa}{\frac{\partial}{\partial a}}
\newcommand{\paa}{\frac{\partial}{\partial \ac22}}
\newcommand{\pbb}{\frac{\partial}{\partial \ac23}}
\newcommand{\pt}{\frac{\p}{\p t}}
\newcommand{\ptv}{\frac{\p}{\p v}}
\newcommand{\ptw}{\frac{\p}{\p w}}

\newcommand{\pv}{\frac{\partial}{\partial v}}
\newcommand{\pw}{\frac{\partial}{\partial w}}
\newcommand{\pz}{\frac{\partial}{\partial z}}
\newcommand{\pp}[1]{\frac{\partial^2}{\partial {#1}^2}}
\newcommand{\R}{\mathbb{R}} 
\newcommand{\Rf}{\R^4} 
\newcommand{\ricLC}{\widehat{\operatorname{Ric}}}
\newcommand{\rp}{\right>}
\newcommand{\sign}{\operatorname{sign}}

\newcommand{\SO}{\mathrm{SO}}
\newcommand{\Span}{\mathrm{span}} 

\newcommand{\third}{{\tfrac{1}{3}}}
\newcommand{\trace}{\operatorname{trace}} 

\newcommand{\tT}{\tilde{T}}
\newcommand{\tV}{\tilde{V}}
\newcommand{\tW}{\tilde{W}}
\newcommand{\twothird}{{\tfrac{2}{3}}}

\newcommand{\wN}{\widehat{\nabla}}
\newcommand{\VW}{(V+\I W)}
\newcommand{\Z}{\mathbb{Z}}

\newtheorem{thm}{Theorem}
\newtheorem{lem}{Lemma}

\theoremstyle{remark}
\newtheorem*{rem}{Remark}

\newcounter{rom}
\renewcommand{\therom}{(\roman{rom})}
\newenvironment{romanlist}{\begin{list}{\therom}
                {\setlength{\leftmargin}{2em}\usecounter{rom}}}%
{\end{list}}
\begin{document}
\title{Indefinite affine hyperspheres admitting a pointwise symmetry}
\author{C. Scharlach\thanks{Partially supported by the DFG-Project
'Geometric Problems and Special PDEs'}} \maketitle

\noindent

\bibliographystyle{ADG-Berlin}

\begin{abstract}
An affine hypersurface $M$ is said to admit a pointwise symmetry, if
there exists a subgroup $G$ of $\Aut(T_p M)$ for all $p\in M$, which
preserves (pointwise) the affine metric $h$, the difference tensor $K$
and the affine shape operator $S$. Here, we consider 
3-dimensional indefinite affine hyperspheres, i.~e.\ $S= H\Id$ (and
thus $S$ is trivially preserved). First we solve an algebraic
problem. We determine the non-trivial stabilizers $G$ of a traceless
cubic form on a Lorentz-Minkowski space $\Min$ under the action of
the isometry group $SO(1,2)$ and find a representative of each
$SO(1,2)/G$-orbit. Since the affine cubic form is defined by $h$ and
$K$, this gives us the possible symmetry groups $G$ and for each $G$ a
canonical form of $K$. Next, we classify hyperspheres admitting a
pointwise $G$-symmetry for all non-trivial stabilizers $G$ (apart from 
$\Z_2$). Besides well-known hyperspheres (for $\Z_2\times \Z_2$
resp.\ $\R$ the hyperspheres have constant sectional curvature and Pick
invariant $J<0$ resp.\ $J=0$) we obtain rich classes of new examples
e.g. warped product structures of two-dimensional affine spheres
(resp.\ quadrics) and curves. Moreover, we find a way to construct
indefinite affine hyperspheres out of 2-dimensional quadrics or
positive definite affine spheres.
\end{abstract}

\medskip\noindent 
{\bfseries Subject class: } 53A15 (15A21, 53B30)
 

\section{Introduction}\label{sec:intro}

Let $M^n$ be a connected, oriented manifold. Consider an immersed
        hypersurface with relative normalization, i.e., an immersion
        $\hs\colon M^n \rightarrow \R^{n+1}$ together with a
        transverse vector field $\xi$ such that $D \xi$ has its image
        in $\hs_*T_xM$. Equi-affine geometry studies the properties of 
        such immersions under equi-affine transformations, i.~e. 
        volume-preserving linear transformations ($SL(n+1,\R)$) and 
        translations.
        
	In the theory of nondegenerate equi-affine hypersurfaces there
	exists a canonical choice of transverse vector field $\xi$
	(unique up to sign), called the affine (Blaschke) normal,
	which induces a connection $\nabla$, a nondegenerate symmetric
	bilinear form $h$ and a 1-1 tensor field $S$ by
	\begin{align} 
	  &D_X Y =\nabla_X Y +h(X,Y)\xi,\label{strGauss}\\
	  &D_X \xi =-SX,\label{strWeingarten}
	\end{align}
	for all $X,Y \in {\cal X}(M)$. The connection $\nabla$ is
	called the induced affine connection, $h$ is called the affine
	metric or Blaschke metric and $S$ is called the affine shape
	operator.  In general $\nabla$ is not the Levi Civita
	connection $\hat\nabla$ of $h$. The difference tensor $K$ is
	defined as
	\begin{equation}\label{defK}
	  K(X,Y)=\nabla_X Y-\hat\nabla_X Y,
	\end{equation}
	for all $X,Y \in {\cal X}(M)$. Moreover the form $h(K(X,Y),Z)$
	is a symmetric cubic form with the property that for any fixed
	$X\in {\cal X}(M)$, $\trace K_X$ vanishes.  This last property
	is called the apolarity condition. The difference tensor $K$,
	together with the affine metric $h$ and the affine shape
	operator $S$ are the most fundamental algebraic invariants for
	a nondegenerate affine hypersurface (more details in
	Sec.~\ref{sec:basics}). We say that $M^n$ is indefinite,
	definite, etc. if the affine metric $h$ is indefinite,
	definite, etc. For details of the basic theory of
	nondegenerate affine hypersurfaces we refer to \cite{LSZ93}
	and \cite{NS94}.
	
	Here we will restrict ourselves to the case of affine
	hyperspheres, i.~e.\ the shape operator will be a (constant)
	multiple of the identity ($S= H \Id$). Geometrically this
	means that alle affine normals pass through a fixed point or
	they are parallel. The abundance of affine hyperspheres dwarfs
	any attempts at a complete classification. Even with the
	restriction to locally strongly convex hyperspheres (i.~e.\ $h$
	is positive definite) and low dimensions the class is simply
	too large to classify. In order to obtain detailed information
	one has therefore to revert to sub-classes such as the class
	of complete affine hyperspheres (see \cite{LSZ93} and the
	references contained therein or i.~e.\ \cite{JL05} for a very
	recent result). Various authors have also imposed curvature
	conditions. In the case of constant curvature the
	classification is nearly finished (see \cite{Vr00} and the
	references contained therein). In analogy to Chen's work,
	\cite{Che93}, a new curvature invariant for positive definite
	affine hyperspheres was introduced in \cite{SSVV97}. A lower
	bound was given and, for $n=3$, the classification of the
	extremal class was started. This classification was completed
	in \cite{SV96}, \cite{KSV01} and \cite{KV99}. The special
	(simple) form of the difference tensor $K$ for this class is
	remarkable, actually it turns out that the hyperspheres admit
	a certain pointwise group symmetry \cite{Vr04}.

      A hypersurface is said to admit a pointwise group symmetry if at
      every point the affine metric, the affine shape operator and the
      difference tensor are preserved under the group
      action. Necessarily the possible groups must be subgroups of the
      isometry group. The study of submanifolds which admit a
      pointwise group symmetry was initiated by Bryant in \cite{Bry01}
      where he studied $3$-dimensional Lagrangian submanifolds of
      $\mathbb C^3$, i.~e.\ the isometry group is $\SO(3)$. Because of
      the similar basic invariants, Vrancken transferred the problem
      to $3$-dimensional positive definite affine hyperspheres. A
      classification of $3$-dimensional positive definite affine
      hyperspheres admitting pointwise symmetries was obtained in
      \cite{Vr04} and then extended to positive definite hypersurfaces
      in \cite{LS04} (here the affine shape operator is non-trivial
      and thus no longer trivially preserved by isometries). Now, for
      the first time, we will consider the indefinite case, namely
      $3$-dimensional indefinite affine hyperspheres.
      
      We can assume that the affine metric has index two, i.~e.\ the
      corresponding isometry group is the (special) Lorentz group
      $\SO(1,2)$. Our question is the following: {\em What can we say
      about a three-dimensional indefinite affine hypersphere, for
      which there exists a non-trivial subgroup $G$ of $SO(1,2)$ such
      that for every $p\in M$ and every $L\in G$:
$$K(LX_p,L Y_p)= L(K(X_p, Y_p))\quad \forall X_p, Y_p\in T_p M.$$} In
      Section~\ref{sec:basics} we will state the basic formulas of
      (equi-)affine hypersurface-theory needed in the further
      classification. We won't need hypersurface-theory in Section
      \ref{sec:nf} and \ref{sec:max}, were we consider the group
      structure of $\SO(1,2)$ and its action on cubic forms. In
      Section \ref{sec:nf} we show that there exist six different
      normalforms of elements of $\SO(1,2)$, depending on the
      eigenvalues and eigenspaces. We can always find an oriented
      basis of $\Min$ such that every $L\in \SO(1,2)$ has one (and
      only one) of the following matrix representations:
      $$\left(\begin{smallmatrix} 1 & 0 & 0 \\ 0 & \cos t & -\sin t\\
      0 & \sin t & \cos t \end{smallmatrix}\right), t\in (0,
      2\pi)\setminus \{\pi\},\quad \left(\begin{smallmatrix} 1 & 0 & 0
      \\ 0 & -1 & 0\\ 0 & 0 & -1 \end{smallmatrix}\right),\quad \Id,
      \quad \left(\begin{smallmatrix} -1 & 0 & 0 \\ 0 & 1 & 0\\ 0 & 0
      & -1
      \end{smallmatrix}\right),$$ all of these with respect to an
      ONB $\{ \bt, \bv, \bw\}$, $\bt$ timelike, $\bv,
      \bw$ spacelike, or 
     $$\left(\begin{smallmatrix} l & 0 & 0 \\ 0 & 1
      & 0\\ 0 & 0 & \frac{1}{l}
    \end{smallmatrix}\right),\; l\neq\pm1, \quad \left(\begin{smallmatrix} 
    1 & -1 & -\frac12 \\ 0 & 1 & 1\\ 0 & 0 & 1
    \end{smallmatrix}\right),$$ with respect to a (LV)basis $\{
    \be,\bv,\bof\}$, $\be, \bof$ lightlike, $\bv$ spacelike
    (Theorem~\ref{thm:NF}).

Since we are interested in pointwise group symmetry, in
    Sec.~\ref{sec:max} we study the nontrivial stabilizer of a
    traceless cubic form $\tilde{K}$ under the $\SO(1,2)$-action $\rho
    (L)(\tilde{K}) = \tilde{K}\circ L$.
(cp. \cite{Bry01} for the classification of the $\SO (3)$-action). It
turns out that the $\SO(1,2)$-stabilizer of a nontrivial traceless
cubic form is isomorphic to either $\SO (2)$, $\SO(1,1)$, $\R$, the
group $S_3$ of order 6, $\Z_2 \times \Z_2$, $\Z_3$, $\Z_2$ or it is
trivial (Theorem~\ref{thm:Subgroups}). 

In the following we classify the indefinite affine hyperspheres which
admit a pointwise $\Z_2 \times \Z_2$-symmetry
(Section~\ref{sec:type6}), $\R$-symmetry (Section~\ref{sec:type9}),
$\SO(2)$- or $\Z_3$-symmetry (Section~\ref{subsec:type2,4}),
$S_3$-symmetry (Section~\ref{subsec:type3}) or $\SO(1,1)$-symmetry
(Section~\ref{sec:type8}). In case of $\Z_2 \times \Z_2$ resp.\ $\R$,
we get the indefinite affine hyperspheres of constant sectional
curvature with negative resp.\ vanishing Pick invariant $J$; those
with $J>0$ are examples for $\Z_2$-symmetry. The other classes are
very rich, most of them are warped products of two-dimensional affine
spheres ($\Z_3$) resp. quadrics ($\SO(2)$, $\SO(1,1)$), with a
curve. Thus we get many new examples of 3-dimensional indefinite
affine hyperspheres. Furthermore, we show how one can construct
indefinite affine hyperspheres out of two-dimensional quadrics or
positive definite affine spheres.

\section{Basics of affine hypersphere theory}
\label{sec:basics}

First we recall the definition of the affine normal $\xi$
(cp. \cite{NS94}). In equi-affine hypersurface theory on the ambient
space $\R^{n+1}$ a fixed volume form $\det$ is given. A transverse
vector field $\xi$ induces a volume form $\theta$ on $M$ by
$\theta(X_1,\ldots,X_n)=\det(\hs_* X_1,\ldots,\hs_* X_n,\xi)$. Also
the affine metric $h$ defines a volume form $\omega_h$ on $M$, namely
$\omega_h=|\det h|^{1/2}$. Now the affine normal $\xi$ is uniquely
determined (up to sign) by the conditions that $D \xi$ is everywhere
tangential (which is equivalent to $\nabla \theta =0$) and that
\begin{equation}\label{BlaschkeNormal}
\theta = \omega_h.
\end{equation}
Since we only consider 3-dimensional indefinite hyperspheres, i.~e. 
\begin{equation}\label{affHypSphere}
S= H \Id,\quad H=\text{const.}
\end{equation}
we can fix the orientation of the affine normal $\xi$ such that the
affine metric has signature one. Then the sign of $H$ in the
definition of an affine hypersphere is an invariant.

Next we state some of the fundamental equations, which a
nondegenerate hypersurface has to satisfy, see also \cite{NS94} or
\cite{LSZ93}. These equations relate $S$ and $K$ with amongst others
the curvature tensor $R$ of the induced connection $\nabla$ and the
curvature tensor $\hat R$ of the Levi Civita connection
$\widehat\nabla$ of the affine metric $h$.  There are the
Gauss equation for $\nabla$, which states that:
\begin{equation*}
R(X,Y)Z =h(Y,Z)SX -h(X,Z) SY,
\end{equation*}
and the Codazzi equation 
\begin{equation*}
(\nabla_X S) Y =(\nabla_Y S) X.
\end{equation*}
Also we have the total symmetry of the affine cubic form
\begin{equation}\label{defC}
	C(X,Y,Z)= (\nabla_X h) (Y,Z) = -2 h(K(X,Y),Z).
\end{equation} 
The fundamental existence and uniqueness theorem, see \cite{Dil89} or
\cite{DNV90}, states that given $h$, $\nabla$ and $S$ such that the
difference tensor is symmetric and traceless with respect to $h$, on a
simply connected manifold $M$ an affine immersion of $M$ exists if and
only if the above Gauss equation and Codazzi equation are satisfied.

	From the Gauss equation and Codazzi equation above the Codazzi
equation for $K$ and the Gauss equation for $\widehat\nabla$ follow:
\begin{equation*}\begin{split}
(\widehat{\nabla}_X K)(Y,Z)- (\widehat{\nabla}_Y K)(X,Z)=& \half( h(Y,Z)SX-
h(X,Z)SY\\& - h(SY,Z)X + h(SX, Z)Y),\end{split}
\end{equation*}
and
\begin{equation*}\begin{split} \hat R (X,Y)Z =&
\tfrac{1}{2} (h(Y,Z) SX - h(X,Z) SY\\& + h(SY,Z)X - h(SX,Z)Y ) - [K_X
    ,K_Y] Z
\end{split}\end{equation*} 
If we define the Ricci tensor of the Levi-Civita connection $\widehat
\nabla$ by:
\begin{equation}\label{def:RicLC}
\ricLC(X,Y)=\trace\{ Z \mapsto \hat R(Z,X)Y\}.
\end{equation}
and the Pick invariant by:
\begin{equation}\label{def:Pick}
J= \frac{1}{n (n-1)} h(K,K),
\end{equation}
then from the Gauss equation we immediately get for the scalar
curvature $\hat{\kappa}=\frac{1}{n (n-1)}(\sum_{i,j} h^{ij}
\ricLC_{ij})$:
\begin{equation}\label{TE}
\hat \kappa= H+J.
\end{equation}
For an affine hypersphere the Gauss and Codazzi equations have the form:
\begin{align}\label{gaussInd}
&R(X,Y)Z=H(h(Y,Z)X -h(X,Z) Y),\\
\label{codazziS} &(\nabla_X H)Y=(\nabla_Y H)X, \quad \text{i. e.}
\quad H= const.,\\
\label{CodK} &(\widehat{\nabla}_X K)(Y,Z)=(\widehat{\nabla}_Y K)(X,Z),\\
\label{gaussLC} &\hat R (X,Y)Z = H(h(Y,Z)X -h(X,Z) Y)- [K_X ,K_Y] Z.
\end{align}
Since $H$ is constant, we can rescale $\hs$ such that $H\in \{-1,0,1\}$.

\section{Normalforms in $\SO(1,2)$}\label{sec:nf}

We denote\footnote{for the notation cp. \cite{Gr67}} by $\Min$ the
Pseudo-Euclidean vector space in which a non-degenerate indefinite
bilinear form of index two is given. The bilinear form is called
the inner product and denoted by $\lp\;,\;\rp$.  A basis $\{
\bt,\bv,\bw\}$ is called orthonormal (ONB) if
\begin{equation}\label{ONB}\begin{split} \lp\bt,\bt\rp &=-1, 
\quad \lp\bv,\bv\rp =1=\lp\bw,\bw\rp, \\ 
0&= \lp\bt,\bv\rp=\lp\bt,\bw\rp=\lp\bv,\bw\rp. \end{split}
\end{equation}
For a chosen ONB the inner product of two vectors is given by 
\begin{equation}\label{IP:ONB}\lp\bx,\by\rp=-x_t y_t + x_v y_v + x_w y_w, 
\quad \bx,\by \in \Min.\end{equation}
A basis $\{ \be,\bv,\bof\}$ is called a light-vector basis (LVB) if  
\begin{equation}\label{LVB} \begin{split} \lp\be,\be\rp=0=\lp\bof,\bof\rp, 
& \lp\be,\bof\rp=1,\\ \lp\be,\bv\rp=0=\lp\bof,\bv\rp, &
\lp\bv,\bv\rp=1. \end{split}\end{equation} 
For a chosen LVB the inner product of two vectors is given by
\begin{equation}\label{IP:LVB}\lp\bx,\by\rp= x_e y_f + x_f y_e + x_v y_v, 
\quad \bx,\by \in \Min.\end{equation}
We want
to consider the special Pseudo-Euclidean rotations $\SO(1,2)$,
i.~e. the linear transformations $\map$ of $\Min$ which preserve
the inner product and have determinant equal to one:
$$ \lp\map \bx, \map \by\rp= \lp \bx,\by\rp, \quad \bx,\by \in \Min,
\quad \det L =1.$$ Depending on the eigenvalues and eigenspaces we get
the following normalforms of the elements of $SO(1,2)$.

\begin{thm} \label{thm:NF} There exists a choice of an oriented basis
of $\R^3_1$ such that every $L \in \SO(1,2)$ is
of one (and only one) of the following types: 
\begin{enumerate} \item
  \begin{description} \item [\it(a)] 
    $A_t:= \left(\begin{smallmatrix} 1 & 0 & 0 \\ 0 & \cos t & -\sin
      t\\ 0 & \sin t & \cos t \end{smallmatrix}\right)$, $t\in [0,
      2\pi)$, $t\neq 0,\pi$, \\ for an ONB $\{ \bt, \bv, \bw\}$, $\bt$
      timelike, $\bv, \bw$ spacelike,\\ eigenvalues: $\lambda_1=1$,
      eigenspaces $E(1)=\Span\{t\}$ timelike.
  \item[\it(b)] $A_\pi:= \left(\begin{smallmatrix} 1 & 0 & 0 \\ 0 & -1 &
    0\\ 0 & 0 & -1 \end{smallmatrix}\right)$, for an ONB $\{ \bt,
    \bv, \bw\}$, as above,\\ eigenvalues: $\lambda_1=1$,
    $\lambda_{2,3}=-1$, \\ eigenspaces $E(1)=\Span\{\bt\}$ timelike,
    $E(-1)=\Span\{\bv,\bw\}$, spacelike.
  \end{description}
\item 
	\begin{description} \item [\it(a)] 
	$A_0:= \Id$,
    for an ONB $\{ \bt, \bv, \bw\}$ as above or a LVB $\{ \be,\bv,\bof\}$,\\
    eigenvalues: $\lambda_{1,2,3}=1$, eigenspaces $E(1)=\Min$.
  \item[\it(b)] $B:=\left(\begin{smallmatrix} -1 & 0 & 0 \\ 0 & 1 & 0\\ 0
	& 0 & -1 \end{smallmatrix}\right)$, for an ONB $\{ \bt,
	\bv, \bw\}$, as above, or a LVB $\{ \be,\bv,\bof\}$,\\
	eigenvalues: $\lambda_{1,3}=-1$, $\lambda_2=1$, \\
	eigenspaces $E(-1)=\Span\{\bt,\bw\}=\Span\{\be,\bof\}$
	timelike, $E(1)=\Span\{\bv\}$, spacelike.
\end{description}
\item 
  \begin{description} \item [\it(a)] $C_l:=\left(\begin{smallmatrix} 
    l & 0 & 0 \\ 0 & 1 & 0\\ 0 & 0 & \frac{1}{l}
    \end{smallmatrix}\right)$, $l\neq\pm1$, for a LVB $\{ \be,\bv,\bof\}$,\\ 
    eigenvalues: $\la_1=l$, $\la_2=1$, $\la_3=\frac{1}{l}$, \\
    eigenspaces $E(l)=\Span\{\be\}$ lightlike,
    $E(1)=\Span\{\bv\}$, spacelike, $E(\frac{1}{l})=\Span\{\bof\}$
    lightlike.
    \item [\it(b)] $C_1:=\left(\begin{smallmatrix} 
    1 & -1 & -\frac12 \\ 0 & 1 & 1\\ 0 & 0 & 1 \end{smallmatrix}\right)$,
    for a LVB $\{ \be,\bv,\bof\}$,\\
    eigenvalues: $\la_{1,2,3}=1$, 
    eigenspaces $E(1)=\Span\{\be\}$ lightlike.
  \end{description}
\end{enumerate}
\end{thm} 
\noindent For the proof we will use the following

\begin{lem}\label{lem:LVBtoLVB}
Let $L \in \SO(1,2)$ with $L(\be)= l \be$, $l\neq 0$, for
a LVB $\{ \be,\bv,\bof\}$. Since under $L$ a LVB will be mapped
to a LVB, the corresponding matrix must have the following form :
$$ C_{l,m}=\mat{l}{-l m}{-l\frac{m^2}{2}}{0}{1}{m}
{0}{0}{\frac{1}{l}}. $$ 
\end{lem}
\begin{proof} We will use the 
Notation: $\be'=\map(\be)$, $\bv'=\map(\bv)$, $\bof'=\map(\bof)$, and
compute \eqref{LVB}, using \eqref{IP:LVB}:
\begin{align*} 0=\lp\be',\bv'\rp=l v'_f &\folgt v'_f=0,\\
1=\lp\bv',\bv'\rp=(v'_v)^2 &\folgt v'_v=\eps,\quad \eps^2=1,\\ 
1=\lp\be',\bof'\rp=l f'_f &\folgt
f'_f=\frac{1}{l},\\ 1=det \map = l \eps \frac{1}{l} &\folgt
\eps=1,\\ 0=\lp\bof',\bv'\rp=v'_e \frac{1}{l}+f'_v &\folgt f'_v=
-\frac{v'_e}{l},\\ 0=\lp\bof',\bof'\rp= 2 f'_e \frac{1}{l} +
\frac{(v'_e)^2}{2} &\folgt f'_e=-\frac{1}{l}\frac{(v'_e)^2}{2}.
\end{align*}
With $v'_e=-l m$ we obtain $C_{l,m}$.
\end{proof}

\begin{proof}[Proof of Thm.~\ref{thm:NF}]
Every $L \in \SO(1,2)$ must have a real eigenvalue since it is an 
automorphism of a three dimensional vector space. A
corresponding eigenvector will be either timelike, spacelike or
lightlike.

Let's consider first the case that we have (at least) a timelike
eigenvector.  We can choose an ONB $\{ \bt, \bv, \bw\}$, such that $\bt$ is
the timelike eigenvector with eigenvalue $\eps = \pm 1$. Then $\bv,
\bw$ are spacelike and $L$ restricted to $\bt^{\bot}=\Span\{\bv,\bw\}$
is an isometry of $\R^2$, i.~e. an Euclidean rotation.  

If $\eps = 1$, then $L$ restricted to $\bt^{\bot}$ is in
$\SO(2)$ (proper Euclidean rotation): In general we get no more real
eigenvalues (case {\it 1.}(a)). If the restriction is a rotation by an
angle of $\pi$, then we get the second eigenvalue $-1$ of
multiplicity two. The eigenspace is $\bt^{\bot}$ (and thus
spacelike) (case {\it 1.}(b)).  Finally, if the restriction (and thus
$L$) is the identity map $\Id$, every vector is an eigenvector
and we also can choose a LVB (case {\it 2.}(a)).

If $\eps = -1$, then
$L$ restricted to $\bt^{\bot}$ is an improper Euclidean
rotation of $\R^2$, thus it has eigenvalues $1$ and $-1$. We get the
eigenvalue $1$ with eigenspace $\Span\{\bv\}$ (spacelike) and the
eigenvalue $-1$ of multiplicity two with eigenspace $\Span\{\bt,\bw\}$
(timelike). Since $L$ restricted to $\Span\{\bt,\bw\}$ is equal to
$-\Id$, we also can choose a basis of two lightlike eigenvectors
(case {\it 2.(b)}).

Next we will consider the case that we have (at least) a lightlike
eigenvector. We can choose a LVB $\{\be,\bv,\bof\}$, such that $\be$ is this
lightlike eigenvector with eigenvalue $l \in \R$, $l\neq 0$. Since under
$L$ a LVB will be mapped to a LVB, the corresponding matrix must
have the following form (cp. Lem.~\ref{lem:LVBtoLVB}):
$$ C_{l,m}=\mat{l}{-l m}{-l\frac{m^2}{2}}{0}{1}{m}
{0}{0}{\frac{1}{l}}. $$ 

$C_{l,m}$ has the eigenvalues $\la_1=l$, $\la_2=1$ and
$\la_3=\frac{1}{l}$. If $l\neq\pm 1$,
we get three distinct real eigenvalues and the corresponding (1-dim.)
eigenspaces are either lightlike ($\Eig(l)$ and
$\Eig(\frac{1}{l})$) or spacelike ($\Eig(1)$) (case {\it
3.}(a)). (Since we take an eigenvector basis, the LVB is up to the
length of the lightlike eigenvectors uniquely determined.) If
$l=-1$, $\la_1=\la_3=-1$, thus we have an eigenvalue of
multiplicity two with eigenspace $\Span\{\be,\bof\}$ (case {\it
2.}(b)). The case $l=1$ is left, i.~e. only one eigenvalue of
multiplicity three: For $m=0$ we obtain the identity map (case {\it
2.}(a)). For $m\neq 0$, the eigenspace $\Eig(1)=\Span\{\be\}$ only is
one-dimensional. To get a normalform of $\map$ we compute how
$C_{1,m}$ changes if we choose another LVB $\{\be',\bv',\bof'\}$ (with the same orientation) where
$\be'$ is an eigenvector of $\map$ ($\be'=a \be$). If we express
$\{\be',\bv',\bof'\}$ in terms of $\{\be,\bv,\bof\}$ and compute
\eqref{LVB}, we obtain (cp. Lem.~\ref{lem:LVBtoLVB}):
\begin{align*} 
\be'&=a \be,\\
\bv'&= v'_e \be +\bv,\\
\bof'&= -\frac{(v'_e)^2}{2 a} \be - \frac{v'_e}{a}\bv+\frac{1}{a} \bof.
\end{align*}
The matrix representing $\map$ with respect to the new LVB,
$C'_{1,m}$, has the form $C'_{1,m}=\mat{1}{-\frac{m}{a}} {-\frac{m^2}{2}\frac{1}{a^2}}{0}{1}{\frac{m}{a}}{0}{0}{1}$. Thus $C'_{1,m}=C_{1,\frac{m}{a}}$, and we can
choose a LVB such that $\frac{m}{a}=1$ (case {\it 3.}(b)). From
the above computations we see that this LVB still isn't completly
determined, we can choose $v'_e$ arbitrary.

Finally the case is left that we have (at least) a spacelike eigenvector $\bv$. The 
corresponding eigenvalue must be $\eps =\pm 1$ and $L$ restricted 
to $\bv^{\bot}$ is an isometry of $\R^2_1$, i.~e. a Pseudo-Euclidean rotation 
(boost). Thus it always has two real eigenvalues with one-dimensional 
eigenspaces. We can choose an eigenvector basis, which will be either an 
ONB or a LVB, and we get one of the following cases: {\it 1.}(b), {\it 2.}(a),
{\it 2.}(b) or  {\it 3.}(a) (cp. \cite{Gr67}, p.~273).

\end{proof}

\begin{rem} \label{rem:NF} The choice of basis for the above
  normalforms is unique up to: 
\begin{enumerate} \item
  \begin{description} \item [\it(a)]
     the ONB is unique up
      to $\bt \to \eps \bt, \bv, \bw$ up to a proper
      ($\eps =1$) or improper ($\eps=-1$) Euclidean rotation in
      $\R^2$. 
  \item[\it(b)]  the ONB is unique up to
    $\bt \to \eps \bt, \bv, \bw$ up to a proper ($\eps =1$) or
    improper ($\eps=-1$) Euclidean rotation in $\R^2$.
  \end{description}
\item 
	\begin{description} \item [\it(a)] 
	 every ONB or LVB.
  \item[\it(b)] the ONB is unique
	up to $\bv \to \eps \bv, \bt, \bw$ up to a proper ($\eps =1$)
	or improper ($\eps=-1$) Pseudo-Euclidean rotation in
	$\R^2_1$,

	the LVB is unique up to $\begin{cases} \bv \to & \bv,\\ \be \to
	& a \be,\\ \bof \to & \frac1a \bof, \end{cases}$ or $\begin{cases} \bv
	\to & - \bv,\\ \be \to & a \bof,\\ \bof \to & \frac1a \be,
	\end{cases}\quad a\in \R$,
\end{description}
\item 
  \begin{description} \item [\it(a)] the LVB is unique up to $\begin{cases} \be \to &a \be,\\
    \bof \to &\frac1a \bof,
    \end{cases} \quad a\in \R.$
    Under $\begin{cases} \bv
	\to & - \bv,\\ \be \to & a \bof,\\ \bof \to & \frac1a \be,
	\end{cases}$, $C_\la$ goes to $C_\frac1\la$.
    \item [\it(b)] 
    the LVB is unique up to 
    $\begin{cases} \be \to & \be,\\ \bv \to & b \be + \bv,\\ 
      \bof \to & -\frac{b^2}{2} \be - b \bv +\bof, \end{cases} \quad
    b\in \R.$
    Under $\begin{cases} \be \to & a\be,\\ \bv \to & b \be + \bv,\\ 
      \bof \to & -\frac{b^2}{2a} \be - \frac{b}{a} \bv +\frac1{a} \bof,
    \end{cases}$, $C_1$ goes to $C_{1,\frac1{a}}$.
  \end{description}
\end{enumerate}
\end{rem} 

\section{Non-trivial $SO(1,2)$-stabilizers}
\label{sec:max}

Since we are interested in pointwise group
    symmetry, we study the nontrivial stabilizer of a traceless cubic
    form $\tilde{K}$ under the $\SO(1,2)$-action
    $\rho (L)(\tilde{K}) = \tilde{K}\circ L $ resp.
$$\rho (L)(h (K(\,.\, ,\,.\,),\,.\,)) = h (K(L\,.\,,L\,.\,) ,L\,.))$$
(cp. \cite{Bry01} for the classification under the $\SO (3)$-action).
We will use the following notation for the coefficients of the 
difference tensor $K$ with respect to an ONB $\{ \bt, \bv, \bw\}$:
\begin{equation}\label{ONBcoeffK}
\begin{split}
K_{\bt}&=\begin{pmatrix} -a_1 & -a_2 & -a_3\\
                        a_2 & a_4 & a_5\\
                        a_3 & a_5 & a_1-a_4  \end{pmatrix},\quad
K_{\bv}=\begin{pmatrix} -a_2 & -a_4 & -a_5\\
                        a_4 & a_6 & a_7\\
                        a_5 & a_7 & a_2-a_6  \end{pmatrix},\\
K_{\bw}&=\begin{pmatrix} -a_3 & -a_5 & -(a_1-a_4)\\
                        a_5 & a_7 & a_2-a_6\\
                        a_1-a_4 & a_2-a_6 & a_3-a_7  \end{pmatrix},
\end{split}\end{equation}

resp. with respect to a LVB $\{\be,\bv,\bof\}$:

\begin{equation}\label{LVBcoeffK}
K_{\be}=\begin{pmatrix} b_1 & b_4 & b_5\\
                         b_2 & -2 b_1 & b_4\\
                         b_3 & b_2 & b_1 \end{pmatrix},
K_{\bv}=\begin{pmatrix}  b_4 & -2 b_5 & b_6\\
                        -2 b_1 & -2 b_4 & -2 b_5\\
                         b_2 & -2 b_1 & b_4  \end{pmatrix},
K_{\bof}=\begin{pmatrix}  b_5 & b_6 & b_7 \\
                         b_4 & -2 b_5 & b_6\\
                         b_1 & b_4 & b_5 \end{pmatrix}.
\end{equation}

\noindent We will prove the following theorem, stating not only
the non-trivial stabilizers, but also give a normal form of $K$ for
each stabilizer.

\begin{thm}\label{thm:Subgroups}
Let $p\in M$ and assume that there exists a non-trivial element of
$SO(1,2)$ which preserves $K$. Then there exists an ONB resp. a LVB of
$T_p M$ such that either
\begin{enumerate}
\item $K=0$, and this form is preserved by every isometry, or
\item $a_1=2 a_4, a_4>0$, and all other coefficients vanish, this
form is preserved by the subgroup $\{A_t, t\in \R\}$, 
isomorphic to $SO(2)$, or
\item $a_6 >0$, and all other coefficients vanish, this form is
preserved by the subgroup with generators $< A_\frac{2 \pi}{3}, B>$,
isomorphic to $S_3$, or
\item $a_1=2 a_4, a_4> 0$, $a_6 >0$, and all other
coefficients vanish, this form is preserved by the subgroup with
generator $< A_\frac{2 \pi}{3}>$, isomorphic to $\Z_3$, or
\item $a_2, a_5 \in\R$, $a_6 \geq 0$, where
$(a_2,a_6)\neq 0$, and all other coefficients vanish, this form is
preserved by the subgroup with generator $< B>$, isomorphic to $\Z_2$, or
\item $a_5> 0$, and all other coefficients vanish, this form is
preserved by the subgroup with generators $< A_\pi, B>$, isomorphic to
$\Z_2 \times \Z_2$, or
\item $a_1> 0$ or $a_4> 0$, $a_1\neq 2 a_4$,
and all other coefficients vanish, this form is preserved by the
subgroup with generator $< A_\pi>$, isomorphic to $\Z_2$, or
\item $b_4 >0$, and all other coefficients vanish, this form is
preserved by the subgroup $\{C_l, l\in \R\setminus\{0\}\}$, isomorphic
to $SO(1,1)$, or
\item $b_7>0$, and all other coefficients vanish, this form is
preserved by the subgroup $\{\left(\begin{smallmatrix} 1 & -m &
-\frac{m^2}{2} \\ 0 & 1 & m\\ 0 & 0 & 1 \end{smallmatrix}\right),
m\in \R\}$, isomorphic to $\R$.
\end{enumerate}
\end{thm}

\noindent To get ready for the proof, first we will find out what it
means for $K$ to be invariant under one element of $\SO(1,2)$. Some of
the computations were done with the CAS
Mathematica\footnote{http://www.math.tu-berlin.de/$\sim$schar/IndefSym\_Stabilizers.html}.

\begin{lem}\label{lem:CoeffA}
Let $p\in M$ and assume, that $K$ is invariant under the
transformation $A_t \in \SO(1,2)$, $t\in (0, 2\pi)$. Then we get for the
coefficients of $K$ with respect to the corresponding ONB of
$T_pM$:
\begin{enumerate}
\item if $t\neq \frac{2\pi}{3},\pi,\frac{4\pi}{3}$, then $a_1=2
a_4$, $a_4\in\R $, and all other coefficients vanish,
\item if $t= \frac{2\pi}{3}$ or $t= \frac{4\pi}{3}$, then $a_1=2 a_4$,
$a_4, a_6, a_7 \in \R$, and all other coefficients vanish,
\item if $t=\pi$, then $a_1, a_4, a_5 \in \R$, and all other
coefficients vanish,
\end{enumerate}
\end{lem}

\begin{proof} The proof is a straight forward computation, evaluating 
the equations $h(K(X,Y),Z)=h(K(A_t X,A_t Y),A_t Z)$ for $X,Y,Z\in
\{\bt,\bv,\bw\}$. The computations were done with the CAS
Mathematica. For all $t\in (0,2\pi)$ we obtain from eq2 ($X,Y=\bt$,
$Z=\bv$) and eq3 ($X,Y=\bt$, $Z=\bw$) that $a_2=0$ and $a_3=0$. If
$t=\pi$, then eq7 ($X,Y,Z=\bv$) and eq8 ($X,Y=\bv$, $Z=\bw$) give
$a_6=0$ and $a_7=0$. If $t\neq \pi$, then eq4 ($X=\bt$, $Y,Z=\bv$) and
eq5 ($X=\bt$, $Y=\bv$, $Z=\bw$) lead to $a_5=0$ and $a_1=2 a_4$. Now,
for $t=\frac{2\pi}{3}$ or $t= \frac{4\pi}{3}$, all equations are
true. Otherwise, only $a_6=0$ and $a_7=0$ solve eq7 and eq8.
\end{proof}

\begin{lem}\label{lem:CoeffB}
Let $p\in M$ and assume, that $K$ is invariant under the
transformation $B \in \SO(1,2)$. Then we get for the coefficients of
$K$ with respect to the corresponding ONB of $T_pM$ that $a_2, a_5,
a_6 \in \R$, and all other coefficients vanish.
\end{lem}

\begin{proof} The computations were done with the CAS
Mathematica, too. We obtain from eq1 ($X,Y,Z=\bt$), eq3 ($X,Y=\bt$,
$Z=\bw$), eq4 ($X=\bt$, $Y,Z=\bv$) and eq8 ($X,Y=\bv$, $Z=\bw$) that
$a_1=0$, $a_3=0$, $a_4=0$ and $a_7=0$.
\end{proof}

\begin{lem}\label{lem:CoeffClLVB}
Let $p\in M$ and assume, that $K$ is invariant under the
transformation $C_l \in \SO(1,2)$, $l\in \R\setminus\{0, 1\}$. Then
we get for the coefficients of $K$ with respect to the corresponding
LVB of $T_pM$:
\begin{enumerate}
\item if $l\neq -1$, then $b_4\in \R$, and all other coefficients vanish,
\item if $l=-1$, then $b_2, b_4, b_6\in \R$, and all other
coefficients vanish.
\end{enumerate}
\end{lem}

\begin{proof} The computations were done with the CAS
Mathematica, too. We obtain from eq1 ($X,Y,Z=\be$), eq3 ($X,Y=\be$,
$Z=\bof$), eq6 ($X=\be$, $Y,Z=\bof$) and eq10 ($X,Y,Z=\bof$) that
$b_1=0$, $b_3=0$, $b_5=0$ and $b_7=0$. If $l\neq -1$, then eq2
($X,Y=\be$, $Z=\bv$) and eq9 ($X=\bv$, $Y,Z=\bof$) additionally give
that $b_2=0$ and $b_6=0$.
\end{proof}

\begin{lem}\label{lem:CoeffC1m}
Let $p\in M$ and assume, that $K$ is invariant under the
transformation $C_{1,m} \in \SO(1,2)$, $m\in \R\setminus\{0\}$. Then
we get for the coefficients of $K$ with respect to the corresponding
LVB of $T_pM$ that $b_7\in \R$, and all other coefficients vanish.
\end{lem}

\begin{proof} The computations were done with the CAS
Mathematica, too. We obtain successively from eq2 ($X,Y=\be$,
$Z=\bv$), eq3 ($X,Y=\be$, $Z=\bof$), eq5 ($X=\be$, $Y=\bv$, $Z=\bof$),
eq6 ($X=\be$, $Y,Z=\bof$), eq9 ($X=\bv$, $Y,Z=\bof$) and eq10
($X,Y,Z=\bof$) that $b_3=0$, $b_2=0$, $b_1=0$, $b_4=0$, $b_5=0$ and
$b_6=0$.
\end{proof}

\noindent In the following $U$ denotes an arbitrary subgroup of $\SO(1,2)$,
which leaves $K$ invariant. We want to find out to which extend $K$
determines the properties of the elements of $U$.

\begin{lem}\label{lem:PropA}
If there exists $t\in (0,2 \pi)$, $t\neq \pi$, with $A_t \in U$, and
$K\neq 0$, then we get for the timelike eigenvector $\bt$ of $A_t$:
\begin{enumerate}
\item for $t\neq \frac{2\pi}{3}$ and $t\neq \frac{4\pi}{3}$:
$M\bt=\bt$ for all $M\in U$,
\item for $t= \frac{2\pi}{3}$ or $t= \frac{4\pi}{3}$: $M\bt=\eps\bt$
for all $M\in U$,
\end{enumerate}
\end{lem}

\begin{proof} Let $M \in U$. From Lem.~\ref{lem:CoeffA} we know that 
$$h(K(M\bt,M \bt), MY) = h(K(\bt, \bt), Y)= \begin {cases} -2a_4, &
Y=\bt,\\ 0, &Y=\bv \;\text{or}\; Y=\bw. \end{cases}$$ Thus
$K(M\bt,M\bt)=-2 a_4 M\bt$, furthermore $h(
M\bt,M\bt)=h(\bt,\bt)=-1$. Now assume that $X=x \bt +y \bv +z\bw \in
T_p M$ has the same properties ($K(X,X)=-2 a_4 X$ and $h( X,X)=
-1$). This is equivalent to (cp. Lem.~\ref{lem:CoeffA}):
\begin{eqnarray}\label{KXX1}  (-2 x^2-y^2-z^2) a_4 &= -2 a_4 x,\\
\label{KXX2} 2xy a_4 +(y^2-z^2) a_6+ 2yz a_7&= -2 a_4 y,\\
\label{KXX3} 2xz a_4 - 2yz a_6 +(y^2-z^2)a_7 &= -2 a_4 z, \\
\label{Xtimelike} -x^2+y^2+z^2 &= -1. 
\end{eqnarray}
If $a_4\neq 0$, then \eqref{KXX1} and \eqref{Xtimelike} imply that
$3x^2-2x-1=0$ and $x^2 \geq1$, this means $x=1$ and $y=z=0$. Thus
$X=\bt$ and $M\bt=\bt$.

If $a_4=0$, then \eqref{KXX2} and \eqref{KXX3} imply that 
\begin{eqnarray*} (y^2 -z^2)a_6 + 2yz a_7 & =0,\\
-2yz a_6 +(y^2-z^2) a_7 &=0.
\end{eqnarray*}
The two equations, linear in $a_6$ and $a_7$, only have a non-trivial
solution if $y=z=0$. With \eqref{Xtimelike} we obtain that $X=\eps
\bt$ and thus $M\bt=\eps\bt$.
\end{proof}

\begin{lem}\label{lem:PropCl}
If there exists $l\in \R$, $l\neq 0,\pm1$, with $C_l \in U$, and
$K\neq0$, then we get for the spacelike eigenvector $\bv$ of $C_l$:
$M\bv=\bv$ for all $M\in U$.
\end{lem}

\begin{proof} Let $M \in U$. From Lem.~\ref{lem:CoeffClLVB} we know that 
$$h( K(M\bv,M \bv), MY) = h( K(\bv, \bv), Y)= \begin {cases} -2a_4, &
Y=\bv,\\ 0, &Y=\be \;\text{or}\; Y=\bof. \end{cases}$$ Thus
$K(M\bv,M\bv)=-2 a_4 M\bv$, furthermore $h(
M\bv,M\bv)=h(\bv,\bv)=1$. Now assume that $X=x \bt +y \bv +z\bw \in
T_p M$ has the same properties ($K(X,X)=-2 a_4 X$ and $h( X,X)=
1$). This is equivalent to (cp. Lem.~\ref{lem:CoeffClLVB}):
\begin{eqnarray}\label{KXX_LVB1}  2xy b_4 &= -2 b_4 x,\\
\label{KXX_LVB2} 2(-y^2+xz) b_4 &= -2 b_4 y,\\
\label{KXX_LVB3} 2yz b_4 &= -2 b_4 z, \\
\label{Xspacelike} 2xz + y^2 &= 1. 
\end{eqnarray}
Since $b_4\neq 0$, \eqref{KXX_LVB1} is equivalent to $x=0$ or $y=-1$,
and \eqref{KXX_LVB3} is equivalent to $z=0$ or $y=-1$. Now $y=-1$ in
\eqref{KXX_LVB2} gives $xy=2$, which is a contradiction to
\eqref{Xspacelike}. Thus $x=0$ and $z=0$. With \eqref{KXX_LVB2} we
obtain that $X= \bv$ and thus $M\bv=\bv$.
\end{proof}

\begin{lem}\label{lem:PropC1m}
If there exists $m\in \R$, $m\neq 0$, with $C_{1,m} \in U$, and $K\neq
0$, then we get for the lightlike eigenvector $\be$ of $C_{1,m}$:
$M\be=\be$ for all $M\in U$.
\end{lem}

\begin{proof} 
Let $M \in U$. From Lem.~\ref{lem:CoeffC1m} we know that \[K(X,X)=b_7
z^2 \be \;\text{for all}\; X=x \bt +y \bv +z\bw\in T_pM,\] i.~e. $\be$
is determined by $K$ up to length: $\be=\frac{1}{h(
K(\bof,\bof),\bof)} K(\bof,\bof)$. From the invariance of $K$ under
$M$ it follows that $M\be= \frac{1}{h( K(\bof,\bof),\bof)}
K(M\bof,M\bof) = (h( M\bof,\be))^2 \be$. Since $1= h( M\bof,M\be)$, we
get now: $1=h( M\bof,h( M\bof,\be)^2 \be) =h( M\bof,\be)^3$, thus $h(
M\bof,\be)=1$.
\end{proof}

\noindent Now we are ready for the proof of Thm.~\ref{thm:Subgroups}.

\begin{proof}[Proof of Theorem~\ref{thm:Subgroups}]
For the proof we will consider several cases which are supposed to be
exclusive. Let $U$ be a maximal subgroup of $\SO(1,2)$ which leaves
$K$ invariant.

{\bf 1. Case} We assume that there exists $t\in (0,2 \pi)$, $t\neq
\frac{2 \pi}{3}, \pi, \frac{4 \pi}{3}$, with $A_t \in U$. Thus there
exists an ONB $\{\bt,\bv,\bw\}$ such that $K$ has the form
(Lem.~\ref{lem:CoeffA}): $a_1=2 a_4$, $a_4\in\R$ and all other
coefficients vanish. If $a_4=0$, then $K=0$ and $U=\SO(1,2)$. If
$a_4\neq 0$, then we know that $M\bt=\bt$ for all $M\in U$
(Lem.~\ref{lem:PropA}). We have seen in Sec.~\ref{sec:nf} that $M\in
U$ must be of type $\it 1 (a)$, $\it 1 (b)$ or $\it 2 (a)$. Now let
$N\in \SO(1,2)$ be of type $\it 1 (a)$, $\it 1 (b)$ or $\it 2 (a)$
with eigenvector $\bt$. We can normalize simultaneously (i.~e. find
an ONB such that both $A_t$ and $N$ have normalform) and we see that
$N$ leaves $K$ invariant (Lem.~\ref{lem:CoeffA}). Thus $U=\{A_t,
t\in\R\}$. Finally, if $a_4<0$, we can change the ONB
$\{\bt,\bv,\bw\}$ and take instead$\{-\bt,\bw,\bv\}$.

{\bf 2. Case} We assume that $A_\frac{2 \pi}{3} \in U$ or $A_\frac{4
\pi}{3} \in U$. Thus there exists an ONB $\{\bt,\bv,\bw\}$ such that
$K$ has the form (Lem.~\ref{lem:CoeffA}): $a_1=2 a_4$, $a_4, a_6, a_7
\in\R$ and all other coefficients vanish. Without loss of generality
$a_6$ and $a_7$ will not vanish both. Furthermore we know that
$M\bt=\eps \bt$ for all $M\in U$ (Lem.~\ref{lem:PropA}). We have seen
in Sec.~\ref{sec:nf} that $M\in U$ must be of type $\it 1 (a)$, $\it 1
(b)$, $\it 2 (a)$ or $\it 2 (b)$. Now let $N\in \SO(1,2)$ be of type
$\it 1 (a)$, $\it 1 (b)$, $\it 2 (a)$ or $\it 2 (b)$ with eigenvector
$\bt$. We can normalize simultaneously.

If $N= A_\frac{2 \pi}{3}$, $N= A_\frac{4 \pi}{3}$ or $N= \Id$ (type
$\it 1 (a)$ or $\it 2 (a)$), then it leaves $K$ invariant
(Lem.~\ref{lem:CoeffA}). If $N=A_\pi$ (type $\it 1 (b)$), then by
Lem.~\ref{lem:CoeffA} $a_6=0=a_7$, which gives a contradiction. If
$N=B$ (type $\it 2 (b)$), then by Lem.~\ref{lem:CoeffB} $a_4=0=a_7$
and $a_6$ is the only non-vanishing coefficient of $K$.

We get two possibilities for $K$ and the corresponding maximal
subgroup $U$. Either $a_6 \in \R\setminus\{0\}$ and all other
coefficients of $K$ vanish, and $U=< A_\frac{2 \pi}{3}, B>$, if
necessary by a change of basis ($\{-\bt,-\bv,\bw\}$) we can make sure
that $a_6>0$. Or $a_1=2 a_4$, $a_4, a_6, a_7 \in\R$ and all other
coefficients vanish, and $ U=< A_\frac{2 \pi}{3}>$. As before we can
choose $\bt$ such that $a_4\geq 0$. A computation gives that under a
change of ONB $\{\bt^*,\bv^*,\bw^*\} = \{\bt,\cos s \bv +\sin s \bw,
-\sin s \bv +\cos s \bw\}$ we obtain for $K$ (cp. \eqref{ONBcoeffK}):
$a_6^*= a_6 \cos(3s) + a_7 \sin(3s)$, $a_7^*= a_7 \cos(3s) - a_6
\sin(3s)$, i.~e. there exists $s\in\R$ such that $a_7^*=0$. We can
change the sign of $a_6$ by switching from $\{\bt,\bv,\bw\}$ to
$\{\bt,-\bv,-\bw\}$. Finally we see that $a_4\neq 0$.

{\bf 3. Case} We assume that there exists $l\neq 0, \pm 1$, with $C_l
\in U$.  There exists a LVB $\{\be,\bv,\bof\}$ such that $K$ has the
form (Lem.~\ref{lem:CoeffClLVB}): $b_4\in\R$ and all other
coefficients vanish. If $b_4=0$, then $K=0$ and $U=\SO(1,2)$.

If $b_4\neq 0$, then we know that $M\bv=\bv$ for all $M\in U$
(Lem.~\ref{lem:PropCl}). We have seen in Thm.~\ref{thm:NF} that $M\in
U$ must be of type $\it 2 (a)$, $\it 2 (b)$ or $\it 3 (a)$. Now let
$N\in \SO(1,2)$ be of type $\it 2 (a)$, $\it 2 (b)$ or $\it 3 (a)$
with eigenvector $\bv$. We can normalize simultaneously and we see
that $N$ leaves $K$ invariant (Lem.~\ref{lem:CoeffClLVB},
$B=C_{-1}$). Thus $U=\{C_l, l\in\R\setminus\{0\}\}$. Finally, if
$b_4<0$, we can change the LVB $\{\be,\bv,\bof\}$ to
$\{\bof,-\bv,\be\}$ (cp. Rem.~\ref{rem:NF}).

{\bf 4. Case} We assume that $C_1 \in U$.  There exists a LVB
$\{\be,\bv,\bof\}$ such that $K$ has the form
(Lem.~\ref{lem:CoeffC1m}): $b_7\in\R$ and all other coefficients
vanish. If $b_7=0$, then $K=0$ and $U=\SO(1,2)$.

If $b_7\neq 0$, then we know that $M\be=\be$ for all $M\in U$
(Lem.~\ref{lem:PropC1m}). We have seen in Thm.~\ref{thm:NF} that $M\in
U$ must be of type $\it 3 (b)$. Now let $N\in \SO(1,2)$ be of type
$\it 3 (b)$ with eigenvector $\be$, i.~e. $N$ has the form $C_{1,m}$,
$m\in \R$. We see that $N$ leaves $K$ invariant
(Lem.~\ref{lem:CoeffC1m}). Thus $U=\lp C_1 \rp$ ($C_{1,m}
C_{1,n}=C_{1,m+n}$). Finally, if $b_7<0$, we can change the LVB
$\{\be,\bv,\bof\}$ to $\{-\be,\bv,-\bof\}$ (cp. Rem.~\ref{rem:NF}).

{\bf 5. Case} We assume that $A_\pi \in U$.  There exists an ONB
$\{\bt,\bv,\bw\}$ such that $K$ has the form (Lem.~\ref{lem:CoeffA}):
$a_1, a_4, a_5 \in\R$ and all other coefficients vanish. Every $M\in
U$ must be of type $\it 1 (b)$ or $\it 2 (b)$, otherwise we are in one
of the foregoing cases.

a) Let $N\in \SO(1,2)$ be of type $\it 1 (b)$. There exists a
spacelike eigenvector $\bw$ such that $A_\pi \bw= -\bw$ and $N \bw=
-\bw$, and we can choose an ONB $\{\bt,\bv,\bw\}$ such that $A_\pi$
has normalform and $N= \left(\begin{smallmatrix} \cosh s & -\sinh s &0
\\ \sinh s& -\cosh s & 0\\0 & 0 & -1 \end{smallmatrix}\right)$. If we
assume that $K$ is invariant under $N$ then we obtain for $s\neq 0$
that $K=0$.  The computations were done with the CAS Mathematica,
too. We obtain from eq3 ($X,Y=\bt$, $Z=\bw$) and eq9 ($X=\bv$,
$Y,Z=\bw$) that $a_5=0$ and $a_1=a_4$, and from eq1 ($X,Y,Z=\bt$) that
$a_4=0$.

b) Let $N\in \SO(1,2)$ be of type $\it 2 (b)$. There exists a
spacelike eigenvector $\bw$ such that $A_\pi \bw= -\bw$ and $N \bw=
-\bw$, and we can choose an ONB $\{\bt,\bv,\bw\}$ such that $A_\pi$
has normalform and $N= \left(\begin{smallmatrix} -\cosh s & -\sinh s
&0 \\ \sinh s& \cosh s & 0\\0 & 0 & -1 \end{smallmatrix}\right)$. If
we assume that $K$ is invariant under $N$ then we obtain that $a_1=0
=a_4$. If $s\neq 0$, also $a_5=0$, i.~e. $K=0$.  The computations were
done with the CAS Mathematica, too. If $s\neq 0$, we get from eq3
($X,Y=\bt$, $Z=\bw$) and eq9 ($X=\bv$, $Y,Z=\bw$) that $a_5=0$ and
$a_1=a_4$, and from eq1 ($X,Y,Z=\bt$) that $a_4=0$. If $s=0$, we get
from eq1 ($X,Y,Z=\bt$) and eq4 ($X=\bt$, $Y,Z=\bv$) that $a_1=0=a_4$.

Summarized we got two different forms of $K$ with
corresponding maximal subgroups $U$: a) Either there exists an ONB
such that $a_1, a_4, a_5 \in\R$, where $a_1\neq 2 a_4$ or $a_5\neq 0$,
and all other coefficients vanish, this form is preserved by $U=<
A_\pi>$. Since $K_\bt$ is a symmetric operator on the positive
definite space $\bt^{\bot}$, we can diagonalize, then $a_5=0$. If
necessary, we still can take $\{-\bt,-\bv,\bw\}$ to get $a_1>0$ or
$a_4>0$. b) In the other case there exists an ONB such that $a_5\in
\R$, $a_5\neq 0$, and all other coefficients vanish, this form is
preserved by $U=< A_\pi, B>$. If $a_5<0$, we switch to the ONB
$\{-\bt,\bw,\bv\}$.

{\bf 6. Case} We assume that $B \in U$. There exists an ONB
$\{\bt,\bv,\bw\}$ such that $K$ has the form (Lem.~\ref{lem:CoeffB}):
$a_2, a_5, a_6 \in\R$ and all other coefficients vanish. Every $M\in
U$ must be of type $\it 2 (b)$, otherwise we are in one of the
foregoing cases. Now let $N\in \SO(1,2)$ be of type $\it 2 (b)$. We
can't normalize simultaneously. We only know that $B$ and $N$ both
have two-dimensional timelike eigenspaces, which intersect in a
line. This line $g$ can be space-, time- or lightlike.

a) If $g$ is spacelike, we can choose an ONB $\{\bt,\bv,\bw\}$ such
that $B$ has normalform and $N= \left(\begin{smallmatrix} -\cosh s &
-\sinh s &0 \\ \sinh s& \cosh s & 0\\0 & 0 & -1
\end{smallmatrix}\right)$ (cp. case 5). If we assume that $K$ is
invariant under $N$ then we obtain for $s\neq 0$ that $K=0$. The
computations were done with the CAS Mathematica, too. If $s\neq 0$, we
get from eq3 ($X,Y=\bt$, $Z=\bw$) and eq9 ($X=\bv$, $Y,Z=\bw$) that
$a_5=0$ and $a_2=a_6$, and from eq1 ($X,Y,Z=\bt$) that $a_6=0$.

b) If $g$ is timelike, we can choose an ONB $\{\bt,\bv,\bw\}$ such
that $B$ has normalform and $N= \left(\begin{smallmatrix} -1 & 0 &0 \\
0& \cos s & \sin s\\0 & \sin s & -\cos s \end{smallmatrix}\right)$. If
we assume that $K$ is invariant under $N$ then we obtain for $s\neq 0,
\frac{2\pi}{3}, \pi, \frac{4\pi}{3}$ that $K=0$. For
$s=\frac{2\pi}{3}, \pi, \frac{4\pi}{3}$ we are in one of the foregoing
cases. The computations were done with the CAS Mathematica, too. If
$s\neq 0, \pi$, we get from eq2 ($X,Y=\bt$, $Z=\bv$) and eq4 ($X=\bt$,
$Y,Z=\bv$) that $a_2=0$ and $a_5=0$. If also $s\neq \frac{2\pi}{3},
\frac{4\pi}{3}$, then we get from eq10 ($X,Y,Z=\bw$) that $a_6=0$.

c) If $g$ is lightlike, we can choose a LVB $\{\be,\bv,\bof\}$ such
that $B$ has normalform and $N= \left(\begin{smallmatrix} -1 & -2 m &
2 m^2 \\ 0& 1 & -2 m\\0 & 0 & -1 \end{smallmatrix}\right)$ (use
Lem.~\ref{lem:LVBtoLVB}). If we assume that $K$ is invariant under $N$
then we obtain for $m\neq 0$ that $K=0$. The computations were done
with the CAS Mathematica, too. If $m\neq 0$, we get from eq5 ($X=\bt$,
$Y=\bv$, $Z=\bw$) and eq2 ($X,Y=\bt$, $Z=\bv$) that $b_6=0$ and
$b_4=0$, and from eq1 ($X,Y,Z=\bt$) that $b_2=0$.

Therefore we have that $U=\lp B\rp$. If $a_6<0$, we can switch to
$\{-\bt,-\bv,\bw\}$. Since $K_\bv$ is a symmetric operator on an
indefinite space we can't always diagonalize. Thus we can't simplify
$K$ in general.  \end{proof} \begin{rem} In the proof we only have
used multilinear algebra. Thus the theorem stays true for an arbitrary
$(1,2)-$tensor $K$ on $\Min$ with $\lp K(X,Y),Z\rp$ totally symmetric and
vanishing $\trace K_X$.  \end{rem} 

\section{Pointwise $\Z_2\times\Z_2$-symmetry}
\label{sec:type6}

Let $\M$ be a hypersphere admitting a pointwise $\Z_2 \times
\Z_2$-symmetry. According to Thm.~\ref{thm:Subgroups}, there exists
for every $p\in \M$ an ONB $\{ \bt, \bv, \bw\}$ of $T_p \M$ such that
\begin{align} K(\bt,\bt)&= 0, & K(\bt,\bv)&=a_5 \bw, & K(\bt,\bw)&= a_5 \bv,\\
K(\bv,\bv)&= 0, & K(\bv,\bw)&= -a_5 \bt, & K(\bw,\bw)&= 0.
\end{align}
Substituting this in Eq.~\eqref{gaussLC}, we obtain
\begin{equation}
\hat R (X,Y)Z =(H- a_5^2) (h(Y,Z) X - h(X,Z) Y).
\end{equation}
Schur's Lemma implies that $\M$ has constant sectional curvature, by
the affine theorema egregium \eqref{TE} we obtain that $\hat
\kappa=H-a_5^2$ and $J=-a_5^2<0$. Affine hyperspheres with constant
affine sectional curvature and nonzero Pick invariant were classified
by Magid and Ryan \cite{MR92}. They show in their main theorem that an
affine hypersphere with Lorentz metric of constant curvature and
nonzero Pick invariant is equivalent to an open subset of either
$(x_1^2+x_2^2)(x_3^2+x_4^2)=1$ or $(x_1^2+x_2^2)(x_3^2-x_4^2)=1$. In
both cases $\hat\kappa =0$, i.~e. $H=-J$. In the proof of the main
theorem they explicitly show that only $(x_1^2+x_2^2)(x_3^2+x_4^2)=1$
has negative Pick invariant and that $K$ has normalform. This
proves:
\begin{thm}\label{thm:type6}
An affine hypersphere admits a pointwise $\Z_2 \times \Z_2$-symmetry
if and only if it is affine equivalent to an open subset of
$$(x_1^2+x_2^2)(x_3^2+x_4^2)=1.$$
\end{thm}
For $(x_1^2+x_2^2)(x_3^2-x_4^2)=1$ they compute that the only
non-vanishing coefficient of $K$ is $a_2$. Thus it follows
(Thm.~\ref{thm:Subgroups}):
\begin{rem} 
The affine hypersphere ($x_1^2+x_2^2)(x_3^2-x_4^2)=1$ admits a
pointwise $\Z_2$-symmetry.
\end{rem}

\section{Pointwise $\R$-symmetry}
\label{sec:type9}

Let $\M$ be a hypersphere admitting a pointwise
$\R$-symmetry. According to Thm.~\ref{thm:Subgroups}, there exists for
every $p\in \M$ a LVB $\{ \be, \bv, \bof\}$ of $T_p \M$ such that
\begin{align} K(\be,\be)&= 0, & K(\be,\bv)&=0, & K(\be,\bof)&= 0,\\
K(\bv,\bv)&= 0, & K(\bv,\bof)&= 0, & K(\bof,\bof)&= b_7 \be.
\end{align}
Substituting this in Eq.~\eqref{gaussLC}, we obtain
\begin{equation}
\hat R (X,Y)Z =H (h(Y,Z) X - h(X,Z) Y).
\end{equation}
Schur's Lemma implies that $\M$ has constant sectional curvature, by
the affine theorema egregium \eqref{TE} we obtain that $\hat \kappa=H$
and $J=0$. Affine hyperspheres with constant affine sectional
curvature and zero Pick invariant were classified in \cite{DMV00} (see
Thm. 6.2 ($H=0$), Thm. 7.2 ($H=1$) and Thm. 8.2 ($H=-1$)). They are
determined by a null curve in resp. $\R^3_1$, $S^3_1$, $H^3_1$, and a
function along this curve (note that in the notion of \cite{DMV00} (2)
holds).
\begin{thm}\label{thm:type9}
Let $\M$ be an affine hypersphere admitting a pointwise
$\R$-sym\-me\-try. Then $\M$ has constant sectional curvature $\hat
\kappa=H$ and zero Pick invariant $J=0$.
\end{thm}
\begin{rem}
A study of \cite{DMV00} shows that an affine hypersphere admits a
pointwise $\R$-symmetry if and only if (2) holds (in their notations).

If (3) holds for an affine hypersphere with constant sectional
curvature and zero Pick invariant, then it admits a pointwise
$\Z_2$-symmetry (cp. Thm.~\ref{thm:Subgroups}).
\end{rem}

\section{Pointwise $\SO(2)$-, $S_3$- or $\Z_3$-symmetry}
\label{sec:type2,3,4}

Let $\M$ be a hypersphere admitting a $\SO(2)$-, $S_3$- or
$\Z_3$-symmetry. According to Thm.~\ref{thm:Subgroups}, there exists
for every $p\in \M$ an ONB $\{ \bt, \bv, \bw\}$ of $T_p \M$ such that
\begin{align*} K(\bt,\bt)&= -2a_4 \bt, & K(\bt,\bv)&=a_4 \bv, & K(\bt,\bw)&= 
a_4 \bw,\\
K(\bv,\bv)&= -a_4 \bt + a_6 \bv, & K(\bv,\bw)&= -a_6 \bw, &
K(\bw,\bw)&= -a_4 \bt -a_6 \bv,
\end{align*}
where $a_4 > 0$ and $a_6=0$ in case of $SO(2)$-symmetry, $a_4=0$ and
$a_6>0$ for $S_3$, and $a_4>0$ and $a_6>0$ for $\Z_3$.

We would like to extend the ONB locally. It is well known that
$\ricLC$ (cp. \eqref{def:RicLC}) is a symmetric operator and we
compute (some of the computations in this section are done with the
CAS
Mathematica\footnote{http://www.math.tu-berlin.de/$\sim$schar/IndefSym\_typ234.html}):
\begin{lem} Let $p \in M$ and $\{\bt,\bv,\bw\}$ the basis constructed
  earlier. Then 
\begin{alignat*}{2}
  &\ricLC(\bt,\bt) =-2(H- 3a_4^2) , \qquad\quad &&\ricLC(\bt,\bv)=0, \\
  &\ricLC(\bt,\bw)=0,\qquad\quad &&
\ricLC(\bv,\bv)=2(H-a_4^2+a_6^2) ,\\
  &\ricLC(\bv,\bw)=0,&&\ricLC(\bw,\bw)=2(H-a_4^2+a_6^2) .
\end{alignat*}
\end{lem}

\begin{proof}
The proof is a straight-forward computation using the Gauss
equation~\eqref{gaussLC}.  It follows e.~g. that
\begin{align*}
  \hat R(\bt,\bv)\bt&= H \bv -K_{\bt}(a_4\bv)+K_{\bv}(-2 a_4 \bt) = H
  \bv -a_4^2 \bv -2 a_4^2\bv \\ 
  &= (H- 3 a_4^2)\bv,\\ 
  \hat
  R(\bt,\bw)\bt&= H \bw -K_{\bt}(a_4 \bw) +K_{\bw}(-2 a_4 \bt) =H \bw
  - a_4^2 \bw -2 a_4^2 \bw\\ &= (H- 3 a_4^2)\bw,\\ 
  \hat
  R(\bt,\bv)\bw&=-K_{\bt}(-a_6 \bw)+K_{\bv}(a_4\bw)=0.
\end{align*}
From this it immediately follows that 
$$\ricLC(\bt,\bt) = -2(H-3a_4^2)$$
and
$$\ricLC(\bt,\bw)=0.$$
The other equations follow by similar computations.
\end{proof}

We want to show that the basis, we have constructed at each point $p$,
can be extended differentiably to a neighborhood of the point $p$ such
that, at every point, $K$ with respect to the frame $\{T,V,W\}$ has
the previously described form.
\begin{lem}\label{lem:KfT234} 
Let $\M$ be an affine hypersphere in $\mathbb R^4$ which admits a
  pointwise $\SO(2)$-, $S_3$- or $\Z_3$-symmetry. Let $p \in M$. Then
  there exists an orthonormal frame $\{T,V,W\}$ defined in a neighborhood of the point $p$ such that $K$ is given by:
\begin{align*} 
K(T,T)&= -2a_4 T, & K(T,V)&=a_4 V, & K(T,W)&= 
a_4 W,\\
K(V,V)&= -a_4 T + a_6 V, & K(V,W)&= -a_6 W, &
K(W,W)&= -a_4 T -a_6 V,
\end{align*}
where $a_4 > 0$ and $a_6=0$ in case of $\SO(2)$-symmetry, $a_4=0$ and
$a_6>0$ in case of $S_3$-symmetry, and $a_4>0$ and $a_6>0$ in case of
$\Z_3$-symmetry.
\end{lem}

\begin{proof}First we want to show that at every point the vector
  $\bt$ is uniquely defined (up to sign) and differentiable. We
  introduce a symmetric operator $\hat A$ by:
\begin{equation*}
\ricLC(Y,Z)= h(\hat A Y,Z).
\end{equation*}
Clearly $\hat A$ is a differentiable operator on $M$. Since $2(H-
3a_4^2) \neq 2(H-a_4^2+a_6^2)$, the operator has two distinct
eigenvalues. A standard result then implies that the
eigendistributions are differentiable. We take $T$ a local unit
vectorfield spanning the 1-dimensional eigendistribution, and local
orthonormal vectorfields $\tilde{V}$ and $\tilde{W}$ spanning the
second eigendistribution. If $a_6=0$, we can take $V=\tilde{V}$ and
$W= \tilde{W}$.

As $T$ is (up to sign) uniquely determined, for $a_6\neq 0$ there
exist differentiable functions $a_4$, $c_6$ and $c_7$,
$c_6^2+c_7^2\neq 0$, such that
\begin{align*} 
K(T,T)&= -2a_4 T, & K(\tilde{V},\tilde{V})&= -a_4 T + c_6 \tilde{V}+
c_7 \tilde{W}, \\ K(T,\tilde{V})&=a_4 \tilde{V},&
K(\tilde{V},\tilde{W})&= c_7 \tilde{V} -c_6 \tilde{W}, \\
K(T,\tilde{W})&= a_4 \tilde{W},& K(\tilde{W},\tilde{W})&= -a_4 T -c_6
\tilde{V} -c_7 \tilde{W}.
\end{align*}
As we have shown in the proof of Thm.~\ref{thm:Subgroups} (Case 2), we
can always rotate $\tilde{V}$ and $\tilde{W}$ such that we obtain the
desired frame.
\end{proof}
\begin{rem} It actually follows from the proof of the previous lemma
  that the vector field $T$ is (up to sign) invariantly defined on
  $M$, and therefore the function $a_4$, too. Since the Pick invariant
  \eqref{def:Pick} $J= \frac{1}{3}(-5 a_4^2 + 2 a_6^2)$, the function
  $a_6$ also is invariantly defined on the affine hypersphere $\M$.
\end{rem}
In this section we always will work with the local frame constructed
in the previous lemma. We denote the coefficients of the Levi-Civita
connection with respect to this frame by:
\begin{align*}
   \widehat{\nabla}_T T &= \ac12 V + \ac13 W,&
   \widehat{\nabla}_T V &= \ac12 T - \bc13 W,&
   \widehat{\nabla}_T W &= \ac13 T + \bc13 V, \\
   \widehat{\nabla}_V T &= \ac22 V + \ac23 W, &
   \widehat{\nabla}_V V &= \ac22 T - \bc23 W, &
   \widehat{\nabla}_V W &= \ac23 T + \bc23 V, \\
   \widehat{\nabla}_W T &= \ac32 V + \ac33 W, &
   \widehat{\nabla}_W V &= \ac32 T - \bc33 W, &
   \widehat{\nabla}_W W &= \ac33 T + \bc33 V.
 \end{align*}
We will evaluate first the Codazzi and then the Gauss equations
(\eqref{CodK} and \eqref{gaussLC}) to obtain more informations.

\begin{lem}\label{lem:CodK}
Let $\M$ be an affine hypersphere in $\mathbb R^4$ which admits a
pointwise $\SO(2)$-, $S_3$- or $\Z_3$-symmetry and $\{T,V,W\}$ the
corresponding ONB. If the symmetry group is
\begin{description}
\item[$\mathbf{SO(2)}$,] then $0= \ac12 =\ac13 =\ac23 =\ac32$,
$\ac33=\ac22$ and \\ $T(a_4)=-4\ac22 a_4$, $0=V(a_4) = W(a_4)$,
\item[$\mathbf{S_3}$,] then $0= \ac12 =\ac13$, $ \ac23=-3 \bc13=
-\ac32$, $\ac33=\ac22$ and \\ $T(a_6)=-\ac22 a_6$, $V(a_6) = 3\bc33
a_6$, $W(a_6)=-3 \bc23 a_6$,
\item[$\mathbf{\Z_3}$ and $\mathbf{a_6^2\neq 4 a_4^2}$,] then $0=
\ac12 =\ac13 =\ac23 =\ac32$, $\ac33=\ac22$, $\bc13=0$,\\
$T(a_4)=-4\ac22 a_4$, $0=V(a_4) = W(a_4)$, and \\ $T(a_6)=-\ac22 a_6$,
$V(a_6) = 3\bc33 a_6$, $W(a_6)=-3 \bc23 a_6$,
\item[$\mathbf{\Z_3}$ and $\mathbf{a_6=2 a_4}$,] then $\ac12 =2 \ac22=
-2 \ac33=- \bc33$,\\ $\ac13 =- 2\ac23 = -2 \ac32= \bc23$, $\bc13=0$,
and \\ $T(a_4)=0$, $V(a_4)=-4 \ac22 a_4$, $W(a_4)= 4 \ac23 a_4$,
\end{description}
\end{lem}

\begin{proof}
An evaluation of the Codazzi equations \eqref{CodK} with the help of
the CAS
Mathematica leads to the following equations (they relate
to eq1--eq6 and eq8--eq9 in the Mathematica notebook):
\begin{eqnarray}
   V(a_4)=- 2\ac12 a_4, \quad T(a_4)=-4 \ac22 a_4 + \ac12 a_6, \quad
   0=4 \ac23 a_4 + \ac13 a_6,\label{CodKeq1} \\ W(a_4)=- 2\ac13
   a_4,\quad 0=4 \ac32 a_4 + \ac13 a_6,\quad T(a_4)=-4 \ac33 a_4 -
   \ac12 a_6, \label{CodKeq2}\\ T(a_6)-V(a_4)=3\ac12 a_4 -\ac22
   a_6,\quad 0=\ac13 a_4 + (\ac23 + 3 \bc13)a_6, \label{CodKeq3}\\
   \begin{split}W(a_4)=(\ac23+ \ac32) a_6,\quad  W(a_6)= 
(-\ac23+3 \ac32)a_4 - \bc23 a_6, \\ V(a_6)=(-\ac22+\ac33) a_4+ 3 \bc33
   a_6,
   \label{CodKeq4}\end{split}\\ T(a_6)=-\ac12 a_4 -\ac33 a_6,\quad
   W(a_4)=- 3\ac13 a_4+ (-\ac32+3 \bc13)a_6,\label{CodKeq5}\\
   V(a_4)=(- \ac22+\ac33) a_6, \quad W(a_6)= (3\ac23 - \ac32)a_4 -3
   \bc23 a_6,\label{CodKeq6}\\ 0=(\ac23-\ac32) a_4,\label{CodKeq8}\\
   W(a_4)=-\ac13 a_4 +(\ac32-3 \bc13) a_6.\label{CodKeq9}
 \end{eqnarray}
 
 From the first equation of \eqref{CodKeq2} (we will use the notation
 \eqref{CodKeq2}.1) and \eqref{CodKeq4}.1 resp. \eqref{CodKeq1}.3 and
 \eqref{CodKeq2}.2 we get:
 \begin{align}\label{a13,a23+a32}
 0&=2 \ac13 a_4 +(\ac23+\ac32)a_6\\
 0&=2 (\ac23+\ac32)a_4 + \ac13 a_6.
 \end{align}
  From \eqref{CodKeq6}.1) and \eqref{CodKeq1}.1
  resp. \eqref{CodKeq1}.2 and \eqref{CodKeq2}.3 we get:
 \begin{align}\label{a12,a22+a33}
 0&=-2 \ac12 a_4 +2 (\ac22-\ac33)a_6\\
 0&=2 (-\ac22+\ac33)a_4 + \ac12 a_6.
 \end{align}
 We consider first the case, that $\mathbf{a_6^2\neq 4 a_4^2}$. Then
 we obtain from the foregoing equations that $\ac13=0$,
 $\ac32=-\ac23$, $\ac12=0$ and $\ac33=\ac22$. Furthermore it follows
 from \eqref{CodKeq1}.1 that $V(a_4)=0$, from \eqref{CodKeq1}.2 that
 $T(a_4)=-4 \ac22 a_4$ and from \eqref{CodKeq1}.3 that $\ac23
 a_4=0$. Equation \eqref{CodKeq2}.1 becomes $W(a_4)=0$, equation
 \eqref{CodKeq3}.2 $T(a_6)=-\ac22 a_6$ and \eqref{CodKeq3}.3 $(\ac23+3
 \bc13)a_6=0$. Finally equation \eqref{CodKeq4}.2 resp. 3 gives
 $W(a_6)=-3 \bc23 a_6$ and $V(a_6)= 3 \bc33 a_6$.
 
 In case of $SO(2)$-symmetry ($a_4>0$ and $a_6=0$) it follows that
 $\ac23=0$ and thus the statement of the theorem.
 
 In case of $S_3$-symmetry ($a_4=0$ and $a_6>0$) it follows that
 $\ac23=-3\bc13$ and thus the statement of the theorem.
 
 In case of $\Z_3$-symmetry ($a_4>0$ and $a_6>0$) it follows that
 $\ac23=0$ and $\bc13=0$ and thus the statement of the theorem.
 
 In case that $a_6= \pm 2 a_4$ ($\neq 0$), we can choose $V, W$ such
 that $\mathbf{a_6 = 2a_4}$. Now equations \eqref{CodKeq8},
 \eqref{CodKeq1}.3 and \eqref{CodKeq3}.3 lead to $\ac23=\ac32$,
 $\ac13=-2 \ac23$ and $\bc13=0$. A combination of \eqref{CodKeq1}.2
 and \eqref{CodKeq2}.3 gives $\ac12=(\ac22-\ac33)$, and then by
 equations \eqref{CodKeq3}.2, \eqref{CodKeq1}.1 and \eqref{CodKeq1}.2
 that $\ac33=-\ac22$. Thus $T(a_4)=0$ by \eqref{CodKeq1}.2,
 $V(a_4)=- 4 \ac22 a_4$ by \eqref{CodKeq1}.1 and $W(a_4)= 4 \ac22
 a_4$ by \eqref{CodKeq2}.1. Finally \eqref{CodKeq4}.2 and
 \eqref{CodKeq2}.1 resp. \eqref{CodKeq4}.3 and \eqref{CodKeq1}.1 imply
 that $\bc23=-\ac23$ resp. $\bc33=-\ac22$.
 \end{proof}
 
An evaluation of the Gauss equations \eqref{gaussLC} with the help of
the CAS Mathematica leads to the following :

\begin{lem}\label{lem:GaussLC}
Let $\M$ be an affine hypersphere in $\mathbb R^4$ which admits a
pointwise $\SO(2)$-, $S_3$- or $\Z_3$-symmetry and $\{T,V,W\}$ the
corresponding ONB. Then
\begin{align} 
T(\ac22)&= -\ac22^2 +\ac23^2 +H-3 a_4^2,\label{Gauss1.1}\\
T(\ac23)&=-2\ac22\ac23,\label{Gauss1.2}\\ W(\ac22) + V(\ac23)
&=0,\label{Gauss1.3}\\ W(\ac23) - V(\ac22) &=0,\label{Gauss1.4}\\
V(\bc13) - T(\bc23)&= \ac22\bc23 + (\ac23 +
\bc13)\bc33,\label{Gauss1.5}\\ T(\bc33) - W(\bc13)&=
(\ac23+\bc13)\bc23 - \ac22\bc33 ,\label{Gauss1.6}\\ V(\bc33) -
W(\bc23)&= -\ac22^2-\ac23^2 +2 \ac23\bc13 +\bc23^2 +\bc33^2 +H +a_4^2
+2 a_6^2,\label{Gauss1.7}
\end{align}
If the symmetry group is $\Z_3$, then $a_6^2\neq 4 a_4^2$.
\end{lem}

\begin{proof}
The equations relate to eq11--eq13 and eq16 in the Mathematica
notebook. If $a_6^2= 4 a_4^2 (\neq 0)$, then we obtain by equations
eq11.1 and eq12.3 resp. eq15.3 and eq12.3 that $2 V(\ac22)=-4 \ac22^2
-H+ 3 a_4^2$ resp. $2 W(\ac23)=4 \ac23^2+H-3 a_4^2$, thus $
V(\ac22)-W(\ac23)= -2 \ac22^2 -2 \ac23^2 -H +3 a_4^2$. This gives a
contradiction to eq13.3, namely $ V(\ac22)-W(\ac23)= -2 \ac22^2 -2
\ac23^2 -H -9 a_4^2$.
\end{proof}

\subsection{Pointwise $\Z_3$- or $\SO(2)$-symmetry} \label{subsec:type2,4}
As the vector field $T$ is globally defined, we can define the
distributions $L_1=\Span\{T\}$ and $L_2=\Span\{V,W\}$. In the
following we will investigate these distributions. For the terminology
we refer to \cite{Noe96}.
\begin{lem}\label{L1}
The distribution $L_1$ is autoparallel (totally geodesic) with respect
to $\widehat\nabla$.
\end{lem}
\begin{proof} From $\widehat{\nabla}_{T} T = \ac12 V + \ac13 W=0$ 
(cp. Lemma~\ref{CodK}) the claim follows immediately.
\end{proof} 
\begin{lem}\label{L2}
  The distribution $L_2$ is spherical with mean curvature normal
  $U_2=\ac22 T$.
\end{lem}
\begin{proof} For $U_2=\ac22 T\in L_1=L_2^{\perp}$ we have
  $h(\widehat{\nabla}_{E_a} E_b, T)= h(E_a, E_b) h(U_2,T)$ for $E_a,
E_b\in \{V,W\}$, and $h(\widehat{\nabla}_{E_a} U_2, T)= h(E_a(\ac22) T
+ \ac22 \widehat{\nabla}_{E_a} T, T)=0$ (cp. Lemma \ref{lem:CodK} and
\eqref{Gauss1.3}, \eqref{Gauss1.4}.
\end{proof}
\begin{rem} $\ac22$ is independent of the
  choice of ONB $\{V,W\}$. It therefore is a globally defined
  function on $M$. 
\end{rem}
We introduce a coordinate function $t$ by $\pt:=T$. Using the
previous lemma, according to \cite{PR93}, we get:
\begin{lem}\label{warped} $(\M,h)$ admits a warped product structure
  $\M=I \times_{e^f}N^2$ with $f: I \to \mathbb R$
  satisfying
\begin{equation}\label{deff}
\frac{\partial f}{\partial t}=\ac22.
\end{equation}
\end{lem}
\begin{proof} Prop. 3 in \cite{PR93} gives the warped product structure with 
warping function $\lambda_2:I \to \mathbb R$. If we introduce $f=\ln
\lambda_2$, following the proof we see that $\ac22
T=U_2=-\grad(\ln\lambda_2)=-\grad f$.
\end{proof}
\begin{lem}\label{lem:curvN2} The curvature of $N^2$ is 
${}^NK(N^2)=e^{2f}(H+2 a_6^2 + a_4^2-\ac22^2)$.
\end{lem}
\begin{proof} From Prop. 2 in \cite{PR93} we get the following relation 
between the curvature tensor $\hat R$ of the warped product $\M$ and
the curvature tensor $\tilde R$ of the usual product of
pseudo-Riemannian manifolds ($X,Y,Z\in {\cal X}(M)$ resp. their
appropriate projections):
\begin{equation*}\begin{split}\hat R(X,Y)Z&= \tilde{R} (X,Y)Z \\
&+ h(Y,Z)(\widehat{\nabla}_X U_2 - h(X,U_2)U_2) - h(\widehat{\nabla}_X
U_2- h(X,U_2)U_2,Z)Y \\ &- h(X,Z)(\widehat{\nabla}_Y U_2 -
h(Y,U_2)U_2) + h(\widehat{\nabla}_Y U_2- h(Y,U_2)U_2,Z)X \\
&+h(U_2,U_2)(h(Y,Z)X-h(X,Z)Y)
\end{split}\end{equation*}
Now $\tilde{R}(X,Y)Z={}^N\hat R(X,Y)Z$ for all $X,Y,Z\in TN^2$ and
otherwise zero (cp. \cite{O'N83}, pg. 89, Corollary 58) and
$K(N^2)=K(V,W)= \frac{h(-\hat{R}(V,W)V,W)}{h(V,V)h(W,W)-h(V,W)^2}$
(cp. \cite{O'N83}, pg. 77, the curvature tensor has the opposite
sign). Since $h(X,Y)=e^{2f} {}^Nh(X,Y)$ for $X,Y\in TN^2$, it follows
that $${}^NK(N^2)=e^{2f} h(-{}^N\hat{R}(V,W)V,W).$$ Finally we obtain by
the Gauss equation~\eqref{gaussLC} the last
ingredient for the computation: $\hat R(V,W)V =-(H+2 a_6^2 + a_4^2)
W$ (cp. the Mathematica notebook).
\end{proof}
Summarized we have obtained the following structure equations
(cp. \eqref{strGauss}, \eqref{strWeingarten} and \eqref{defK}), where
$a_6=0$ in case of $\SO(2)$-symmetry resp. $\bc13=0$ in case of
$\Z_3$-symmetry:
\begin{alignat}{4}
&D_T T =& && -2a_4 T&- \xi, \label{D11}\\
&D_T V =& +a_4 V & - \bc13 W,&& \label{D12}\\
&D_T W =& +\bc13 V & + a_4 W,&& \label{D13}\\
&D_V T =&+(\ac22 +a_4) V, &&& \label{D21}\\
&D_W T =& &+(\ac22 + a_4)W, && \label{D31}\\
&D_V V =&+ a_6 V &-\bc23 W & +(\ac22 - a_4)T &+\xi, \label{D22}\\
&D_V W =& +\bc23 V &- a_6 W,&& \label{D23}\\
&D_W V =&&-(\bc33 +a_6) W, && \label{D32}\\
&D_W W =& +(\bc33- a_6) V & & +(\ac22 - a_4) T &+\xi, \label{D33}
\end{alignat}
\begin{equation}\label{Dxi}
D_{X} \xi= -H X,
\end{equation}
The Codazzi and Gauss equations (\eqref{CodK} and \eqref{gaussLC})
have the form (cp. Lem.~\ref{lem:CodK} and \ref{lem:GaussLC}):
\begin{align}
 T(a_4)&=-4\ac22 a_4, \quad 0=V(a_4) = W(a_4)\label{Da4}\\
 T(a_6)&=-\ac22 a_6,\quad V(a_6) = 3\bc33 a_6,\quad W(a_6)=-3 \bc23
 a_6,\label{Da6}\\ 
 T(\ac22)&= -\ac22^2 +H-3 a_4^2, \;
 V(\ac22)=0, \; W(\ac22) =0,\label{Da22}\\ 
 V(\bc13) - T(\bc23)&= \ac22\bc23 + \bc13\bc33,\label{Db1}\\ 
 T(\bc33) - W(\bc13)&=\bc13\bc23 - \ac22\bc33 ,\label{Db2}\\ 
 V(\bc33) - W(\bc23)&=
 -\ac22^2+\bc23^2 +\bc33^2 +H +a_4^2 +2 a_6^2,\label{Db3}
\end{align}
where $a_6=0$ in case of $\SO(2)$-symmetry resp. $\bc13=0$ in case of
$\Z_3$-symmetry.

Our first goal is to find out how $N^2$ is immersed in $\Rf$, i.~e. to
find an immersion independent of $t$. A look at the structure
equations \eqref{D11} - \eqref{Dxi} suggests to start with a linear
combination of $T$ and $\xi$.

We will solve the problem in two steps. First we look for a vector
field $X$ with $D_T X=\alpha X$ for some funtion $\alpha$: We define
$X:=A T +\xi$ for some function $A$ on $\M$. Then $D_T X=\alpha X$ iff
$\alpha=-A$ and $\pt A= -A^2 +2a_4 A+ H$, and $A:=\ac22- a_4$ solves
the latter differential equation. Next we want to multiply $X$ with
some function $\beta$ such that $D_T (\beta X)=0$: We define a positive
function $\beta$ on $\R$ as the solution of the differential equation:
\begin{equation}\label{dtbeta}
\tfrac{\partial}{\partial t} \beta = (\ac22- a_4)\beta 
\end{equation}
with initial condition $\beta(t_0)>0$. Then $D_T(\beta X)=0$ and by
\eqref{D21}, \eqref{Dxi} and \eqref{D31} we get (since
$\beta$, $\ac22$ and $a_4$ only depend on $t$):
\begin{align}
D_{T}(\beta((\ac22- a_4)T +\xi))&=0,\label{eq31}\\
D_{V}(\beta((\ac22- a_4)T +\xi))&=\beta(\ac22^2-a_4^2-H)V ,\label{eq32}\\ 
D_{W}(\beta((\ac22- a_4)T +\xi))&=\beta(\ac22^2-a_4^2-H)W.\label{eq33} 
\end{align}
To obtain an immersion we need that $\nu:=\ac22^2-a_4^2-H$ vanishes
nowhere, but we only get:
\begin{lem}\label{nu}
  The function $\nu=\ac22^2-a_4^2-H$ is globally defined,
  $\pt(e^{2f} \nu)=0$ and $\nu$ vanishes identically or nowhere on $\R$.
\end{lem}
\begin{proof} Since $0=\pt {}^NK(N^2) = \pt(e^{2f}(2a_6^2-\nu))$ 
  (Lem.~\ref{lem:curvN2}) and $\pt(e^{2f}2 a_6^2)=0$ (cp.~\eqref{Da6} and
  \eqref{deff}), we get that $\pt(e^{2f} \nu)=0$. Thus $\pt\nu=-2 (\pt
  f)\nu= -2\ac22\nu$.
\end{proof}

\subsubsection{The first case: $\nu \neq 0$ on $\M$}
\label{sec:case1}

We may, by translating $f$, i.e. by replacing $N^2$ with a homothetic
copy of itself, assume that $e^{2f} \nu =\eps_1$, where $\eps_1 =\pm 1$.

\begin{lem}\label{defphi}
$\varPhi:=\beta ((\ac22- a_4)T +\xi)\colon
M^3 \to \R^4$ induces a proper affine sphere structure, say
$\tilde{\phi}$, mapping $N^2$ into a 3-dimensional linear subspace
of $\R^4$. $\tilde{\phi}$ is part of a quadric iff $a_6 =0$.
\end{lem}
\begin{proof} 
By \eqref{eq32} and \eqref{eq33} we have $\varPhi_*(E_a)= \beta
\nu E_a$ for $E_a\in \{V,W\}$. A further differentiation, using \eqref{D22}
($\beta$ and $\nu$ only depend on $t$), gives:
\begin{align*}
D_{V} \varPhi_*(V)& = \beta \nu D_{V} V\\
&= \beta \nu ((\ac22- a_4)T +a_6 V - \bc23 W +\xi)\\
&=a_6\varPhi_*(V)-\bc23 \varPhi_*(W) +\nu \varPhi\\
&=a_6\varPhi_*(V)-\bc23 \varPhi_*(W) +\eps_1 e^{-2f}\varPhi.
\end{align*}
Similarly, we obtain the other derivatives, using \eqref{D23} -
\eqref{D33}, thus:
\begin{alignat}{3}
D_{V} \varPhi_*(V)&= & a_6\varPhi_*(V) &-\bc23 \varPhi_*(W)&+
e^{-2f}\eps_1 \varPhi, \label{Dphivv}\\ D_{V} \varPhi_*(W)&=& \bc23
\varPhi_*(V) &-a_6 \varPhi_*(W),& \label{Dphivw}\\ D_{W} \varPhi_*(V)&=&
&-(\bc33+ a_6)\varPhi_*(W),& \label{Dphiwv}\\ D_{W} \varPhi_*(W)&=&
(\bc33-a_6)\varPhi_*(V) &&+e^{-2f}\eps_1 \varPhi, \label{Dphiww}\\ D_{E_a}
\varPhi &=& \beta e^{-2f}\eps_1 E_a.\label{Dphiea} & &
\end{alignat}
The foliation at $f=f_0$ gives an immersion of $N^2$ to $M^3$, say
$\pi_{f_0}$. Therefore, we can define an immersion of $N^2$ to $\R^4$
by $\tilde{\phi}:=\varPhi\circ\pi_{f_0}$, whose structure equations
are exactly the equations above when $f=f_0$. Hence, we know that
$\tilde{\phi}$ maps $N^2$ into
$\Span\{\varPhi_*(V),\varPhi_*(W),\varPhi\}$, an affine hyperplane of
$\R^4$ and $\pt\varPhi=0$ implies $\varPhi(t,v,w)=\tilde{\phi}(v,w)$.

We can read off the coefficients of the difference tensor
$K^{\tilde{\phi}}$ of $\tilde{\phi}$ (cf. \eqref{strGauss} and
\eqref{defK}): $K^{\tilde{\phi}}(\tilde{V},\tilde{V})=a_6 \tilde{V}$,
$K^{\tilde{\phi}}(\tilde{V},\tilde{W})=-a_6
\tilde{W}$,$K^{\tilde{\phi}}(\tilde{W},\tilde{W})=- a_6 \tilde{V}$, and
see that $\trace (K^{\tilde{\phi}})_X$ vanishes. The affine metric
introduced by this immersion corresponds with the metric on $N^2$.
Thus $\eps_1 \tilde{\phi}$ is the affine normal of $\tilde{\phi}$ and
$\tilde{\phi}$ is a proper affine sphere with mean curvature
$\eps_1$.  Finally the vanishing of the difference tensor
characterizes quadrics.
\end{proof}
Our next goal is to find another linear combination of $T$ and $\xi$,
this time only depending on $t$. (Then we can express $T$ in terms
of $\phi$ and some function of $t$.)
\begin{lem}\label{defdelta}
  Define $\delta := H T +(\ac22+a_4) \xi$. Then there exist a constant
  vector $C \in \R^4$ and a function $a(t)$ such that
$$ \delta(t)= a(t) C.$$
\end{lem}
\begin{proof} Using \eqref{D21} resp. \eqref{D31} and
  \eqref{Dxi} we obtain that $D_{V}\delta = 0=D_{W} \delta$. Hence
  $\delta$ depends only on the variable $t$. Moreover, we get by
  \eqref{D11}, \eqref{Da22},\eqref{Da4} and \eqref{Dxi} that
\begin{align*}
  \pt\delta&=D_{T} (H T+(\ac22+a_4)\xi)\\ 
  &=H(-2a_4 T-\xi) + (-\ac22^2
  +H -3 a_4^2-4 \ac22 a_4)\xi -(\ac22 + a_4) H T\\ 
  &=-(3 a_4+\ac22)(H
  T +(\ac22+a_4)\xi)\\ 
  &=-(3 a_4+\ac22) \delta.
\end{align*}
This implies that there exists a constant vector $C$ in $\Rf$ and a
function $a(t)$ such that $\delta(t)=a(t)C$.
\end{proof}
Notice that for an improper affine hypersphere ($H=0$) $\xi$ is
constant and parallel to $C$. Combining $\tilde{\phi}$ and $\delta$ we
obtain for $T$ (cp. Lem.~\ref{defphi} and \ref{defdelta}) that
\begin{equation}\label{T}
T(t,v,w)= -\frac{a}{\nu}C +\frac{1}{\beta\nu}(\ac22+a_4)\tilde{\phi}(v,w).
\end{equation}

In the following we will use for the partial derivatives the
abbreviation $\hs_x:= \frac{\partial}{\partial x}\hs $, $x=t,v,w$.
\begin{lem}\label{partialF} 
\begin{align*}
&\hs_t = -\frac{a}{\nu}C +\pt(\frac{1}{\beta \nu})\tilde{\phi},\\
&\hs_v = \frac{1}{\beta\nu} \tilde{\phi}_v,\\
&\hs_w = \frac{1}{\beta\nu} \tilde{\phi}_w.
\end{align*}
\end{lem}
\begin{proof} As by \eqref{dtbeta} and Lem.~\ref{nu} $\pt
  \frac{1}{\beta\nu}= \frac{1}{\beta\nu}(\ac22+a_4)$, we obtain the
  equation for $\hs_t =T$ by \eqref{T}. The other equations follow
  from \eqref{eq32} and \eqref{eq33}.
\end{proof}
It follows by the uniqueness theorem of first order differential
equations and applying a translation that we can write
$$\hs(t,v,w)= \tilde{a}(t) C +\frac{1}{\beta\nu}(t)
\tilde{\phi}(v,w)$$ for a suitable function $\tilde{a}$ depending only
on the variable $t$. Since $C$ is transversal to the image of
$\tilde{\phi}$ (cp. Lem.~\ref{defphi} and \ref{defdelta},
$\nu\not\equiv 0$), we obtain that after applying an equiaffine
transformation we can write: $\hs(t,v,w) =(\gamma_1(t), \gamma_2(t)
\phi(v,w))$, in which $\tilde{\phi}(v,w)=(0,\phi(v,w))$.  Thus we have
proven the following:

\begin{thm}\label{thm:ClassC1} Let $\M$ be an indefinite affine hypersphere 
of $\,\Rf$ which admits a pointwise $\mathbb
 Z_3$- or
  $SO(2)$-symmetry. Let $\ac22^2-a_4^2 \neq H$ for some $p\in
  \M$. Then
 $\M$ is affine equivalent to
  $$\hs:I\times N^2\to \Rf:(t,v,w)\mapsto (\gamma_1(t), \gamma_2(t)
  \phi(v,w)),$$ where $\phi: N^2 \to \mathbb R^3$ is a (positive
  definite) elliptic or hyperbolic affine sphere and $\gamma:I\to
  \mathbb R^2$ is a curve.\hfill \newline Moreover, if $\M$
  admits a pointwise $SO(2)$-symmetry then $N^2$ is either an
  ellipsoid or a two-sheeted hyperboloid.
\end{thm}

We want to investigate the conditions imposed on the curve $\ga$. For
this we compute the derivatives of $\hs$:
\begin{alignat}{4}\label{DF}
\hs_t &=(\ga_1',\ga_2' \phi ),\quad & \hs_v &=(0,\ga_2\phi_v ),\quad & \hs_w
&=(0,\ga_2\phi_w ),\nonumber\\ \hs_{tt}&=(\ga_1'',\ga_2'' \phi),\quad &
\hs_{tv}&=(0,\ga_2'\phi_{v}),\quad & \hs_{tw}&=(0,\ga_2'\phi_w),\\
\hs_{vv}&=(0,\ga_2\phi_{vv}),\quad & \hs_{vw}&=(0,\ga_2\phi_{vw}),\quad &
\hs_{ww}&=(0,\ga_2'\phi_{ww}) \nonumber.
\end{alignat}
Furthermore we have to distinguish if $\M$ is proper ($H=\pm 1$) or
improper ($H=0$).

First we consider the case that $\M$ is proper, i.~e. $\xi=-H\hs$. An
easy computation shows that the condition that $\xi$ is a transversal
vector field, namely $ 0\neq \det(\hs_t ,\hs_v, \hs_w, \xi)=-\ga_2^2
(\ga_1\ga_2' - \ga_1' \ga_2) \det(\phi_v, \phi_w, \phi)$, is
equivalent to $\ga_2\neq 0$ and $\ga_1\ga_2' - \ga_1' \ga_2 \neq
0$. To check the condition that $\xi$ is the Blaschke normal
(cp. \eqref{BlaschkeNormal}), we need to compute the Blaschke metric
$h$, using \eqref{strGauss}), \eqref{DF},
\eqref{Dphivv}--\eqref{Dphiww} and the notation $r,s\in\{v,w\}$ and
$g$ for the Blaschke metric of $\phi$:
\begin{align*}
\hs_{tt}&= \ldots \hs_t +\frac{\ga_1'\ga_2'' - \ga_1''
\ga_2'}{H(\ga_1\ga_2' - \ga_1' \ga_2)} \xi,\\ \hs_{tr}&= \text{tang.} \\
\hs_{rs}&= \text{tang.}  -\frac{\ga_1' \ga_2}{H(\ga_1\ga_2' - \ga_1'
\ga_2)}\eps_1 g(\frac{\p}{\p r}, \frac{\p}{\p s})\xi.
\end{align*}
We obtain that $\det h = h_{tt}
(h_{vv}h_{ww}-h_{vw}^2)=\frac{\ga_1'\ga_2'' - \ga_1'' \ga_2'}{H^3
(\ga_1\ga_2' - \ga_1' \ga_2)^3}(\ga_1')^2 \ga_2^2 \det g$. Thus
\eqref{BlaschkeNormal} is equivalent to $\ga_2^4 (\ga_1\ga_2' - \ga_1'
\ga_2)^2 \det(\phi_v,\phi_w,\phi)^2=|\frac{\ga_1'\ga_2'' - \ga_1''
\ga_2'}{(\ga_1\ga_2' - \ga_1' \ga_2)^3} (\ga_1')^2 \ga_2^2 \det
g|$. Since $\phi$ is a definite proper affine sphere with normal
$-\eps_1\phi$, we can again use \eqref{BlaschkeNormal} to obtain
\[\xi=-H\hs \Longleftrightarrow \ga_2^2|\ga_1\ga_2' - \ga_1' \ga_2|^5=
 |\ga_1'\ga_2'' - \ga_1'' \ga_2'| (\ga_1')^2 \neq 0.\] From the
computations above ($g$ is positive definite) also it follows that $\hs$
is indefinite iff either $H\sign(\ga_1\ga_2' - \ga_1' \ga_2) =
\sign(\ga_1'\ga_2'' - \ga_1'' \ga_2')= \sign(\ga_1'\ga_2
\eps_1)$ or $-H\sign(\ga_1\ga_2' - \ga_1' \ga_2) =
\sign(\ga_1'\ga_2'' - \ga_1'' \ga_2')= \sign(\ga_1'\ga_2
\eps_1)$.

Next we consider the case that $\M$ is improper, i.~e. $\xi$ is
constant. By Lem.~\ref{defdelta} $\xi$ is parallel to $C$ and thus
transversal to $\phi$. Hence we can apply an affine transformation to
obtain $\xi=(1,0,0,0)$. An easy computation shows that the condition
that $\xi$ is a transversal vector field, namely $0\neq \det(\hs_t
,\hs_v, \hs_w, \xi)=-\ga_2^2 \ga_2' \det(\phi_v, \phi_w, \phi)$, is
equivalent to $\ga_2\neq 0$ and $\ga_2' \neq 0$. To check the
condition that $\xi$ is the Blaschke normal
(cp. \eqref{BlaschkeNormal}) we need to compute the Blaschke metric
$h$, using \eqref{strGauss}), \eqref{DF},
\eqref{Dphivv}--\eqref{Dphiww} and the notation $r,s\in\{v,w\}$ and
$g$ for the Blaschke metric of $\phi$:
\begin{align*}
\hs_{tt}&= \ldots \hs_t -\frac{\ga_1'\ga_2'' - \ga_1'' \ga_2'}{\ga_2'}
\xi,\\ \hs_{tr}&= \text{tang.} \\ \hs_{rs}&= \text{tang.}  +\frac{\ga_1'
\ga_2}{\ga_2'}\eps_1 g(\frac{\p}{\p r}, \frac{\p}{\p s})\xi.
\end{align*}
We obtain that $\det h = h_{tt}
(h_{vv}h_{ww}-h_{vw}^2)=-\frac{\ga_1'\ga_2'' - \ga_1''
\ga_2'}{(\ga_2')^3}(\ga_1')^2 \ga_2^2 \det g$. Thus
\eqref{BlaschkeNormal} is equivalent to $\ga_2^4 (\ga_2')^2
\det(\phi_v,\phi_w,\phi)^2=|\frac{\ga_1'\ga_2'' - \ga_1''
\ga_2'}{(\ga_2')^3}(\ga_1')^2 \ga_2^2 \det g|$. Since $\phi$ is a
definite proper affine sphere with normal $-\eps_1\phi$, we can again
use \eqref{BlaschkeNormal} to obtain
\[\xi=(1,0,0,0) \Longleftrightarrow \ga_2^2|\ga_2'|^5=
|\ga_1'\ga_2'' - \ga_1'' \ga_2'| (\ga_1')^2 \neq 0 .\] From the
computations above also it follows that $\hs$ is indefinite iff either
$- \sign(\ga_2')= \sign(\ga_1'\ga_2'' - \ga_1'' \ga_2')= \sign(\ga_1'
\ga_2 \eps_1)$ or $ \sign(\ga_2')= \sign(\ga_1'\ga_2'' - \ga_1''
\ga_2')= \sign(\ga_1' \ga_2 \eps_1)$.

Now we are ready for the converse.
\begin{thm}\label{thm:ExC1} 
  Let $\phi:N^2 \to \mathbb R^3$ be a positive definite elliptic or
  hyperbolic affine sphere (with mean curvature $\eps_1=\pm 1$), and
  let $\gamma: I \to \mathbb R^2$ be a curve such that 
  $\hs(t,v,w)=(\gamma_1(t), \gamma_2(t) \phi(v,w ))$ defines a 3-dimensional indefinite
  affine hypersphere. Then $\hs(N^2\times I)$  admits a pointwise $\mathbb Z_3$- or
 $SO(2)$-symmetry.

  \begin{romanlist}
  \item If $\ga=(\ga_1,\ga_2)$ satisfies $\ga_2^2|\ga_1\ga_2' - \ga_1'
  \ga_2|^5= \sign(\ga_1' \ga_2 \eps_1)(\ga_1'\ga_2'' - \ga_1'' \ga_2')
  (\ga_1')^2\neq 0$, then $\hs$ defines a 3-dimensional indefinite
  proper affine hypersphere.
  \item If $\ga=(\ga_1,\ga_2)$ satisfies $\ga_2^2|\ga_2'|^5=
  \sign(\ga_1' \ga_2 \eps_1)(\ga_1'\ga_2'' - \ga_1'' \ga_2')
  (\ga_1')^2\neq 0$, then $\hs$ defines a 3-dimensional indefinite
  improper affine hypersphere.
  \end{romanlist}
\end{thm}
\begin{proof}
We already have shown that $\hs$ defines a 3-dimensional indefinite
proper resp. improper affine hypersphere. To prove the symmetry we
need to compute $K$. By assumption, $\phi$ is an affine sphere with
Blaschke normal $\xi^\phi =-\eps_1 \phi$. For the structure equations
\eqref{strGauss} we use the notation $\phi_{rs} = {}^\phi \Ga_{rs}^{u}
\phi_u - g_{rs} \eps_1 \phi$, $r,s,u \in \{v,w\}$. Furthermore we
introduce the notation $\al = \ga_1 \ga_2' - \ga_1' \ga_2$. Note that
$\al'=\ga_1 \ga_2'' - \ga_1'' \ga_2$. \newline (i) Using \eqref{DF},
we get the structure equations \eqref{strGauss} for $\hs$:
\begin{align*}
\hs_{tt}&= \frac{\al'}{\al} \hs_t +\frac{\ga_1'\ga_2'' - \ga_1''
\ga_2'}{H\al} \xi,\\ \hs_{tr}&= \frac{\ga_2'}{\ga_2} \hs_r, \\ \hs_{rs}&=
{}^\phi \Ga_{rs}^{u} \hs_u - g_{rs} \eps_1 \frac{\ga_1 \ga_2}{\al} \hs_t -
g_{rs} \eps_1 \frac{\ga_1' \ga_2}{H\al}\xi.
\end{align*}
We compute $K$ using \eqref{defC} and obtain:
\begin{align*}
(\nabla_{\hs_t} h)(\hs_r,\hs_s) &= ((\frac{\ga_1 \ga_2}{\al})'
\frac{\al}{\ga_1 \ga_2} -2 \frac{\ga_2'}{\ga_2}) h(\hs_r,\hs_s),\\
(\nabla_{\hs_r} h)(\hs_t,\hs_t) &= 0,
\end{align*}
implying that $K_{\hs_t}$ restricted to the space spanned by $\hs_v$ and
$\hs_w$ is a multiple of the identity. Taking $T$ in direction of $\hs_t$,
we see that $\hs_v$ and $\hs_w$ are orthogonal to $T$. Thus we can
construct an ONB $\{T,V,W\}$ with $V,W$ spanning $\Span\{\hs_v,\hs_w\}$
such that $a_1 = 2 a_4$, $a_2=a_3=a_5=0$. By the considerations in
Sec.~\ref{sec:max} we see that $\hs$ admits a pointwise $\mathbb Z_3$- or
 $SO(2)$-symmetry.
\newline (ii) The proof runs completely analog.
\end{proof}

\subsubsection{The second case: $\nu \equiv 0$ and $H\neq 0$ on $\M$}
\label{sec:case2}

Next, we consider the case that $H =\ac22^2 - a_4^2$ and $H\neq 0$ on
$\M$. It follows that $\ac22\neq \pm a_4$ on $\M$. 

We already have seen that $\M$ admits a warped product structure. The
map $\varPhi$ we have constructed in Lemma~\ref{defphi} will not define
an immersion (cp. \eqref{eq32} and \eqref{eq33}). Anyhow, for a fixed
point $t_0$, we get from \eqref{D22} - \eqref{D33}, \eqref{eq32} and
\eqref{eq33}, using the notation $\tilde{\xi}=(\ac22-a_4) T + \xi$:
\begin{align*}
D_V V &=a_6 V - \bc23 W +\tilde{\xi},\\
D_V W &=\bc23 V - a_6 W,\\
D_W V &=-(\bc33 + a_6) W,\\
D_W W &=(\bc33 - a_6) V +\tilde{\xi},\\
D_{E_a}\tilde{\xi}&=0, \qquad E_a\in \{V,W\}.
\end{align*}
Thus, if $v$ and $w$ are local coordinates which span the second
distribution $L_2$, then we can interprete $\hs(t_0,v,w)$ as a
positive definite improper affine sphere in a $3$-dimensional
linear subspace.

Moreover, we see that this improper affine sphere is a paraboloid
provided that $a_6(t_0, v,w)$ vanishes identically. From the
differential equations \eqref{Da6} determining $a_6$, we see that this
is the case exactly when $a_6$ vanishes identically, i.e.  when $\M$
admits a pointwise $SO(2)$-symmetry.

After applying a translation and a change of coordinates, we may
assume that
\begin{equation*}
\hs(t_0,v,w)=(v,w,f(v,w),0),
\end{equation*}
with affine normal $\tilde{\xi}(t_0,v,w)=(0,0,1,0)$. To obtain $T$
at $t_0$, we consider \eqref{D21} and \eqref{D31} and get that
\begin{equation*}
D_{E_a}(T-(\ac22+a_4) \hs) = 0,\qquad E_a, E_b\in \{V,W\}.
\end{equation*}
Evaluating at $t=t_0$, this means that there exists a constant vector
$C$, transversal to $\Span\{V,W,\xi\}$, such that
$T(t_0,v,w)=(\ac22+a_4)(t_0) \hs(t_0,v,w) +C$. Since $\ac22+a_4\neq 0$
everywhere, we can write:
\begin{equation}\label{Tt0}
  T(t_0,v,w)=\alpha_1 (v,w, f(v,w),\alpha_2),
\end{equation}
where $\alpha_1, \alpha_2\neq 0$ and we applied an equiaffine
transformation so that $C=(0,0,0,\alpha_1\alpha_2)$. To obtain
information about $D_T T$ we have that $D_T T= -2 a_4 T -\xi$
(cp. \eqref{D11}) and $\xi=\tilde{\xi} - (\ac22- a_4)T$ by the
definition of $\tilde{\xi}$.  Also we know that
$\tilde{\xi}(t_0,v,w)=(0,0,1,0)$ and by \eqref{eq31} - \eqref{eq33}
that $D_{X}(\beta\tilde{\xi})=0$, $X\in {\cal X}(M)$. Taking suitable
initial conditions for the function $\beta$ ($\beta(t_0)=1$), we get
that $\beta\tilde{\xi}=(0,0,1,0)$ and finally the following vector
valued differential equation:
\begin{equation*}
D_T T= (\ac22 -3 a_4) T -\frac1{\beta} (0,0,1,0).
\end{equation*}
Solving this differential equation, taking into account the initial
conditions \eqref{Tt0} at $t=t_0$, we get that there exist functions
$\delta_1$ and $\delta_2$ depending only on $t$ such that
\begin{equation*}
  T(t,u,v)= (\delta_1(t) v,\delta_1(t) w, \delta_1(t) (f(v,w)
  +\delta_2(t)), \alpha_2 \delta_1(t)), 
\end{equation*}
where $\delta_1(t_0)=\alpha_1$, $\delta_2(t_0)=0$, $\delta_1'(t)
=(\ac22 -3 a_4) \delta_1(t)$ and $\delta_2'(t) =\delta_1^{-1}(t)
\beta^{-1}(t)$.  As $T(t,v,w) =\tfrac{\partial \hs}{\partial
  t}(t,v,w)$ and $\hs(t_0,v,w) =(v,w,f(v,w),0)$ it follows by
integration that
$$\hs(t,v,w)= (\gamma_1(t) v, \gamma_1(t) w, \gamma_1(t) f(v,w)
+\gamma_2(t) , \alpha_2 (\gamma_1(t)-1)),$$
where $\gamma_1'(t)
=\delta_1(t)$, $\gamma_1(t_0)=1$, $\gamma_2(t_0)=0$ and $\gamma_2'(t)
=\delta_1(t)\delta_2(t)$.  After applying an affine transformation
we have shown:

\begin{thm} \label{thm:ClassC2} Let $\M$ be an indefinite proper affine 
  hypersphere of $\,\Rf$ which admits a pointwise $\mathbb Z_3$- or
 $SO(2)$-symmetry. Let 
  $H= \ac22^2 -a_4^2(\neq 0)$ on $M^3$. Then $\M$ is affine equivalent with
  $$\hs:I\times N^2\to \Rf:(t,v,w)\mapsto (\gamma_1(t) v, \gamma_1(t)
  w, \gamma_1(t) f(v,w) +\gamma_2(t),\gamma_1(t)),$$ where $\psi: N^2
  \to \mathbb R^3:(v,w) \mapsto (v,w,f(v,w))$ is a positive definite
  improper affine sphere with affine normal $(0,0,1)$ and $\gamma:I\to
  \mathbb R^2$ is a curve.\hfill \newline Moreover, if $\M$ admits a
  pointwise $SO(2)$-symmetry then $N^2$ is an elliptic paraboloid.
\end{thm}
We want to investigate the conditions imposed on the curve $\ga$. For
this we compute the derivatives of $\hs$:
\begin{align}\label{DF2}
\hs_t &=(\gamma_1' v,\gamma_1' w,\gamma_1' f(v,w)
+\gamma_2',\gamma_1'),\nonumber\\
\hs_v &=(\gamma_1,0,\gamma_1 f_v,0),\quad
 \hs_w =(0,\gamma_1,\gamma_1 f_w,0),\nonumber\\ 
\hs_{tt}&=(\gamma_1'' v, \gamma_1''w ,\gamma_1'' f(v,w) +
\gamma_2'',\gamma_1''), \\
\hs_{tv}&=\tfrac{\gamma_1'}{\gamma_1}
\hs_v,\quad  \hs_{tw}=\tfrac{\gamma_1'}{\gamma_1} \hs_w,\nonumber\\
\hs_{vv}&=(0,0,f_{vv}\gamma_1,0),\quad 
\hs_{vw}=(0,0,\gamma_1 f_{vw},0),\quad  \hs_{ww}=(0,0,\gamma_1 f_{ww},0)
\nonumber.
\end{align}

$\M$ is a proper hypersphere, i.~e. $\xi=-H\hs$. An easy computation
shows that the condition that $\xi$ is a transversal vector field,
namely $0\neq \det(\hs_t ,\hs_v, \hs_w, \xi)=-H \ga_1^2 (\ga_1\ga_2' -
\ga_1'\ga_2)$, is equivalent to $\ga_1\neq 0$ and $\ga_1\ga_2' -
\ga_1'\ga_2\neq 0$. Since $(0,0,1,0)= \frac{\ga_1}{\ga_1\ga_2' -
\ga_1'\ga_2} \hs_t - \frac{\ga_1'}{\ga_1\ga_2' - \ga_1'\ga_2} \hs$, we
have the following structure equations:
\begin{align}\label{streqC2}
\hs_{tt}&=(\tfrac{\ga_1''}{\ga_1'} + \tfrac{\ga_1'\ga_2'' -
\ga_1''\ga_2'}{\ga_1'}\tfrac{\ga_1}{\ga_1\ga_2' - \ga_1'\ga_2}) \hs_t +
\tfrac{\ga_1'\ga_2'' - \ga_1''\ga_2'}{\ga_1\ga_2' -
\ga_1'\ga_2}\tfrac1H \xi, \\ 
\hs_{tr}&=\tfrac{\gamma_1'}{\gamma_1} \hs_r,\nonumber\\
\hs_{rs}&=\tfrac{\ga_1^2}{\ga_1\ga_2' - \ga_1'\ga_2} f_{rs} \hs_t +
\tfrac{\ga_1\ga_1'}{\ga_1\ga_2' - \ga_1'\ga_2} f_{rs} \tfrac1H \xi
\nonumber.
\end{align}

We obtain that $\det h = h_{tt} (h_{vv}h_{ww}-h_{vw}^2)
=\frac{\ga_1'\ga_2'' - \ga_1'' \ga_2'}{H^3(\ga_1\ga_2' -
\ga_1'\ga_2)^3}\ga_1^2(\ga_1')^2 (f_{vv} f_{ww}-f_{vw}^2)$.  Since
$\psi$ is a positive definite improper affine sphere with affine
normal $(0,0,1)$, we get by \eqref{BlaschkeNormal} that $f_{vv}
f_{ww}-f_{vw}^2 =1$. Now \eqref{BlaschkeNormal} (for $\xi$) is
equivalent to $\ga_1^4 (\ga_1\ga_2' -\ga_1'\ga_2)^2 =
|\frac{\ga_1'\ga_2'' - \ga_1'' \ga_2'}{(\ga_1\ga_2'
-\ga_1'\ga_2)^3}|\ga_1^2(\ga_1')^2 $. It follows that

\[\xi=-H\hs \Longleftrightarrow \ga_1^2|\ga_1\ga_2' -\ga_1'\ga_2|^5=
|\ga_1'\ga_2'' - \ga_1'' \ga_2'| (\ga_1')^2 \neq 0 .\] From the
computations above also it follows that $\hs$ is indefinite iff either
$ \sign(\ga_1'\ga_2'' - \ga_1'' \ga_2')= \sign(H(\ga_1\ga_2'
-\ga_1'\ga_2))= -\sign(\ga_1 \ga_1')$ or $ \sign(\ga_1'\ga_2'' -
\ga_1'' \ga_2')= -\sign(H(\ga_1\ga_2' -\ga_1'\ga_2))= -\sign(\ga_1
\ga_1')$.

Now we can formulate the converse theorem:
\begin{thm}\label{thm:ExC2} 
  Let $\psi: N^2 \to \mathbb R^3:(v,w) \mapsto (v,w,f(v,w))$ be a
  positive definite improper affine sphere with affine normal
  $(0,0,1)$, and let $\gamma: I \to \mathbb R^2$ be a curve such that
  $\hs(t,v,w)=(\gamma_1(t) v, \gamma_1(t) w, \gamma_1(t) f(v,w)
  +\gamma_2(t),\gamma_1(t))$ defines a 3-dimensional indefinite proper affine
  hypersphere. Then $\hs(N^2\times I)$ admits a pointwise $\mathbb
  Z_3$- or $SO(2)$-symmetry. \\ Moreover, if $\ga=(\ga_1,\ga_2)$ satisfies
  $\ga_1^2|\ga_1\ga_2' -\ga_1'\ga_2|^5=
  -\sign(\ga_1\ga_1')(\ga_1'\ga_2'' - \ga_1'' \ga_2') (\ga_1')^2 \neq
  0$, then $\hs$  defines a 3-dimensional indefinite proper affine
  hypersphere.
\end{thm}
\begin{proof}
We already have shown that $\hs$ defines a 3-dimensional indefinite
proper affine hypersphere with affine normal $\xi= -H\hs$. To prove the
symmetry we need to compute $K$. We get the induced connection and the
affine metric from the structure equations \eqref{streqC2}. We
compute $K$ using \eqref{defC} and obtain:
\begin{align*}
(\nabla_{\hs_t} h)(\hs_r,\hs_s) &= (\pt\ln(\frac{\ga_1
\ga_1'}{\ga_1\ga_2' -\ga_1'\ga_2})-2 \frac{\ga_1'}{\ga_1})
h(\hs_r,\hs_s),\\ (\nabla_{\hs_r} h)(\hs_t,\hs_t) &= 0,
\end{align*}
implying that $K_{\hs_t}$ restricted to the space spanned by $\hs_v$ and
$\hs_w$ is a multiple of the identity. Taking $T$ in direction of $\hs_t$,
we see that $\hs_v$ and $\hs_w$ are orthogonal to $T$. Thus we can
construct an ONB $\{T,V,W\}$ with $V,W$ spanning $\Span\{\hs_v,\hs_w\}$
such that $a_1 = 2 a_4$, $a_2=a_3=a_5=0$. By the considerations in
Sec.~\ref{sec:max} we see that $\hs$ admits a pointwise
$\mathbb Z_3$- or $SO(2)$-symmetry.
\end{proof} 

\subsubsection{The third case: $\nu \equiv 0$ and $H=0$ on $\M$}
\label{sec:case3} 

 The final cases now are that $\nu \equiv 0$ and $H=0$ on the
whole of $M^3$ and hence $\ac22=\pm a_4$. 

First we consider the case that $\ac22= a_4 =:a>0$. Again we use that
$M^3$ admits a warped product structure and we fix a parameter
$t_0$. At the point $t_0$, we have by \eqref{D22}--\eqref{Dxi}:
\begin{alignat*}{3}
&D_V V =&+ a_6 V &-\bc23 W   &+\xi, \\
&D_V W =& +\bc23 V &- a_6 W,&\\
&D_W V =&&-(\bc33 +a_6) W, &\\
&D_W W =& +(\bc33- a_6) V & & +\xi,\\
&D_{X} \xi=& 0.
\end{alignat*}

Thus, if $v$ and $w$ are local coordinates which span the second
distribution $L_2$, then we can interprete $\hs(t_0,v,w)$ as a
positive definite improper affine sphere in a $3$-dimensional
linear subspace.

Moreover, we see that this improper affine sphere is a paraboloid
provided that $a_6(t_0, v,w)$ vanishes identically. From the
differential equations \eqref{Da6} determining $a_6$, we see that this
is the case exactly when $a_6$ vanishes identically, i.e.  when $\M$
admits a pointwise $SO(2)$-symmetry.

After applying an affine transformation and a change of coordinates,
we may assume that
\begin{equation}\label{initcond1}
\hs(t_0,v,w)=(v,w,f(v,w),0),
\end{equation}
with affine normal $\xi(t_0,v,w)=(0,0,1,0)$, actually
$$\xi(t,v,w)=(0,0,1,0)$$ ($\xi$ is constant on $\M$ by
assumption). Furthermore we obtain by \eqref{D21} and \eqref{D31},
that $D_U T= 2 a U$ for all $U\in L_2$. We define $\delta:= T- 2a
\hs$, which is transversal to $\Span\{V,W, \xi\}$. Since $a$ is
independent of $v$ and $w$ (cp. \eqref{Da4}), $D_U \delta =0$, and we
can assume that
\begin{equation}\label{initcond2}
T(t_0,v,w) - 2a(t_0) \hs(t_0,v,w)= (0,0,0,1).
\end{equation}
We can integrate \eqref{Da4} ($T(a)= -4 a^2$) and we take $a=
\frac1{4t}$, $t>0$. Thus \eqref{D11} becomes $D_T T=-\frac1{2t} T
-\xi$ and we obtain the following linear second order ordinary
differential equation:
\begin{equation}\label{DTT}  \pp{t}\hs + \frac1{2t} \pt \hs = -\xi.
\end{equation}
The general solution is $\hs(t,v,w)= -\frac{t^2}3 \xi + 2\sqrt{t}
A(v,w)+ B(v,w)$. The initial conditions \eqref{initcond1} and
\eqref{initcond2} imply that $A(v,w)=(\frac{v}{2\sqrt{t_0}},
\frac{w}{2\sqrt{t_0}}, \frac{f(v,w)}{2\sqrt{t_0}}+\frac23 t_0^{3/2},
\sqrt{t_0})$ and $B(v,w)=(0,0,-t_0^2,-2 t_0)$. Obviously we can
translate $B$ to zero. Furthermore we can translate the affine sphere
and apply an affine transformation to obtain
$A(v,w)=\frac{1}{2\sqrt{t_0}}(v, w, f(v,w), 1)$. After a change of
coordinates we get:
\begin{equation}\label{result3.1}
\hs(t,v,w)= (tv, tw, t f(v,w) - c t^4, t), \quad c,t>0.
\end{equation}
Next we consider the case that $-\ac22= a_4 =:a>0$. Again we use that
$M^3$ admits a warped product structure and we fix a parameter
$t_0$. A look at \eqref{D22}--\eqref{Dxi} suggests to define $\tilde{\xi}= -2a T+\xi$, then we get at the point $t_0$:
\begin{alignat*}{3}
&D_V V =&+ a_6 V &-\bc23 W   &+\tilde{\xi}, \\
&D_V W =& +\bc23 V &- a_6 W,&\\
&D_W V =&&-(\bc33 +a_6) W, &\\
&D_W W =& +(\bc33- a_6) V & & +\tilde{\xi},\\
&D_{U} \tilde{\xi}=& 0.
\end{alignat*}

Thus, if $v$ and $w$ are local coordinates which span the second
distribution $L_2$, then we can interprete $\hs(t_0,v,w)$ as a
positive definite improper affine sphere in a $3$-dimensional
linear subspace.

Moreover, we see that this improper affine sphere is a paraboloid
provided that $a_6(t_0, v,w)$ vanishes identically. From the
differential equations \eqref{Da6} determining $a_6$, we see that this
is the case exactly when $a_6$ vanishes identically, i.e.  when $\M$
admits a pointwise $SO(2)$-symmetry.

After applying an affine transformation and a change of coordinates,
we may assume that
\begin{equation}\label{init}
\hs(t_0,v,w)=(v,w,f(v,w),0),
\end{equation}
with affine normal 
\begin{equation}\label{initxi}
\tilde{\xi}(t_0,v,w)=(0,0,1,0). 
\end{equation}
We have considered $\tilde{\xi}$ before. We can solve
\eqref{dtbeta} ($\pt \beta= -2a \beta$) explicity by $\beta=c
\frac1{\sqrt{a}}$ (cp. \eqref{Da4}) and get by
\eqref{eq31}--\eqref{eq33} that $D_X (\frac1{\sqrt{a}} \tilde{\xi})
=0$. Thus $\frac1{\sqrt{a}}(-2a T+\xi)=:C$ for a constant vector $C$,
i.~e.  $T=-\frac1{2a}(\sqrt{a} C - \xi)$. Notice that by \eqref{Dxi}
$\xi$ is a constant vector, too. We can choose $a=\frac1{4|t|}$, $t<0$
(cp. \eqref{Da4}), and we obtain the ordinary differential equation:
\begin{equation}\label{DT}
\pt \hs= -\sqrt{|t|} C -2t \xi, \quad t<0.
\end{equation}
 The solution (after a
translation) with respect to the initial condition \eqref{init} is
$\hs(t,v,w)=\frac23 |t|^{\frac32} C- t^2 \xi + (v,w,f(v,w),0)$. Notice
that $C$ is a multiple of $\tilde{\xi}$ and hence by \eqref{initxi} a
constant multiple of $(0,0,1,0)$. Furthermore $\xi$ is transversal to 
the space spanned by
$\hs(t_0,v,w)$. So we get after an affine transformation and a change
of coordinates:
\begin{equation}\label{result3.2}
\hs(t,v,w)= (v,w, f(v,w) + c t^3, t^4 ),\quad c,t>0.
\end{equation}
Combining both results \eqref{result3.1} and \eqref{result3.2} we have:

\begin{thm} \label{thm:ClassC3} Let $\M$ be an indefinite improper affine 
  hypersphere of $\Rf$ which admits a pointwise $\mathbb
  Z_3$- or $SO(2)$-symmetry. Let $\ac22^2 =a_4^2$ on $M^3$. Then $\M$ is affine
  equivalent with either
  \begin{align*} &\hs:I\times N^2\to \Rf:(t,v,w)\mapsto (t v, t w,t f(v,w) -c
  t^4,t),\quad(\ac22 =a_4)\quad \text{or}\\
  &\hs:I\times N^2\to \Rf:(t,v,w)\mapsto (v, w, f(v,w) +c
  t^3,t^4),\quad(-\ac22 =a_4)
\end{align*}
 where $\psi: N^2 \to \mathbb R^3:(v,w) \mapsto (v,w,f(v,w))$ is a
  positive definite improper affine sphere with affine normal
  $(0,0,1)$ and $c,t\in \R^+$.\hfill \newline 
Moreover, if $\M$ admits a pointwise $SO(2)$-symmetry then $N^2$ is an
  elliptic paraboloid.
\end{thm}
The computations for the converse statement can be done completely
analogous to the previous cases, they even are simpler (the curve is
given parametrized). 
\begin{thm}\label{thm:ExC3} 
  Let $\psi: N^2 \to \mathbb R^3:(v,w) \mapsto (v,w,f(v,w))$ be a
  positive definite improper affine sphere with affine normal
  $(0,0,1)$. Define $\hs(t,v,w)= (t v, t w,t f(v,w) -c t^4,t)$
  or $\hs(t,v,w)= (v, w, f(v,w) +c t^3,t^4)$, where $c,t\in \R^+$.
  Then $\hs$ defines a 3-dimensional indefinite improper affine
  hypersphere, which admits a pointwise $\mathbb
  Z_3$- or $SO(2)$-symmetry.
\end{thm}

\subsection{Pointwise $S_3$-symmetry}
\label{subsec:type3}

In the beginning of the section (cp. Lem.~\ref{lem:KfT234}) we have
shown that in the case of $S_3$-symmetry, there exists a locally
defined orthonormal frame $\{T,V,W\}$ such that $K$ is given by:
\begin{align*} 
K(T,T)&= 0, & K(T,V)&=0, & K(T,W)&= 
0,\\
K(V,V)&= a_6 V, & K(V,W)&= -a_6 W, &
K(W,W)&= -a_6 V,
\end{align*}
(the vector field $T$ (up to sign) and the function $a_6$ are
invariantly defined) and the structure equations (\eqref{strGauss},
\eqref{strWeingarten}) are:
\begin{alignat}{4}
&D_T T =& && &- \xi, \label{D11t3}\\
&D_T V =&  & +\third \ac23 W,&& \label{D12t3}\\
&D_T W =& -\third \ac23 V, & && \label{D13t3}\\
&D_V T =&+\ac22 V &+\ac23 W,&& \label{D21t3}\\
&D_W T =&-\ac23 V &+\ac22 W, && \label{D31t3}\\
&D_V V =&+ a_6 V &-\bc23 W & +\ac22 T &+\xi, \label{D22t3}\\
&D_V W =& +\bc23 V &- a_6 W& +\ac23 T,& \label{D23t3}\\
&D_W V =&&-(\bc33 +a_6) W &-\ac23 T, & \label{D32t3}\\
&D_W W =& +(\bc33- a_6) V & & +\ac22 T &+\xi, \label{D33t3}
\end{alignat}
\begin{equation}\label{Dxit3}
D_{X} \xi= -H X,
\end{equation}
The Codazzi and Gauss equations (\eqref{CodK} and \eqref{gaussLC})
have the form (cp. Lem.~\ref{lem:CodK} and \ref{lem:GaussLC}):
\begin{align}
 T(\ln a_6)&=-\ac22,\quad V(\ln a_6) = 3\bc33,\quad W(\ln a_6)=-3
 \bc23, \label{Da6t3}\\ T(\ac22)&= -\ac22^2 +\ac23^2 +H, \nonumber\\
 V(\ac22)&=W(\ac23), \; W(\ac22) =-V(\ac23),\label{Da22t3}\\
 T(\ac23)&= -2 \ac22 \ac23, \; V(\ac23)=-3(T(\bc23) +\ac22\bc23
 +\tfrac23 \ac23\bc33), \nonumber\\ W(\ac23) &=3 (-T(\bc33)
 +\tfrac23 \ac23\bc23 -\ac22\bc33),\label{Da23t3}\\ V(\bc33) -
 W(\bc23)&= -\ac22^2 -\tfrac53 \ac23^2 +\bc23^2 +\bc33^2 +H +2
 a_6^2.\label{Db1t3}
\end{align}
An evaluation of the equation $[V,W]a_6 =(\widehat\nabla_V W
-\widehat\nabla_W V)a_6$ gives us in addition:
\begin{equation}\label{Db2t3}
V(\bc23)+W(\bc33) = \tfrac23 \ac22\ac23.
\end{equation}
Notice that by \eqref{Da6t3} and \eqref{Da22t3} the functions
$|\ac22|$ and $|\ac23|$ are invariantly defined.

Positive definite affine hyperspheres admitting an $S_3$-symmetry
turn out to satisfy Chen's equality \cite{Vr04}. They were studied in
several articles (\cite{SSVV97}, \cite{SV96} ,\cite{KSV01} ,\cite{KV99},
\cite{Vr04}) and now they are completely classified.  Our structure and
integrability conditions are kind of similar and luckily we can follow
the approach in \cite{KV99}, solving the system of equations
directly. We do expect that the approach of \cite{KSV01} (reduction to
a classification of 2-dimensional minimal surfaces in $S^5_3$ whose
ellipses are circles, done for elliptic affine hyperspheres) also is
generalizable to all cases.

From the structure equations \eqref{D11t3} and \eqref{Dxit3} we get
that $D_T D_T \hs= H \hs$ if $H= \pm 1$, resp. $D_T D_T \hs=
\text{const.}$, if $H=0$. This implies for the integral curves
$\gamma$ of $T$ that $\gamma''=H \gamma$, thus (cp. \cite{KSV01},
Cor. 1):
\begin{rem}\label{rem:ruled}
$\hs$ is locally ruled by arcs of ellipses resp. hyperbolas if it is proper 
and $H=-1$ resp. $H=1$, and by straight lines if it is improper ($H=0$).
\end{rem}
As a consequence of the ruling we can integrate the integrability
conditions in $T$-direction explicitly. We introduce a coordinate
function $t$ by $\pt:=T$.
\begin{lem}\label{lem:Int1}
For any function $t$ with $Tt=1$ there is a dense open subset
$\tilde{M} \subset \M$ and there are functions $h,k,l\colon
\tilde{M} \to \R$ with $Th= Tk= Tl= 0$ such that on $\tilde{M}$ we
have if 
\begin{description} 
\item[$\mathbf{H=-1}$] then either $\ac22=0$, $\ac23=\pm 1$, 
$a_6=\efct^k$ or 
\begin{align*} \ac22=&-\frac{\sin(t+l)\cos(t+l)}{\cos^2(t+l) + \sinh^2(h)}, 
\ac23=-\frac{\sinh(h)\cosh(h)}{\cos^2(t+l) + \sinh^2(h)},\\
a_6=&\frac{\efct^k}{\sqrt{\cos^2(t+l) +\sinh^2(h)}}.
\end{align*}
\item[$\mathbf{H=1}$] then either $\ac22=\pm1$, $\ac23=0$, 
$a_6=\efct^{k\mp t}$ or
\begin{align*} \ac22=&\frac{\sinh(t+l)\cosh(t+l)}{\cos^2(h) + \sinh^2(t+l)},
\quad \ac23= -\frac{\sin(h) \cos(h)}{\cos^2(h) + \sinh^2(t+l)},\\ 
a_6=&\frac{\efct^k}{\sqrt{\cos^2(h) +\sinh^2(t+l)}}.
\end{align*}
\item[$\mathbf{H=0}$] then either $\ac22=0$, $\ac23=0$, $a_6=\efct^{k}$ or 
\begin{align*} \ac22=&\frac{t+l}{(t+l)^2+ h^2},
\quad \ac23= -\frac{h}{(t+l)^2+ h^2},\\ 
a_6=&\frac{\efct^k}{\sqrt{(t+l)^2+ h^2}}.
\end{align*}
\end{description} 
\end{lem}
\begin{proof}
Equations \eqref{Da22t3} and \eqref{Da23t3} can be combined to
$T(\ac22+ \I \ac23) = -(\ac22+ \I \ac23)^2 +H$. Using the notation
$y=-(\ac22+ \I \ac23)$, we have to solve the first order differential
equations $\pt y=y^2+1$ ($H=-1$), $\pt y=y^2-1$ ($H=1$) resp. $\pt
y=y^2$ ($H=0$). For $H=-1$ we get the generic solution $y=\tan (z+c)$
resp. the exceptional solution $y=\pm \I$, for $H=1$ we get $y=\I \tan
(z+c)$ resp. $y=\pm 1$ and for $H=0$ we get $y=-\frac1z$
resp. $y=0$. To obtain the real valued functions $\ac22$ and $\ac23$
along the ruling, we have to split $y(z(t))=y(g(t)+ \I h(t))$ into
real and imaginary part. In the first case ($H=-1$) we get (along the
ruling) for example the generic solution: $y(t)=\tan(g(t)+\I h(t))$
with $g(t)=t+l$ and $\pt h=0$. Now $y(t)=\frac{\sin(g(t)+\I
h)}{\cos(g(t)+\I h)}= \frac{\sin(g(t)) \cos(g(t))+ \I \sinh(h)
\cosh(h)}{\cos^2(g(t)) +\sinh^2(h)}$. Next we integrate \eqref{Da6t3}
along the ruling to obtain $a_6$. In the first case this gives
generically: $-\ln a_6 = \half \ln(\cos^2(t+l) +
\sinh^2(h))+\tilde{k}$. The other results are computed in the same
way.
\end{proof}
In order to solve the rest of the integrability conditions we need to
introduce adapted coordinates. Since $T$ is geometrically
distinguished we would like to find Gaussian vector fields $\tT, \tV,
\tW$ such that $\tT=T$. Analogous to \cite{KV99} we look for functions
$c_1, c_2, b_1, b_2$ such that
\begin{equation}\label{Gbasis}
\left(\! \ba{c} \tV \\ \tW  \ea\! \right)= \begin{pmatrix} c_1& c_2\\-c_2& c_1& \end{pmatrix}\left( \left(\! \ba{c} V \\ W  \ea\! \right) + \left(\! \ba{c} b_1 \\ b_2  \ea\! \right) T \right).
\end{equation}
\begin{lem}\label{lem:Gbasis}
For $\tT, \tV, \tW$ to be (local) Gaussian vector fields it is
necessary and sufficient that the following system of equations is
satisfied:
\begin{alignat}{2}
T(c_1)&=\ac22 c_1- \twothird \ac23 c_2,& T(c_2)&= \twothird \ac23 c_1
+\ac22 c_2, \label{Tc}\\ 
-V(c_1)+ W(c_2)&= \bc33 c_1 + \bc23 c_2,
&V(c_2)+W(c_1)&= \bc23 c_1 -\bc33 c_2,\label{VWc}\\ 
T(b_1)&= -\ac22
b_1 -\twothird \ac23 b_2,& T(b_2)&= \twothird \ac23 b_1 -\ac22
b_2,\label{Tb}\\ 
V(b_2)- W(b_1)&= -2 \ac23 +\bc23 b_1 +\bc33
b_2.&&\label{VWb}
\end{alignat}
\end{lem}
\begin{proof}
We have to show that our system of equations is equivalent to the
vanishing of all commutators. Using $[\alpha X,\beta Y]=\alpha
X(\beta) Y-\beta Y(\alpha) X +\alpha \beta (\wN_X Y-\wN_Y X)$ we get
\begin{equation}\label{LieVWT}
\begin{split}
[V,W]&=\bc23 V + \bc33 W +2 \ac23 T,\quad [T,V]=-\ac22 V- \twothird
\ac23 W, \\
[T,W]&= \twothird \ac23 V -\ac22 W.
\end{split}\end{equation}
Therefore
\begin{equation*} 
\begin{split}[\tT,\tV]=&[T,c_1 V+c_2W+ (c_1 b_1+ c_2 b_2)T] \\
=&(T(c_1) -\ac22 c_1+\twothird \ac23 c_2)V + (T(c_2) -\twothird \ac23
c_1- \ac22 c_2) W \\ &+ (T(c_1)b_1 + c_1 T(b_1) + T(c_2) b_2+ c_2
T(b_2))T.
\end{split}\end{equation*}
Hence $[\tT, \tV]=0$ is equivalent to \eqref{Tc} and the equation 
\begin{equation}\label{Tbsys1}
c_1 T(b_1) + c_2 T(b_2)= c_1(-\ac22 b_1 -\twothird \ac23 b_2)+
c_2(\twothird \ac23 b_1 -\ac22 b_2).
\end{equation}
Similarly we get
\begin{equation*} 
\begin{split}[\tT,\tW]=&[T,-c_2 V+c_1W+ (-c_2 b_1+ c_1 b_2)T] \\
=&(-T(c_2) +\ac22 c_2+\twothird \ac23 c_1)V + (T(c_1) +\twothird \ac23
c_2- \ac22 c_1) W \\ &+ (-T(c_2)b_1 - c_2 T(b_1) + T(c_1) b_2+ c_1
T(b_2))T.
\end{split}\end{equation*}
Hence $[\tT, \tW]=0$ is equivalent to \eqref{Tc} and the equation 
\begin{equation}\label{Tbsys2}
-c_2 T(b_1) + c_1 T(b_2)= -c_2(-\ac22 b_1 -\twothird \ac23 b_2)+
 c_1(\twothird \ac23 b_1 -\ac22 b_2).
\end{equation}
Since $c_1^2+c_2^2\neq 0$, equations \eqref{Tbsys1} and \eqref{Tbsys2}
are equivalent to \eqref{Tb}. Finally, a lengthy but straightforward
calculation gives (we use the abbreviation $\al=c_1 b_1+ c_2 b_2$ and
$\bet=-c_2 b_1+ c_1 b_2$)
\begin{equation*} \begin{split}[\tV,\tW]=&[c_1 V + c_2 W +\al T,-c_2 V+c_1W+ \bet T] \\
=&\{c_1(-V(c_2)-W(c_1))+ c_2(V(c_1)- W(c_2)) +(c_1^2 +c_2^2) \bc23 \\
&-\al (T(c_2) -c_2 \ac22 - c_1 \twothird \ac23 ) - \bet (T(c_1)- c_1
\ac22 +c_2 \twothird \ac23 )\}V \\ 
&+\{-c_2(-V(c_2)-W(c_1))+ c_1(V(c_1)- W(c_2)) +(c_1^2 +c_2^2) \bc33 \\ 
&+\al (T(c_1) -c_1 \ac22+ c_2 \twothird \ac23 ) 
- \bet (T(c_2)- c_2 \ac22 -c_1 \twothird \ac23)\}W \\ 
&+\{c_1 (V(\bet)- W\al)) +c_2 (V(\al )+W(\bet ))+(c_1^2
+c_2^2) 2\ac23 \\ &+\al T(\bet )- \bet T(\al )\}T.
\end{split}\end{equation*}
Using equation \eqref{Tc}, the vanishing of the V- and W-component of
$[\tV,\tW]$ is equivalent to \eqref{VWc}. Now $[\tV,\tW]$ vanishes in
addition to $[\tT,\tV]$ and $[\tT,\tV]$ if and only if its $V$-,
$W$-components vanish and the equation
\[0=(c_1^2 +c_2^2)(2\ac23 -b_1 \bc23 - b_2 \bc33 +V(b_2)- W(b_1))\]
holds. However, $c_1^2+c_2^2\neq 0$, thus the last equation is
equivalent to \eqref{VWb}.
\end{proof}
To obtain adapted coordinates we need to solve the system of equations
\eqref{Tc}--\eqref{VWb}. First we will treat the generic cases.
\begin{lem}\label{lem:solvecb}
Assume that we are in one of the generic cases
(cp. Lem.~\ref{lem:Int1}) and let $f$ be any function with
$Tf=-1$. Then the functions $c_1, c_2, b_1, b_2$, defined by
\begin{equation}\label{solvecb} 
(c_1+\I c_2)^3= \frac1{a_6((\ac22+\I \ac23)^2-H)}, \quad b_1=V(f),
\quad b_2=W(f),
\end{equation}
are a solution of the system of differential equations
\eqref{Tc}--\eqref{VWb}.
\end{lem}
\begin{proof}
Observe first that the system of equations splits into two independent
subsystems, namely the system \eqref{Tc}, \eqref{VWc}, and the system
\eqref{Tb}, \eqref{VWb}. Assume first that $d\ac22$, $d\ac23$ and
$da_6$ are everywhere linearly independent.  We start with the system
\eqref{Tc}, \eqref{VWc}. We reformulate the system with respect to the
coordinate functions $\ac22$, $\ac23$, $a=\ln a_6$, using the
integrability conditions \eqref{Da6t3}--\eqref{Da23t3}:
\begin{align}
\ac22 c_1-\twothird \ac23 c_2 =&T(c_1) \nonumber\\
=&((-\ac22^2+\ac23^2 +H)\paa - 2 \ac22 \ac23 \pbb -\ac22 \pa)c_1, \label{Tc1}\\
\twothird \ac23 c_1 +\ac22 c_2 =&T(c_2)\nonumber\\
=&(-\ac22^2+\ac23^2 +H)\paa - 2 \ac22 \ac23 \pbb -\ac22 \pa)c_2,\label{Tc2}\\
\bc33 c_1 + \bc23 c_2 =& -V(c_1)+ W(c_2)\nonumber\\
=& V(\ac22)(-\paa c_1 +\pbb c_2) +W(\ac22)(\paa c_2 +\pbb c_1) \nonumber\\
&-3\bc33\pa c_1-3 \bc23\pa c_2,\label{VWc1}\\
\bc23 c_1 -\bc33 c_2 =&V(c_2)+W(c_1) \nonumber\\
=& V(\ac22)(\paa c_2 +\pbb c_1) +W(\ac22)(\paa c_1 -\pbb c_2)\nonumber\\
& +3\bc33\pa c_2-3 \bc23\pa c_1.\label{VWc2}
\end{align}
The last two equations are satisfied if $(c_1,c_2)$ satisfy the
Cauchy-Riemann equations with respect to the coordinates $\ac22$ and
$\ac23$ for any fixed $a$, i.~e.
\[\paa c_1=\pbb c_2, \quad \paa c_2=-\pbb c_1\]
and their $a$-dependence is given by
\begin{equation}\label{pa}
\pa c_1=-\frac13 c_1,\quad \pa c_2=-\frac13 c_2.
\end{equation}
Therefore we assume now that $c_1(z)+\I c_2(z)$ is a holomorphic
function with respect to $z=\ac22+\I \ac23$ and that \eqref{pa}
holds. Since $(z^2-H)\pz(c_1+\I c_2)= ((\ac22^2-\ac23^2 -H)\paa + 2
\ac22 \ac23 \pbb)c_1 +\I((+\ac22^2-\ac23^2 -H)\paa + 2 \ac22 \ac23
\pbb)c_2$ and $-\twothird z (c_1+\I c_2)= -\twothird (\ac22 c_1-\ac23
c_2+ \I(\ac22 c_2+\ac23 c_1))$, equations \eqref{Tc1} and \eqref{Tc2} 
are equivalent to
$$(z^2-H)\pz (c_1+ \I c_2)= -\twothird z(c_1+ \I c_2).$$
We can integrate the last equation and obtain that $(c_1+ \I c_2)^{-3} 
= g(a) ((\ac22+ \I \ac23)^2 -H)$. Finally equation \eqref{pa} implies $g(a) 
= \efct^{a}= a_6$. The second subsystem can be rewritten as the integrability 
conditions for a function $f$ with $Tf=-1$. In fact, if we set $b_1=V(f)$, 
$b_2= W(f)$ and $b_3= T(f)$, we get by \eqref{Tb} and \eqref{VWb}:
\begin{align*}
[V,T](f) &= \ac22 V(f) + \twothird \ac23 W(f) = \ac22 b_1 + \twothird \ac23
b_2= -T(b_1),\\ 
[W,T](f) &= -\twothird \ac23 V(f)+\ac22 W(f)=
-\twothird \ac23 b_1+\ac22 b_2= -T(b_2),\\ 
[V,W](f) &= 2 \ac23 b_3 +
\bc23 b_1+ \bc33 b_2 = 2 \ac23 (b_3+1) +V(b_2) -W(b_1).
\end{align*}
Hence these integrability conditions are satisfied for $b_3=-1$ and
arbitrary functions $b_1$ and $b_2$.

Now we deal with the case that $d\ac22$, $d\ac23$ and $da_6$ are
linearly dependent in some subset of $\M$. We can still define
functions $c_1, c_2, b_1, b_2$ by equations \eqref{solvecb}. A direct
computation shows that these functions satisfy the system of equations
\eqref{Tc}--\eqref{VWb}.
\end{proof}
We denote the coordinates provided by Lemma~\ref{lem:solvecb} by
$(t,v,w)$, i.~e. $\pt=\tT=T$, $\ptv=\tV$ and $\ptw=\tW$. We will now
make a special choice for $f$ in order to simplify equations in
Lemma~\ref{lem:Int1}, namely to make $l=0$.
\begin{lem}\label{lem:l0}
There exist a function $f$ with $T(f)=-1$ and a constant $c$ such that
after replacing $t$ by $t+c$ one has $l=0$ in Lemma~\ref{lem:Int1}. 
The function $f$ is given by 
\begin{description}
\item[$\mathbf{H=-1}$:] $f= -\half \arcsin 
\left(\frac{2a}{\sqrt{(-\ac22^2+\ac23^2 -1)^2+ 4\ac22^2 \ac23^2}}\right)$, 
\item[$\mathbf{H=1}$:] $f= -\half \arcsinh 
\left(\frac{2a}{\sqrt{(-\ac22^2+\ac23^2 +1)^2+ 4\ac22^2 \ac23^2}}\right)$,
\item[$\mathbf{H=0}$:]$f= -\frac{a}{\sqrt{(-\ac22^2+\ac23^2 )^2
+ 4\ac22^2 \ac23^2}}$.
\end{description} 
The variable $t$ is given by $t= -f +const.$
\end{lem}
\begin{proof}
A direct calculation shows that our choice of $f$ satisfies $T(f)=-1$
. Since $V(f)=b_1$, $W(f)=b_2$, we get
from \eqref{Gbasis} that
$$\left(\! \ba{c} \ptv \\ \ptw \ea\! \right)f= \begin{pmatrix} c_1&
c_2\\-c_2& c_1& \end{pmatrix}\left( \left(\! \ba{c} V(f) \\ W(f) \ea\!
\right) + \left(\! \ba{c} b_1 \\ b_2 \ea\! \right) T(f)
\right)=\left(\! \ba{c} 0 \\ 0 \ea\! \right).$$ 
Using Lemma~\ref{lem:Int1} we compute that $f=-(t+l)$. The function
$l$ must be a constant, since $Tl=0$ by Lemma~\ref{lem:Int1} and we
have just shown that also $\ptv f=\ptw f=0$.
\end{proof}
In the generic cases we have everything to prove the final results:
\begin{thm}\label{ClassT3hyperbolic}
Let $\M$ be an indefinite affine hypersphere of $\Rf$ with $H=-1$ which
admits a pointwise $S_3$-symmetry and assume, that the invariant
functions everywhere satisfy $(\ac22,\ac23)\neq(0,\pm 1)$. Then $\M$
admits coordinates $(t,v,w)$ and there are functions $h,k\colon (v,w)
\to \R$ such that 
\begin{align*} \ac22=&-\frac{\sin(t)\cos(t)}{\cos^2(t) + \sinh^2(h)}, \quad
\ac23=-\frac{\sinh(h)\cosh(h)}{\cos^2(t) + \sinh^2(h)},\\
a_6=&\frac{\efct^k}{\sqrt{\cos^2(t) +\sinh^2(h)}},
\end{align*}
where $h,k$ satisfy the system of differential equations
\begin{equation}\label{diffeqhkh}
(\frac{\partial^2}{\partial v^2} + \frac{\partial^2}{\partial w^2})h 
= \efct^{-\twothird k} \sinh(2h),\quad 
(\frac{\partial^2}{\partial v^2} + \frac{\partial^2}{\partial w^2})k 
= 3 \efct^{-\twothird k} (2 \efct^{2 k} -\cosh (2h)). 
\end{equation}
\end{thm}
\begin{thm}\label{ClassT3elliptic}
Let $\M$ be an indefinite affine hypersphere of $\Rf$ with $H=1$ which
admits a pointwise $S_3$-symmetry and assume, that the invariant
functions everywhere satisfy $(\ac22,\ac23)\neq(\pm 1,0)$. Then $\M$
admits coordinates $(t,v,w)$ and there are functions $h,k\colon (v,w)
\to \R$ such that 
\begin{align*} \ac22=&\frac{\sinh(t)\cosh(t)}{\cos^2(h) + \sinh^2(t)},
\quad \ac23= -\frac{\sin(h) \cos(h)}{\cos^2(h) + \sinh^2(t)},\\ 
a_6=&\frac{\efct^k}{\sqrt{\cos^2(h) +\sinh^2(t)}},
\end{align*}
where $h,k$ satisfy the system of differential equations
\begin{equation}\label{diffeqhke}
(\frac{\partial^2}{\partial v^2} + \frac{\partial^2}{\partial w^2})h 
= -\efct^{-\twothird k} \sin(2h),\quad 
(\frac{\partial^2}{\partial v^2} + \frac{\partial^2}{\partial w^2})k 
= 3 \efct^{-\twothird k} (2 \efct^{2 k} +\cos (2h)). 
\end{equation}
\end{thm}
\begin{thm}\label{ClassT3parabolic}
Let $\M$ be an indefinite improper affine hypersphere of $\Rf$ which
admits a pointwise $S_3$-symmetry and assume, that the invariant
functions everywhere satisfy $(\ac22,\ac23)\neq(0,0)$. Then $\M$
admits coordinates $(t,v,w)$ and there are functions $h,k\colon (v,w)
\to \R$ such that 
\begin{equation*} \ac22=\frac{t}{t^2+ h^2},
\quad \ac23= -\frac{h}{t^2+ h^2},\quad 
a_6=\frac{\efct^k}{\sqrt{t^2+ h^2}},
\end{equation*}
where $h,k$ satisfy the system of differential equations
\begin{equation}\label{diffeqhkp}
(\frac{\partial^2}{\partial v^2} + \frac{\partial^2}{\partial w^2})h 
= 2 \efct^{-\twothird k} h,\quad 
(\frac{\partial^2}{\partial v^2} + \frac{\partial^2}{\partial w^2})k 
= 3 \efct^{-\twothird k} (2 \efct^{2 k} -1). 
\end{equation}
\end{thm}
\begin{proof}[Proof of the theorems]
Let $\M$ be an indefinite affine hypersphere of $\Rf$ which admits a
pointwise $S_3$-symmetry and assume, that the invariant functions
everywhere satisfy $(\ac22,\ac23)\neq(0,\pm 1)$ ($H=-1$),
$(\ac22,\ac23)\neq(\pm 1,0)$ ($H=1$) resp. $(\ac22,\ac23)\neq(0,0)$
($H=0$). By Lemma~\ref{lem:Int1} we can assume without loss of
generality that the functions $\ac22$, $\ac23$ and $a_6$ are given by the
equations above. We need to show that all integrability conditions
(\eqref{Da6t3}--\eqref{Db2t3}) are satisfied. To simplify the
computations we will work in some situations in $\C$ and use the
identification
$$ x+\I y \longleftrightarrow \begin{pmatrix}x &-y\\
y&x\end{pmatrix}$$ 
The equations \eqref{Da22t3} then are equivalent to
\begin{equation}\label{Cauchy-Riemann}
\VW(\ac23)=\I \VW(\ac22)
\end{equation}
 (and $T(\ac22)= -\ac22^2 +\ac23^2 +H$). To get $V( \ac22)$, $V( \ac23)$, $W( \ac22)$ and $W(\ac23)$, first we use \eqref{Gbasis} and the notation $c=c_1 +\I c_2$, $b=b_1+\I b_2$, $\del=\pv + \I \pw$:
\begin{equation}\label{VW}
\VW =\frac1{\cc} \del-b T.
\end{equation}For
the computation we would like to know $\del \ac22$ and $\del \ac23$
resp. $X(\ac22)$ and $X(\ac23)$ for $X\in
\{\pv,\pw\}$. In the proof of Lemma~\ref{lem:Int1} we have used the
notation $y=-(\ac22+\I \ac23)$, and we have shown that $y=\tan(z+c)=\tan(t
+\I h)$ ($H=-1$), $y=\I \tan(z+c)= \I \tan(h +\I t)$ ($H=1$)
resp. $y=-\frac1z= -\frac1{t+\I h}$ ($H=0$), where $T(h)=0$. Now $\frac{\partial
y}{\partial z} \frac{\partial z}{\partial t} = T(y)=-(T(\ac22)+\I
T(\ac23))$, i.~e. $\frac{\partial y}{\partial z}=-(T(\ac22)+\I T(\ac23))$
($H=-1,0$) resp. $\frac{\partial y}{\partial z}=-T(\ac23)+\I T(\ac22)$
($H=1$), and furthermore $-(X(\ac22)+\I X(\ac23)) = X(y)=\frac{\partial
y}{\partial z} X(z)$ with $X(z)=\I X(h)$ ($H=-1,0$) resp. $X(z)=X(h)$
($H=1$). Thus we get with $\eps_H= \begin{cases} 1& \text{if $H=-1,0$},\\ 
-1 &\text{if $H=1$},\end{cases}$:
\begin{equation} \label{Xab}
 X(\ac22) =-\eps_H T(\ac23) X(h),\quad
 X(\ac23)=\eps_H T(\ac22)X(h),
\end{equation}
Applying \eqref{VW} and \eqref{Xab} to \eqref{Cauchy-Riemann},
we obtain that $\frac1{\cc} \eps_H T(\ac22) \delh -b T(\ac23) =
\VW(\ac23)=\I
\VW(\ac22)=-\frac1{\cc}\I \eps_H T(\ac23) \delh -\I b T(\ac22)$. 
Thus equation \eqref{Cauchy-Riemann} is equivalent to 
\begin{equation}\label{bcomplex} 
b=\eps_H \frac1{\cc} \I \delh.
\end{equation}
Now equations \eqref{Da6t3} are equivalent to:
\begin{equation}\label{bc}\begin{split} 
3(\bc33-\I \bc23)&=\VW(\ln a_6) 
=\frac1{\cc} \del(\ln a_6)-b (-\ac22)\\
&=\frac1{\cc}(\del(\ln a_6)+\eps_H \I \delh \ac22).
\end{split}\end{equation}
In complex notion the integrability conditions \eqref{Da23t3} are
\begin{equation} \label{Int3c}
 \VW(\ac23)= -(3\I T+ 3 \I \ac22+2 \ac23)(\bc33-\I \bc23).
\end{equation}
From the above we already know that
\begin{equation} \label{VWbnew}
 \VW(\ac23)= \frac1{\cc} \eps_H \delh (T(\ac22)-\I T(\ac23)).
\end{equation}
For the right side of equation \eqref{Int3c} we need to know
$T(\cc)$. Using \eqref{Tc} we obtain that $T(\cc)=-\frac1{3 \cc }(3\ac22
-\I 2\ac23)$. With the help of \eqref{Xab} and the equations for
$T(\ac22)$ and $T(\ac23)$ (cp. \eqref{Da22t3} and \eqref{Da23t3}) we
get:
\begin{equation} \label{Tbc}\begin{split}
 3T(\bc33-\I \bc23 )=& \frac1{3\cc}\{ \del ( \ln a_6 )(-3\ac22+\I 2
\ac23 )\\ &+ \eps_H \delh (\I \ac22 (-3\ac22+\I 2 \ac23 )+ 3 T(\ac23 )+\I
3 T(\ac22 ))\}\\ =& \frac1{3\cc}\{\del (\ln a_6)(-3\ac22+\I 2 \ac23)\\
&+ \eps_H \delh (-4 \ac22 \ac23 +\I (-6 \ac22^2 +3 \ac23^2 +3 H))\}.
\end{split}\end{equation}
From \eqref{VWbnew} and \eqref{Tbc} (and \eqref{bc}) it is easy to
check that \eqref{Int3c} is true and thus \eqref{Da23t3}. For the
investigation of \eqref{Db1t3} and \eqref{Db2t3} we first compute the 
left hand sides in complex notion. We use \eqref{VW}, the notion $\lap = 
\frac{\partial^2}{\partial^2 v} + \frac{\partial^2}{\partial^2 w}$, and also 
that by \eqref{bc} $\bc33 +\I \bc23 = \frac1{3c}(\overline{\del (\ln a_6)}- \I 
\ac22 \eps_H \overline{\delh})$, and obtain after a while:
\begin{equation}\label{VWbcbar}
\begin{split}
 \VW(\bc33+\I \bc23 )=& \frac1{3\I |c|^2}\{ (\ac22 \eps_H \overline{\delh}+ \I 
 \overline{\del (\ln a_6)})(-\frac1c) (\del c -\eps_H \I \delh T(c)) \\
 &+ \eps_H \ac22 \lap h + \I \lap (\ln a_6 ) - \I T(\ac22) |\delh|^2 \}.
\end{split}\end{equation}
The complex function $c$ was defined by $c=(\efct^{\ln
a_6}(z^2-H))^{-\frac13}$ (cp. \eqref{solvecb}), where $z=\ac22 +\I
\ac23$. Thus $\del c =-\frac{c}3 ( \del (\ln a_6) +\frac{2 z
\eps_H\del(z) \delh}{z^2-H}) = -\frac{c}3 ( \del (\ln a_6) -2
\eps_H\I(\ac22 +\I \ac23) \delh)$ (cp. \eqref{Xab}). By \eqref{Tc} we
have that $T(c)=(\ac22+\I \twothird \ac23) c$. So we can simplify
$\del c -\eps_H \I \delh T(c) = -\frac{c}3 ( \del (\ln a_6) +\eps_H
\ac22 \I \delh)$, and hence \eqref{VWbcbar} becomes:
\begin{equation} \label{VWbcbar2}\begin{split}
 \VW(\bc33+\I \bc23 )=& \frac1{3 |c|^2} \{-\I\ac22 \eps_H \third 2\I \Im (\overline{\delh} \del ( \ln a_6 )) +\third |\del ( \ln a_6 )|^2\\
&+ |\delh|^2 (\frac43 \ac22^2 -\ac23^2 -H) 
  +\lap (\ln a_6) -\I \eps_H \ac22 \lap h  \}.
\end{split}\end{equation}
Splitting into real and imaginary part we have:
\begin{align} \label{IntVWb1}\begin{split}
 V(\bc33)-W(\bc23) =& \frac1{3 |c|^2} \{2\ac22 \eps_H \third \Im
(\overline{\delh} \del ( \ln a_6 )) +\third |\del ( \ln a_6 )|^2\\ &+
|\delh|^2 (\frac43 \ac22^2 -\ac23^2 -H) +\lap (\ln a_6) \},
\end{split}\\
 \label{IntVWb2}
 V(\bc23)+W(\bc33)=& -\frac1{3 |c|^2}\eps_H \ac22 \lap h.
\end{align}
This means that \eqref{Db2t3} is equivalent to
\begin{equation} \label{Laph}
 \lap h=-2 \eps_H |c|^2 b.
\end{equation}
Furthermore by \eqref{bc} we get $\bc33^2+\bc23^2= |\bc33+\I \bc23|^2=
\frac1{9 |c|^2}(|\del ( \ln a_6 )|^2 + \ac22 |\delh|^2 +
2\eps_H \ac22 \Im (\overline{\delh} \del ( \ln a_6 )))$ and therefore:
\begin{equation}\label{IntVWb1a}
 V(\bc33)-W(\bc23) -(\bc33^2+\bc23^2) = \frac1{3 |c|^2}(|\delh|^2
 (\ac22^2 -\ac23^2 -H) +\lap (\ln a_6) ).
\end{equation}
To simplify the right hand side we need to compute $\lap (\ln a_6)
$. A lengthy but straightforward computation, using the definition of
$a_6$ in each case, leads to:
\begin{equation}\label{IntVWb1b}
 V(\bc33)-W(\bc23) -(\bc33^2+\bc23^2) = \frac1{3 |c|^2}(\lap k +\eps_H
 \ac23 \lap h ).
\end{equation}
If we use \eqref{Laph}, we get from \eqref{IntVWb1b} that
\eqref{Db1t3} is equivalent to
\begin{equation} \label{Lapk}
 \lap k=3 |c|^2(-\ac22^2 - \ac23^2 + H + 2 a_6^2).
\end{equation}
Finally, to prove the equivalence of \eqref{Laph} and \eqref{Lapk}
with the systems of differential equations (cp. \eqref{diffeqhkh},
\eqref{diffeqhke}, \eqref{diffeqhkp}), one needs to juggle with
trigonometric equalities. Explicitly we will give one more step:
\begin{description}
\item[$\mathbf{H=-1:}$] $|c|^2= \efct^{-\twothird k}(\cos^2 (t) + 
\sinh^2 (h))$, 
\item[$\mathbf{H=1}$:] $|c|^2= \efct^{-\twothird k}(\cos^2 (h) + 
\sinh^2 (t))$,
\item[$\mathbf{H=0}$:]$|c|^2= \efct^{-\twothird k}(h^2 + t^2)$.
\end{description}
\end{proof}
We are left with the exceptional cases
(cp. Lemma~\ref{lem:Int1}). Using \eqref{LieVWT} it is easy to show
that the Gauss equations \eqref{Da23t3} and \eqref{Db2t3} are just the
integrability conditions for a function $k$ with
\begin{equation}\label{TVWk} T(k)=0,\quad V(k)=3 \bc33, \quad W(k)=- 3
\bc23.
\end{equation}
 For $a_6= \exp(-\ac22 t +k)$ (cp. Lemma~\ref{lem:Int1}) the latter
 equations are the Codazzi equations \eqref{Da6t3}. To solve the last
 non-trivial integrability condition \eqref{Db1t3} we need to find
 adapted coordinates, i.~e. by Lemma~\ref{lem:Gbasis} we need to solve
 the system of equations \eqref{Tc}--\eqref{VWb}. Notice that for
 $H=0$ or $H=1$ we have $\ac23 =0$, hence by \eqref{LieVWT} the
 distribution $T^\perp =\Span\{V,W\}$ is integrable and we do not need
 the functions $b_1, b_2$ in \eqref{Gbasis}.
\begin{lem}\label{lem:solvecbex}
Assume that we are in one of the exeptional cases
(cp. Lem.~\ref{lem:Int1}) and let $f$ be any function with
$Tf=-1$. Then the functions $c_1, c_2, b_1, b_2$, defined by
\begin{equation}\label{solvecbex} 
(c_1+\I c_2)^3= \frac1{a_6} \efct^{2(\ac22+\I \ac23)t}, \quad
b_1=V(f), \quad b_2=W(f),
\end{equation}
are a solution of the system of differential equations
\eqref{Tc}--\eqref{VWb}.
\end{lem}
\begin{proof}
Observe first that the system of equations splits into two independent
subsystems, namely the system \eqref{Tc}, \eqref{VWc}, and the system
\eqref{Tb}, \eqref{VWb}. We rewrite the first systems, using again the 
complex notation $c=c_1+ \I c_2$: \eqref{Tc} is equivalent to
$T(c)=(\ac22+\I \twothird \ac23 ) c$ and \eqref{VWc} to
\begin{equation}\label{VWccomplex}
(V+\I W)(c) = (-\bc33 +\I \bc23) c = -\third (V+\I W)(k) c.
\end{equation}
Integrating these equations we get that $c=\exp(-\third k+ (\ac22+\I
\twothird \ac23 ) t)$. Since $a_6 = \exp(-\ac22 t +k)$, we get the
expression for $c$. The proof for the second subsystem (only
non-trivial for $H=-1$) is the same as before
(cp. Lemma~\ref{lem:solvecb}).
\end{proof}
As before we denote the coordinates, this time provided by
Lemma~\ref{lem:solvecbex}, by $(t,v,w)$, i.~e. $\pt=\tT=T$, $\ptv=\tV$
and $\ptw=\tW$. Now we are ready to prove the final results for the
exceptional cases:
\begin{thm}\label{ClassT3hyperbolicEx}
Let $\M$ be an indefinite affine hypersphere of $\Rf$ with $H= -1$ which
admits a pointwise $S_3$-symmetry and assume, that the invariant
functions everywhere satisfy $(\ac22,\ac23)=(0,\pm 1)$. Then $\M$
admits coordinates $(t,v,w)$ and there is a function $k\colon (v,w)
\to \R$ such that $a_6=\efct^k$,
where $k$ satisfies the elliptic pde 
\begin{equation}\label{diffeqhEx}
(\frac{\partial^2}{\partial v^2} + \frac{\partial^2}{\partial w^2})k = 2 \efct^{-\twothird k} (3 \efct^{2 k} -4). 
\end{equation}
\end{thm}
\begin{thm}\label{ClassT3ellipticEx}
Let $\M$ be an indefinite affine hypersphere of $\Rf$ with $H= 1$ which
admits a pointwise $S_3$-symmetry and assume, that the invariant
functions everywhere satisfy $(\ac22,\ac23)=(\pm 1,0)$. Then $\M$
admits coordinates $(t,v,w)$ and there is a function $k\colon (v,w)
\to \R$ such that $a_6=\efct^{k\mp t}$,
where $k$ satisfies the elliptic pde 
\begin{equation}\label{diffeqeEx}
(\frac{\partial^2}{\partial v^2} + \frac{\partial^2}{\partial w^2})k = 6 \efct^{\frac43 k}. 
\end{equation}
\end{thm}
\begin{thm}\label{ClassT3parabolicEx}
Let $\M$ be an indefinite improper affine hypersphere of $\Rf$ which
admits a pointwise $S_3$-symmetry and assume, that the invariant
functions everywhere satisfy $(\ac22,\ac23)=(0,0)$. Then $\M$
admits coordinates $(t,v,w)$ and there is a function $k\colon (v,w)
\to \R$ such that $a_6=\efct^k$,
where $k$ satisfies the elliptic pde 
\begin{equation}\label{diffeqpEx}
(\frac{\partial^2}{\partial v^2} + \frac{\partial^2}{\partial w^2})k =
6 \efct^{\frac43 k}.
\end{equation}
\end{thm}
\begin{proof}[Proof of the theorems] 
We already have seen that the integrability conditions
\eqref{Da6t3}--\eqref{Da23t3} and \eqref{Db2t3} are satisfied. Hence
we need to show that equation \eqref{Db1t3} is true. Using
\eqref{Db2t3} and \eqref{TVWk} we get in complex notation $V(\bc33) -
W(\bc23) = (V+\I W) (\bc33+ \I \bc23) = (V+\I W)\third \overline{(V+\I
W)} (k)$. We have found coordinates such that $V+\I W
=\frac{c}{|c|^2}\del - b T$ (cp. \eqref{Gbasis}), here $\del=\pv +\I
\pw = \tV+ \I \tW$ and $b= b_1+\I b_2$. Thus $(V+\I W)\third
\overline{(V+\I W)} (k) = \third (V+\I W)
(\frac{\cc}{|c|^2}\overline{\del (k)} - \bar{b} T(k))=
\third (V+\I W) (\frac1{c}\overline{\del (k)}) $, since $T(k)=0$ by 
\eqref{TVWk}. We know $(V+\I W)(c)$ by equation \eqref{VWccomplex}, 
and because of $T(\overline{\del (k)}) =0$ we finally obtain that
\begin{equation}\label{VWbEx}
V(\bc33) - W(\bc23) = \frac1{9|c|^2} (|\del (k)|^2 + 3 \lap k).
\end{equation}
Also from \eqref{VWccomplex} follows that $\bc22^2 +\bc23^2 = 
\frac19 |(V+\I  W)(k)|^2= \frac1{9|c|^2} |\del (k)|^2$, 
thus \eqref{Db1t3} is equivalent to
\begin{equation*}
\frac1{3|c|^2} \lap k= H-\ac22^2 -\frac53 \ac23^2 +2 a_6^2.
\end{equation*}
To complete the proof we only need to replace $H$, $\ac22$, $\ac23$,
$a_6$ and $c$ in each case according to their definition.
\end{proof}

\section{Pointwise $\SO(1,1)$-symmetry}
\label{sec:type8}

Let $\M$ be a hypersphere admitting a $\SO(1,1)$-symmetry. According to
Thm.~\ref{thm:Subgroups}, there exists for every $p\in \M$ an LVB $\{
\be, \bv, \bof\}$ of $T_p \M$ such that
\begin{align*} K(\bv,\bv)&= -2b_4 \bv, & K(\bv,\be)&=b_4 \be, & K(\bv,\bof)&= 
b_4 \bof,\\
K(\be,\be)&= 0, & K(\be,\bof)&= b_4 \bv, &
K(\bof,\bof)&= 0,
\end{align*}
where $b_4 > 0$.

We would like to extend the LVB locally. It is well known that
$\ricLC$ (cp. \eqref{def:RicLC}) is a symmetric operator and we
compute (some of the computations in this section are done with the
CAS
Mathematica\footnote{http://www.math.tu-berlin.de/$\sim$schar/IndefSym\_typ8.html}):
\begin{lem} Let $p \in \M$ and $\{\be,\bv,\bof\}$ the basis constructed
  earlier. Then 
\begin{alignat*}{2}
  &\ricLC(\bv,\bv) =2(H+ 3b_4^2) , \qquad\quad &&\ricLC(\bv,\be)=0, \\
  &\ricLC(\bv,\bof)=0,\qquad\quad &&
\ricLC(\be,\be)=0 ,\\
  &\ricLC(\be,\bof)=2(H+b_4^2),&&\ricLC(\bw,\bw)=0 .
\end{alignat*}
\end{lem}

\begin{proof}
The proof is a straight-forward computation using the Gauss
equation~\eqref{gaussLC}.  It follows
e.~g. that
\begin{align*}
  \hat R(\bv,\be)\bv&= -H \be -K_{\bv}(b_4\be)+K_{\be}(-2 b_4 \bv) =
  -H \be -b_4^2\be -2 b_4^2\be \\ &= -(H+ 3 b_4^2)\be,\\ \hat
  R(\bv,\bof)\bv&= -H \bof -K_{\bv}(b_4 \bof) +K_{\bof}(-2 b_4 \bv)
  =-H \bof - b_4^2 \bof -2 b_4^2 \bof\\ &= -(H+ 3 b_4^2)\bof,\\ \hat
  R(\bof,\be)\bv&=-K_{\bof}(b_4 \be)+K_{\be}(b_4\bof)= 0.
\end{align*}
From this it immediately follows that 
$$\ricLC(\bv,\bv) = 2(H+3b_4^2)$$
and
$$\ricLC(\bv,\be)=0.$$
The other equations follow by similar computations.
\end{proof}
\begin{rem}\label{rem:sc8}
For the scalar curvature of $\M$ we obtain $\hat{\kappa}=H+\frac53
b_4^2$. Therefore $J=\frac53 b_4^2$ and $b_4$ is a globally defined
function on $\M$.
\end{rem}

We want to show that the basis, we have constructed at each point $p$,
can be extended differentiably to a neighborhood of the point $p$ such
that, at every point, $K$ with respect to the frame $\{E,V,F\}$ has
the previously described form.
\begin{lem}\label{lem:KfT8} 
Let $\M$ be an affine hypersphere in $\mathbb R^4$ which admits a
  pointwise $\SO(1,1)$-symmetry. Let $p \in M$. Then there exists a
  lightvector-frame $\{E,V,F\}$ defined in a neighborhood of the point
  $p$ and a postive function $b_4$ such that $K$ is given by:
\begin{align*} 
K(V,V)&= -2b_4 V, & K(V,E)&=b_4 E, & K(V,F)&= 
b_4 F,\\
K(E,E)&=0, & K(E,F)&= b_4 V, &
K(F,F)&=0.
\end{align*}
\end{lem}

\begin{proof}First we want to show that at every point the vector
  $\bv$ is uniquely defined and differentiable. We introduce a
  symmetric operator $\hat A$ by:
\begin{equation*}
\ricLC(Y,Z)= h(\hat A Y,Z).
\end{equation*}
Clearly $\hat A$ is a differentiable operator on $\M$ with $\hat
A(\be)=2(H+b_4^2)\be$, $\hat A(\bv)=2(H+3b_4^2)\bv$ and $\hat
A(\bof)=2(H+b_4^2)\bof$. Since $2(H+3b_4^2) \neq 2(H+b_4^2)$, the
operator has two distinct eigenvalues. A standard result then implies
that the eigendistributions are differentiable. We take $V$ a local
unit vectorfield spanning the 1-dimensional spacelike
eigendistribution, and local orthonormal vectorfields $\tilde{T}$ and
$\tilde{W}$ spanning the second eigendistribution. We define
$E=\frac{1}{\sqrt{2}}(-\tilde{T}+\tilde{W})$ and
$F=\frac{1}{\sqrt{2}}(\tilde{T}+\tilde{W})$.
\end{proof}
\begin{rem} It actually follows from the proof of the previous lemma
  that the vector field $V$ is globally defined on $\M$. 
\end{rem}
In this section we always will work with the local frame constructed
in the previous lemma. We denote the coefficients of the Levi-Civita
connection with respect to this frame by:
\begin{align*}
   \widehat{\nabla}_E E &= \ac11 E + \bc11 V,&
   \widehat{\nabla}_E V &= \ac12 E - \bc11 F,&
   \widehat{\nabla}_E F &=-\ac12 V - \ac11 F, \\
   \widehat{\nabla}_V E &= \ac21 E + \bc21 V, &
   \widehat{\nabla}_V V &= \ac22 E - \bc21 F, &
   \widehat{\nabla}_V F &=-\ac22 V - \ac21 F, \\
   \widehat{\nabla}_F E &= \ac31 E + \bc31 V, &
   \widehat{\nabla}_F V &= \ac32 E - \bc31 F, &
   \widehat{\nabla}_F F &=-\ac32 V - \ac31 F.
 \end{align*}
We will evaluate first the Codazzi and then the Gauss equations
(\eqref{CodK} and \eqref{gaussLC}) to obtain more informations.

\begin{lem}\label{lem:CodKt8}
Let $\M$ be an affine hypersphere in $\mathbb R^4$ which admits a
pointwise $\SO(1,1)$-symmetry and $\{E,V,F\}$ the
corresponding LVB. Then $0= \bc11 =\bc21 =\ac22 =\ac32$,
$\bc31=-\ac12$ and $V(b_4)=-4\ac12 b_4$, $0=E(b_4) = F(b_4)$.
\end{lem}

\begin{proof}
An evaluation of the Codazzi equations \eqref{CodK} leads to the
following equations (they relate to eq1--eq5 and eq7--eq8 in the
Mathematica notebook):
\begin{eqnarray}
   E(b_4)= -2\bc21 b_4, \quad 0= \bc11 b_4 ,\label{CodKeq1t8} \\
   0=(\ac12+\bc31)b_4,\quad 0=E(b_4),\label{CodKeq2t8}\\ 
   V(b_4)=-4\ac12 b_4 ,\label{CodKeq3t8}\\ 
   0=\ac32 b_4,\quad F(b_4)= 2 \ac22 b_4 \quad V(b_4)=4\bc31 b_4, 
   \label{CodKeq4t8}\\ 
   0=F(b_4),\label{CodKeq5t8}\\ 
   0=\ac22 b_4, \label{CodKeq7t8}\\
   0=\bc21 b_4,\label{CodKeq8t8}
 \end{eqnarray}
this proves the theorem.
\end{proof}
 
An evaluation of the Gauss equations \eqref{gaussLC} leads to the following :

\begin{lem}\label{lem:GaussLCt8}
Let $\M$ be an affine hypersphere in $\mathbb R^4$ which admits a
pointwise $\SO(1,1)$-symmetry and $\{E,V,F\}$ the
corresponding LVB. Then
\begin{align} 
V(\ac12)&= -\ac12^2 -3 b_4^2-H,\label{Gauss1.1t8}\\
F(\ac12)&=0,\label{Gauss1.2t8}\\ E(\ac12)&=0,\label{Gauss1.3t8}\\ 
E(\ac21) - V(\ac11) &=\ac11(\ac12-\ac21),\label{Gauss1.4t8}\\
E(\ac31)-F(\ac11)&= \ac12^2 -2 \ac11\ac31 -b_4^2+H,\label{Gauss1.5t8}\\ 
V(\ac31)-F(\ac21)&=-\ac31(\ac12 +\ac21).\label{Gauss1.6t8}\end{align}
\end{lem}

\begin{proof}
The equations relate to eq11--eq12 and eq14--eq16 in 
the Mathematica notebook. 
\end{proof}

As the vector field $V$ is globally defined, we can define the
distributions $L_1=\Span\{V\}$ and $L_2=\Span\{E,F\}$. In the
following we will investigate these distributions. 

\begin{lem}\label{lem:L1t8}
The distribution $L_1$ is autoparallel (totally geodesic) with respect
to $\widehat\nabla$.
\end{lem}
\begin{proof} From $\widehat{\nabla}_{V} V = \ac22 E - \bc21 F=0$ 
(cp. Lemma~\ref{lem:CodKt8}) the claim follows immediately.
\end{proof} 
\begin{lem}\label{lem:L2t8}
  The distribution $L_2$ is spherical with mean curvature normal
  $U_2=-\ac12 V$.
\end{lem}
\begin{proof} For $U_2=-\ac12 V\in L_1=L_2^{\perp}$ we have
  $h(\widehat{\nabla}_{E_a} E_b, V)= h(E_a, E_b) h(U_2,V)$ for $E_a,
E_b\in \{E,F\}$, and $h(\widehat{\nabla}_{E_a} U_2, V)= h(-E_a(\ac12)
V + \ac12 \widehat{\nabla}_{E_a}V, V)=0$ (cp. Lemma \ref{lem:CodKt8}
and \eqref{Gauss1.3t8}, \eqref{Gauss1.4t8}.
\end{proof}
\begin{rem} $\ac12$ depends of $V$ but not of the
  choice of basis $\{E,F\}$ of $L_2$. It therefore is a globally defined
  function on $\M$. 
\end{rem}
We introduce a coordinate function $v$ by $\ptv:=V$. Using the
previous lemma, according to \cite{PR93}, we get:
\begin{lem}\label{lem:warpedt8} $(\M,h)$ admits a warped product structure
  $\M=I \times_{e^f}N^2$ with $f: I \to \mathbb R$
  satisfying
\begin{equation}\label{defft8}
\frac{\partial f}{\partial v}=\ac12.
\end{equation}
\end{lem}
\begin{proof} Prop. 3 in \cite{PR93} gives the warped product structure with 
warping function $\lambda_2:I \to \mathbb R$. If we introduce 
$f=\ln \lambda_2$, following the proof we see that 
$-\ac12 V=U_2=-\grad(\ln\lambda_2)=-\grad f$. 
\end{proof}
\begin{lem}\label{lem:curvN2t8} The curvature of $N^2$ is 
${}^NK(N^2)=e^{2f}(- b_4^2+\ac12^2+H)$.
\end{lem}
\begin{proof} For the formulas cp. Lem.~\ref{lem:curvN2}. Using the Gauss
equation~\eqref{gaussLC} (cp. the Mathematica notebook) we get:
\begin{equation*}\begin{split}
(H-b_4^2)E = &\hat R(E,F)E = {}^N \hat{R}(E,F)E -(E(\ac12) V + \ac12(\ac12
E))+ h(E(\ac12) V \\
&+ \ac12 (\ac12 E),E)F - h(F(\ac12) V + \ac12 (\ac12
F),E)E +\ac12^2 E.
\end{split}\end{equation*}
 Thus ${}^NK(N^2)= \frac{{}^N
h(-{}^N\hat{R}(E,F)E,F)}{{}^N h(E,E){}^N h(F,F)-{}^N h(E,F)^2} =
e^{2f} h({}^N \hat{R}(E,F)E,F) = e^{2f}(- b_4^2+\ac12^2+H)$.
\end{proof}
Summarized we have obtained the following structure equations
(cp. \eqref{strGauss}, \eqref{strWeingarten} and \eqref{defK}):
\begin{alignat}{4}
&D_V V =& && -2b_4 V&+ \xi, \label{D11t8}\\
&D_V E =& +(\ac21 +b_4) E, &&& \label{D12t8}\\
&D_V F =&  & +(-\ac21 +b_4) F,&& \label{D13t8}\\
&D_E V =&+(\ac12 +b_4) E, &&& \label{D21t8}\\
&D_F V =& &+(\ac12 + b_4)F, && \label{D31t8}\\
&D_E E =&+ \ac11 E, &&  & \label{D22t8}\\
&D_E F =&  &-\ac11 F& +(-\ac12 + b_4)V& +\xi, \label{D23t8}\\
&D_F E =&+\ac31 E&&+(-\ac12 + b_4)V& +\xi, \label{D32t8}\\
&D_F F =&  & -\ac31 F,& & \label{D33t8}
\end{alignat}
\begin{equation}\label{Dxit8}
D_{X} \xi= -H X,
\end{equation}
The Codazzi and Gauss equations (\eqref{CodK} and \eqref{gaussLC})
have the form (cp. Lem.~\ref{lem:CodKt8} and \ref{lem:GaussLCt8}):
\begin{align}
 V(b_4)&=-4\ac12 b_4, \quad 0=E(b_4) = F(b_4),\label{Db4t8}\\
 V(\ac12)&= -\ac12^2 -3 b_4^2-H,\; F(\ac12)=0,\; E(\ac12)=0,\label{Da12t8}\\
 E(\ac21) - V(\ac11) &=\ac11(\ac12-\ac21),\label{Da1t8}\\
 V(\ac31)-F(\ac21)&=-\ac31(\ac12 +\ac21),\label{Da2t8}\\ 
 E(\ac31)-F(\ac11)&= \ac12^2 -2 \ac11\ac31 -b_4^2+H.\label{Da3t8}
\end{align}

Our first goal is to find out how $N^2$ is immersed in $\Rf$, i.~e. to
find an immersion independent of $v$. A look at the structure
equations \eqref{D11t8} - \eqref{Dxit8} suggests to start with a linear
combination of $V$ and $\xi$.

We will solve the problem in two steps. First we look for a vector
field $X$ with $D_V X=\alpha X$ for some funtion $\alpha$: We define
$X:=A V +\xi$ for some function $A$ on $\M$. Then $D_V X=\alpha X$ iff
$\alpha=A$ and $\ptv A= A^2 +2b_4 A+ H$, and $A:=-\ac12+ b_4$ solves
the latter differential equation. Next we want to multiply $X$ with
some function $\beta$ such that $D_V (\beta X)=0$: We define a positive
function $\beta$ on $\R$ as the solution of the differential equation:
\begin{equation}\label{dtbetat8}
\tfrac{\partial}{\partial v} \beta = (\ac12- b_4)\beta 
\end{equation}
with initial condition $\beta(v_0)>0$. Then $D_V(\beta X)=0$ and by
\eqref{D21t8}, \eqref{Dxit8} and \eqref{D31t8} we get (since
$\beta$, $\ac12$ and $b_4$ only depend on $v$):
\begin{align}
D_{V}(\beta((b_4- \ac12)V +\xi))&=0,\label{eq31t8}\\
D_{E}(\beta((b_4- \ac12)V +\xi))&=\beta(b_4^2-\ac12^2-H)E ,\label{eq32t8}\\ 
D_{F}(\beta((b_4- \ac12)V +\xi))&=\beta(b_4^2-\ac12^2-H)F.\label{eq33t8} 
\end{align}
To obtain an immersion we need that $\nu:=b_4^2-\ac12^2-H$ vanishes
nowhere, but we only get:
\begin{lem}\label{lem:nut8}
  The function $\nu=b_4^2-\ac12^2-H$ is globally defined,
  $\ptv(e^{2f} \nu)=0$ and $\nu$ vanishes identically or nowhere on $\R$.
\end{lem}
\begin{proof} Since $0=\ptv {}^NK(N^2) = \ptv(e^{2f}(-\nu))$ 
  (Lem.~\ref{lem:curvN2t8}), we get by \eqref{defft8} $\ptv\nu= -2\ac12\nu$.
\end{proof}

\subsection{The first case: $\nu \neq 0$ on $\M$}
\label{sec:case1t8}

We may, by translating $f$, i.e. by replacing $N^2$ with a homothetic
copy of itself, assume that $e^{2f} \nu =\eps_1$, where $\eps_1 =\pm 1$.

\begin{lem}\label{lem:defphit8}
$\varPhi:=\beta ((b_4 -\ac12)V +\xi)\colon
M^3 \to \R^4$ induces the structure of an indefinite proper affine quadric, say
$\tilde{\phi}$, mapping $N^2$ into a 3-dimensional linear subspace
of $\R^4$.
\end{lem}
\begin{proof} 
By \eqref{eq32t8} and \eqref{eq33t8} we have $\varPhi_*(E_a)= \beta
\nu E_a$ for $E_a\in \{E,F\}$. A further differentiation, using \eqref{D22t8}
($\beta$ and $\nu$ only depend on $v$), gives:
\begin{equation*}
D_{E} \varPhi_*(E) = \beta \nu D_{E} E = \beta \nu \ac11 E =\ac11 \varPhi_*(E).
\end{equation*}
Similarly, we obtain the other derivatives, using \eqref{D23t8} -
\eqref{D33t8}, thus:
\begin{alignat}{3}
D_{E} \varPhi_*(E)&= & \ac11\varPhi_*(E) &&, \label{Dphivvt8}\\ D_{E}
\varPhi_*(F)&=& &-\ac11 \varPhi_*(F)& +e^{-2f}\eps_1
\varPhi,\label{Dphivwt8}\\ D_{F} \varPhi_*(E)&=& \ac31 \varPhi_*(E)
&&+e^{-2f}\eps_1 \varPhi, \label{Dphiwvt8}\\ D_{F} \varPhi_*(F)&=&
&-\ac31 \varPhi_*(F) & \label{Dphiwwt8}\\ D_{E_a} \varPhi &=& \beta
e^{-2f}\eps_1 E_a.\label{Dphieat8} & &
\end{alignat}
The foliation at $f=f_0$ gives an immersion of $N^2$ to $\M$, say
$\pi_{f_0}$. Therefore, we can define an immersion of $N^2$ to $\R^4$
by $\tilde{\phi}:=\varPhi\circ\pi_{f_0}$, whose structure equations
are exactly the equations above when $f=f_0$. Hence, we know that
$\tilde{\phi}$ maps $N^2$ into
$\Span\{\varPhi_*(E),\varPhi_*(F),\varPhi\}$, an affine hyperplane of
$\R^4$ and $\ptv\varPhi=0$ (cp. \eqref{eq31t8} implies
$\varPhi(v,x,y)=\tilde{\phi}(x,y)$.  We can read off the coefficients
of the difference tensor $K^{\tilde{\phi}}$ of $\tilde{\phi}$
(cf. \eqref{strGauss} and \eqref{defK}) and see that it vanishes. The
affine metric introduced by this immersion corresponds with the metric
on $N^2$.  Thus $\eps_1 \tilde{\phi}$ is the affine normal of
$\tilde{\phi}$ and $\tilde{\phi}$ is an indefinite proper affine
quadric with mean curvature $\eps_1$.
\end{proof}
Our next goal is to find another linear combination of $V$ and $\xi$,
this time only depending on $v$. (Then we can express $V$ in terms
of $\phi$ and some function of $v$.)
\begin{lem}\label{lem:defdeltat8}
  Define $\delta := H V +(b_4+ \ac12) \xi$. Then there exist a constant vector
  $C \in \R^4$ and a function $a(v)$ such that
$$ \delta(v)= a(v) C.$$
\end{lem}
\begin{proof} Using \eqref{D21t8} resp. \eqref{D31t8} and
  \eqref{Dxit8} we obtain that $D_{E}\delta = 0=D_{F} \delta$. Hence $\delta$
  depends only on the variable $v$. Moreover, we get by
  \eqref{D11t8}, \eqref{Da12t8},\eqref{Db4t8} and \eqref{Dxit8} that 
\begin{align*}
  \ptv\delta&=D_{V} (H V +(b_4+ \ac12)\xi)\\
  &=H(-2b_4 V+\xi) + (-4 \ac12 b_4-\ac12^2 -3 b_4^2 -H)\xi -(b_4+\ac12) H V\\
  &=(-3 b_4-\ac12)(H V +(b_4 +\ac12)\xi)\\
  &=-(3 b_4+\ac12) \delta.
\end{align*}
This implies that there exists a constant vector $C$ in $\Rf$ and a
function $a(v)$ such that $\delta(v)=a(v)C$.
\end{proof}
Notice that for an improper affine hypersphere ($H=0$) $\xi$ is
constant and parallel to $C$. Combining $\tilde{\phi}$ and $\delta$ we
obtain for $V$ (cp. Lem.~\ref{lem:defphit8} and \ref{lem:defdeltat8})
that
\begin{equation}\label{V}
V(v,x,y)= -\frac{a}{\nu}C +\frac{1}{\beta\nu}(b_4+\ac12)\tilde{\phi}(x,y).
\end{equation}
In the following we will use for the partial derivatives the
abbreviation $\hs_u:= \frac{\partial}{\partial u}\hs $, $u=v,x,y$.
\begin{lem}\label{lem:partialFt8} 
\begin{align*}
&\hs_v = -\frac{a}{\nu}C +\ptv(\frac{1}{\beta \nu})\tilde{\phi},\\
&\hs_x = \frac{1}{\beta\nu} \tilde{\phi}_x,\\
&\hs_y = \frac{1}{\beta\nu} \tilde{\phi}_y.
\end{align*}
\end{lem}
\begin{proof} As by \eqref{dtbetat8} and Lem.~\ref{lem:nut8} $\ptv
  \frac{1}{\beta\nu}= \frac{1}{\beta\nu}(b_4+\ac12)$, we obtain the
  equation for $\hs_v =V$ by \eqref{V}. The other equations follow
  from \eqref{eq32t8} and \eqref{eq33t8}.
\end{proof}
It follows by the uniqueness theorem of first order differential
equations and applying a translation that we can write
$$\hs(v,x,y)= \tilde{a}(v) C +\frac{1}{\beta\nu}(v)
\tilde{\phi}(x,y)$$ for a suitable function $\tilde{a}$ depending only
on the variable $v$. Since $C$ is transversal to the image of
$\tilde{\phi}$ (cp. Lem.~\ref{lem:defphit8} and \ref{lem:defdeltat8},
$\nu\not\equiv 0$), we obtain that after applying an equiaffine
transformation we can write: $\hs(v,x,y) =(\gamma_1(v), \gamma_2(v)
\phi(x,y))$, where $\tilde{\phi}(x,y)=(0,\phi(x,y))$.  Thus we have
proven the following:

\begin{thm}\label{thm:ClassC1t8} Let $\M$ be an indefinite affine hypersphere 
of $\Rf$ which admits a pointwise $SO(1,1)$-symmetry. Let
$b_4^2-\ac12^2 \neq H$ for some $p\in \M$. Then $\M$ is affine
equivalent to
  $$\hs:I\times N^2\to \Rf:(v,x,y)\mapsto (\gamma_1(v),
  \gamma_2(v) \phi(x,y)),$$ 
  where
  $\phi: N^2 \to \mathbb R^3$ is a one-sheeted hyperboloid
  and $\gamma:I\to \mathbb R^2$ is a curve.
\end{thm}

Since the immersion $\hs$ has the same form as in the case of
$SO(2)$-symmetry, most of the computations for the converse theorem
are the same (cp. Theorem~\ref{thm:ClassC1} and the computations
thereafter). We just have to keep track of the different coordinates
and that now the affine metric $g$ of $\phi$ is indefinite.  If $H\neq
0$ then $\hs$ is indefinite and $V$ is spacelike iff
$H\sign(\ga_1\ga_2' - \ga_1' \ga_2) =\sign(\ga_1'\ga_2'' - \ga_1''
\ga_2')$. If $H=0$, then $\hs$ is indefinite and $V$ is spacelike iff
$\sign(\ga_2')= \sign(\ga_1'\ga_2'' - \ga_1'' \ga_2')$. Also the
computation of $K$ is completely the same as before in the proof of
Theorem~\ref{thm:ExC1}, implying that $K_{\hs_v}$ restricted to the
space spanned by $\hs_x$ and $\hs_y$ is a multiple of the
identity. Taking $V$ in direction of $\hs_v$, we can construct a LVB
$\{E,V,F\}$ with $E, F$ spanning $\Span\{\hs_x, \hs_y\}$ such that
$b_4\neq 0$ and $b_1= b_2 =b_5 = b_6= 0$. By the considerations in
Sec.~\ref{sec:max} we see that $\hs$ admits a pointwise
$SO(1,1)$-symmetry.
\begin{thm}\label{thm:ExC1t8} 
  Let $\phi:N^2 \to \mathbb R^3$ be a one-sheeted hyperboloid and
  let $\gamma: I \to \mathbb R^2$ be a curve, such that
  $\hs(t,v,w)=(\gamma_1(t), \gamma_2(t) \phi(v,w ))$ defines
  a 3-dimensional indefinite affine hypersphere. Then
  $\hs(N^2\times I)$ admits a pointwise $SO(1,1)$-symmetry.
  \begin{romanlist}
  \item If $\ga=(\ga_1,\ga_2)$ satisfies $\ga_2^2|\ga_1\ga_2' - \ga_1'
  \ga_2|^5= |\ga_1'\ga_2'' - \ga_1'' \ga_2'|
  (\ga_1')^2\neq 0$, then $\hs$ defines a 3-dimensional indefinite
  proper affine hypersphere.
  \item If $\ga=(\ga_1,\ga_2)$ satisfies $\ga_2^2|\ga_2'|^5=
  |\ga_1'\ga_2'' - \ga_1'' \ga_2'|
  (\ga_1')^2\neq 0$, then $\hs$ defines a 3-dimensional indefinite
  improper affine hypersphere.
  \end{romanlist}
\end{thm}

\subsubsection{The second case: $\nu \equiv 0$ and $H\neq 0$ on $\M$}
\label{sec:case2t8}

Next, we consider the case that $H =b_4^2- \ac12^2$ and $H\neq 0$ on
$\M$. It follows that $b_4\neq \pm \ac12$ on $\M$.

We already have seen that $\M$ admits a warped product structure. For a fixed
point $v_0$, we get from \eqref{D22t8} - \eqref{D33t8}, \eqref{eq32t8} and
\eqref{eq33t8}, using the notation $\tilde{\xi}=(b_4-\ac12) V + \xi$:
\begin{align}
D_E E =&+ \ac11 E, \label{D22c82} \\
D_E F =&  -\ac11 F +\tilde{\xi},\label{D23c82}\\
D_F E =&+\ac31 E+ \tilde{\xi},\label{D32c82}\\
D_F F =&  -\ac31 F,\label{D33c82} \\
D_{E_a} \tilde{\xi} =&0,\quad E_a\in\{E,F\}.
\end{align}

Thus, if $x$ and $y$ are local coordinates which span the second
distribution $L_2$, then we can interprete $\hs(v_0,x,y)$ as a
indefinite improper affine quadric, i.~e. a hyperbolic paraboloid, in a 
$3$-dimensional linear subspace.

After applying a translation and a change of coordinates, we may
assume that
\begin{equation}\label{phiinit}
\hs(v_0,x,y)=(x,y,f(x,y),0),
\end{equation}
with affine normal $\tilde{\xi}(v_0,x,y)=(0,0,1,0)$. Notice that by
\eqref{eq31t8} - \eqref{eq33t8} we have that
$D_{X}(\beta\tilde{\xi})=0$, $X\in {\cal X}(M)$. Taking suitable
initial conditions for the function $\beta$ ($\beta(v_0)=1$), we get
that $\beta\tilde{\xi}=(0,0,1,0)$ and thus
\begin{equation}\label{xibeta}
\xi= \frac1{\beta} (0,0,1,0)- (b_4-\ac12) V.
\end{equation}
Applied to \eqref{D11t8} we get for fixed $x,y$ the following second
order linear ODE:
\begin{equation*} \label{DVV}
D_V V = (\ac12- 3 b_4) V + \frac1{\beta} (0,0,1,0).
\end{equation*}
The solution has the form
\begin{equation} \label{Vintt8}
V(v,x,y)= \delta_1(v)(\delta_2(v) (0,0,1,0)+ C(x,y)), 
\end{equation}
where $\delta_1'(v)= (\ac12- 3 b_4) \delta_1(v)$ and $\delta_2'(v)=
\frac1{\beta \delta_1}(v)$.  To obtain $V$ at $v_0$, we consider
\eqref{D21t8} and \eqref{D31t8} and get that
\begin{equation*}
D_{E_a}(V-(\ac12+b_4) \hs) = 0,\qquad E_a \in \{E,F\}.
\end{equation*}
Evaluating at $v=v_0$, this means that there exists a constant vector
$\tilde{C}$ such that 
\begin{equation}\label{Vv0}
V(v_0,x,y)=(\ac12+b_4)(v_0) \hs(v_0,x,y) +\tilde{C}. 
\end{equation}
Since $\tilde{C}$ must be transversal to $\Span\{E,F,\hs\} = 
\Span\{E,F,\xi\}$, we can write:
\begin{equation*}
  V(v_0,x,y)=\alpha_1 (x,y,f(x,y),\alpha_2),
\end{equation*}
where $\alpha_1, \alpha2\neq 0$ and we applied an equiaffine
transformation so that $\tilde{C}=(0,0,0,\alpha_1\alpha_2)$. If we
compare with \eqref{Vintt8}, then we see that $C(x,y)=
\frac{\alpha_1}{\delta_1 (v_0)} (x,y,f(x,y),\alpha_2) -\delta_2(v_0)
(0,0,1,0)$, and therefore:
\begin{equation} \label{Vfinal}
V(v,x,y)= \delta_1(v)(x,y,f(x,y) + \delta_2(v), \alpha_2), 
\end{equation}
where $\delta_1'(v)= (\ac12- 3 b_4) \delta_1(v)$, $\delta_1(v_0)=
\alpha_1$, and $\delta_2'(v)= \frac1{\beta \delta_1}(v)$,
$\delta_2(v_0)=0$.  Respecting the initial condition \eqref{phiinit},
after integration we obtain that
\begin{equation*} 
\hs(v,x,y)= (\gamma_1(v) x,\gamma_1(v) y,\gamma_1(v) f(x,y) +
\gamma_2(v), \alpha_2 (\gamma_1(v)-1)),
\end{equation*}
where $\gamma_1'(v)= \delta_1(v)$, $\gamma_1(v_0)= 1$, and
$\gamma_2'(v)= \delta_1(v) \delta_2(v)$, $\gamma_2(v_0)=0$. After
applying an affine transformation we have shown:

\begin{thm} \label{thm:ClassC2t8} Let $\M$ be an indefinite proper affine 
  hypersphere of $\Rf$ which admits a pointwise 
 $SO(1,1)$-symmetry. Let 
  $H= b_4^2-\ac12^2 (\neq 0)$ on $\M$. Then $\M$ is affine equivalent with
  $$\hs:I\times N^2\to \Rf:(v,x,y)\mapsto (\gamma_1(v) x,\gamma_1(v)
  y,\gamma_1(v) f(x,y) + \gamma_2(v), \gamma_1(v)),$$ where $\psi: N^2
  \to \mathbb R^3:(x,y) \mapsto (x,y,f(x,y))$ is a hyperbolic
  paraboloid with affine normal $(0,0,1)$ and $\gamma:I\to \mathbb
  R^2$ is a curve.
\end{thm}
Since the immersion $\hs$ has the same form as in the case of
$SO(2)$-symmetry, most of the computations for the converse theorem
are the same (cp. Theorem~\ref{thm:ClassC2} and the computations
thereafter). We just have to keep track of the different coordinates
and that $\psi$ is indefinite, i.~e. $f_{xx} f_{yy}-f_{xy}^2=-1$. Thus
$\hs$ is indefinite and $V$ is spacelike iff $H\sign(\ga_1\ga_2' -
\ga_1' \ga_2) =\sign(\ga_1'\ga_2'' - \ga_1'' \ga_2')$. Also the
computation of $K$ is completely the same as before in the proof of
Theorem~\ref{thm:ExC2}, implying that $K_{\hs_v}$ restricted to the
space spanned by $\hs_x$ and $\hs_y$ is a multiple of the
identity. Taking $V$ in direction of $\hs_v$, we can construct a LVB
$\{E,V,F\}$ with $E, F$ spanning $\Span\{\hs_x, \hs_y\}$ such that
$b_4\neq 0$ and $b_1= b_2 =b_5 = b_6= 0$. By the considerations in
Sec.~\ref{sec:max} we see that $\hs$ admits a pointwise
$SO(1,1)$-symmetry.

\begin{thm}\label{thm:ExC2t8} 
  Let $\psi:N^2 \to \mathbb R^3$ be a hyperbolic paraboloid with
  affine normal $(0,0,1)$, and let $\gamma: I \to \mathbb R^2$ be a
  curve, such that $\hs(v,x,y)= (\gamma_1(v) x,\gamma_1(v)
  y,\gamma_1(v) f(x,y) + \gamma_2(v), \gamma_1(v))$ defines a
  3-dimensional indefinite proper affine hypersphere. Then
  $\hs(N^2\times I)$ admits a pointwise $SO(1,1)$-symmetry.\\
  Moreover, if $\ga=(\ga_1,\ga_2)$ satisfies $\ga_1^2|\ga_1\ga_2' -
  \ga_1' \ga_2|^5= |\ga_1'\ga_2'' - \ga_1'' \ga_2'| (\ga_1')^2\neq 0$,
  then $\hs$ defines a 3-dimensional indefinite proper affine
  hypersphere.
\end{thm}

\subsubsection{The third case: $\nu \equiv 0$ and $H=0$ on $\M$}
\label{sec:case3t8} 

 The final cases now are that $\nu \equiv 0$ and $H=0$ on the
whole of $\M$ and hence $\ac12=\pm b_4$. 

First we consider the case that $\ac12 = b_4 =:b>0$. Again we use that
$M^3$ admits a warped product structure and we fix a parameter
$v_0$. At the point $v_0$, we have by \eqref{D22t8}--\eqref{Dxit8}:
\begin{alignat*}{3}
&D_E E =&+ \ac11 E, &&  \\
&D_E F =&  &-\ac11 F&  +\xi, \\
&D_F E =&+\ac31 E&& +\xi, \\
&D_F F =&  & -\ac31 F,&  \\
&D_{X} \xi=& 0.
\end{alignat*}

Thus, if $x$ and $y$ are local coordinates which span the second
distribution $L_2$, then we can interprete $\hs(v_0,x,y)$ as an
indefinite improper affine quadric, i.~e. a hyperbolic paraboloid, 
in a $3$-dimensional linear subspace.

Completely analogous to the case of $SO(2)$-symmetry we obtain by
\eqref{Db4t8} that we can take $b=\frac1{4v}$, $v>0$, and arrive at
the following linear second order ordinary differential equation
(cp. \eqref{DTT}):
$$\pp{v}\hs + \frac1{2v} \pv \hs = -\xi,$$
where $\xi(v,x,y)=(0,0,1,0)$. Thus we end up with (cp. \eqref{result3.1}):
\begin{equation}\label{result3.1t8}
\hs(v,x,y)= (vx, vy, v f(x,y) - c v^4, v), \quad c,v>0.
\end{equation}

Next we consider the case that $-\ac12 =b_4= :b>0$. Again we use that
$M^3$ admits a warped product structure and we fix a parameter
$v_0$. As in the second case we can interpret $\hs(v_0,x,y)$ as an
indefinite improper affine quadric
(cp. \eqref{D22c82}--\eqref{D33c82}) with constant affine normal
$\tilde{\xi} =2b V +\xi$, i.~e. we get a hyperbolic paraboloid in a
$3$-dimensional linear subspace.  Also, in analogy to the case of
$SO(2)$-symmetry (just notice the different sign in the definition of
$\tilde{\xi}$), we can take $b=\frac1{4|v|}$, $v<0$, and arrive at the
following ordinary differential equation:
\begin{equation*}
\pv \hs= \sqrt{|v|} C +2v \xi, \quad v<0.
\end{equation*}
As before we get that the hypersurface is affine equivalent to:
\begin{equation}\label{result3.2t8}
\hs(v,x,y)= (x,y, f(x,y) + c v^3, v^4 ),\quad c,v>0.
\end{equation}
Combining both results \eqref{result3.1t8} and \eqref{result3.2t8} we have:

\begin{thm} \label{thm:ClassC3t8} Let $\M$ be an indefinite improper affine 
  hypersphere of $\Rf$ which admits a pointwise $SO(1,1)$-symmetry. 
  Let $\ac12^2 =b_4^2$ on $M^3$. Then $\M$ is affine
  equivalent with either
  \begin{align*} &\hs:I\times N^2\to \Rf:(v,x,y)\mapsto (vx, vy, v f(x,y) 
  - c v^4, v),\quad(\ac12 =b_4)\quad \text{or}\\
  &\hs:I\times N^2\to \Rf:(v,x,y)\mapsto (x,y, f(x,y) + c v^3, v^4 )
  ,\quad(-\ac12 =b_4)
\end{align*}
 where $\psi: N^2 \to \mathbb R^3:(v,w) \mapsto (v,w,f(v,w))$ is a
  hyperbolic paraboloid with affine normal
  $(0,0,1)$ and $c,t\in \R^+$.\hfill \newline 
\end{thm}

The computations for the converse statement can be done completely
analogous to the previous cases, they even are simpler (the curve is
given parametrized). 
\begin{thm}\label{thm:ExC3t8} 
  Let $\psi: N^2 \to \mathbb R^3:(v,w) \mapsto (v,w,f(v,w))$ be a
  hyperbolic paraboloid with affine normal
  $(0,0,1)$. Define $\hs(t,v,w)= (t v, t w,t f(v,w) -c t^4,t)$
  or $\hs(t,v,w)= (v, w, f(v,w) +c t^3,t^4)$, where $t\in \R^+$, $c\neq 0$.
  Then $\hs$ defines a 3-dimensional indefinite improper affine
  hypersphere, which admits a pointwise $SO(1,1)$-symmetry.
\end{thm}

\bibliography{/homes/geometer/schar/tex/bibs/csbib}

\vspace{2ex}
\noindent Authors address: C. Scharlach, Technische Universit\"at Berlin,
Fakult\"at II, Inst. f. Mathematik, MA 8-3, 10623 Berlin, Germany\\
E-mail: \texttt{schar@math.tu-berlin.de}


\end{document}